\documentclass{article}

\usepackage[latin1]{inputenc}    
\usepackage[OT1]{fontenc}
\usepackage[english]{babel}

\usepackage{fullpage}

\usepackage{amsmath}
\usepackage{amssymb}
\usepackage{amsthm}
\usepackage{mathrsfs}
\usepackage{graphicx}
\usepackage{dsfont}
\usepackage{wrapfig}
\usepackage{comment}
\usepackage{rotating}
\usepackage{hyperref} 
\usepackage[all]{xy}

\newcommand{\R}{\mathbb{R}}
\newcommand{\Rop}{\mathbb{R}^{\text{op}}}
\newcommand{\Rext}{\mathbb{R}_{\text{Ext}}}
\newcommand{\N}{\mathbb{N}}
\newcommand{\Z}{\mathbb{Z}}
\newcommand{\X}{\mathbb{X}}

\newcommand{\T}{\mathbb{T}}
\newcommand{\U}{\mathcal{U}}
\newcommand{\V}{\mathcal{V}}

\newcommand{\M}{\mathcal{M}}
\newcommand{\I}{\mathcal{I}}
\newcommand{\J}{\mathcal{J}}
\newcommand{\Reeb}{\mathrm{R}}
\newcommand{\Ord}{\mathrm{Ord}}
\newcommand{\Ext}{\mathrm{Ext}}
\newcommand{\Rel}{\mathrm{Rel}}
\newcommand{\Dg}{\mathrm{Dg}}
\newcommand{\LZZ}{\mathrm{LZZ}}
\newcommand{\DgZZ}{\mathrm{LBc}}
\newcommand{\CZZ}{\mathrm{CZZ}}
\newcommand{\DgCZZ}{\mathrm{CBc}}
\newcommand{\Mapper}{\mathrm{M}}
\newcommand{\MMapper}{\overline{\mathrm{M}}}
\newcommand{\im}{\mathrm{im}}
\newcommand{\Shift}{\mathrm{Shift}}
\newcommand{\Merge}{\mathrm{Merge}}
\newcommand{\Split}{\mathrm{Split}}
\newcommand{\Crit}{\mathrm{Crit}}
\newcommand{\Rips}{\mathrm{Rips}}
\newcommand{\distb}{d_{\rm b}}
\newcommand{\distfd}{d_{\rm FD}}
\newcommand{\card}{{\rm card}}
\newcommand{\diststair}{d_{\rm{b},\Theta}}
\newcommand{\Diag}{\Delta}
\newcommand{\distcov}{d_{\mathcal I}}
\newcommand{\End}{\mathrm{End}}
\newcommand{\Xset}{X}

\newcommand{\funcPL}{f^{\rm PL}}
\newcommand{\freeb}[1]{\tilde{#1}}
\newcommand{\frips}[1]{{#1}^{\rm PL}}
\newcommand{\Xs}{X}

\newcommand{\e}{\varepsilon}
\newcommand{\disth}{d_{\rm H}}
\newcommand{\symdiff}{{\scriptstyle\triangle}}
\newcommand{\Mapperv}{\Mapper^\bullet}
\newcommand{\MMapperv}{\MMapper^\bullet}
\newcommand{\Mappere}{\Mapper^\triangle}
\newcommand{\MMappere}{\MMapper^\triangle}
\newcommand{\MMtel}{\bar T_{\mathcal I}}
\newcommand{\Mtel}{T_{\mathcal I}}
\newcommand{\fMMtel}{\bar f_{\mathcal I}}
\newcommand{\fRMMtel}{\tilde{\bar f}_{\mathcal I}}
\newcommand{\fMtel}{f_{\mathcal I}}
\newcommand{\fRMtel}{\tilde{f}_{\mathcal I}}
\newcommand{\fPLMMtel}{\bar{f}^{\rm PL}_{\mathcal I}}
\newcommand{\fPLMtel}{f^{\rm PL}_{\mathcal I}}
\newcommand{\fRPLMMtel}{\bar{\sf{m}}^{\rm PL}_{\mathcal I}}
\newcommand{\fRPLMtel}{{\sf{m}^{\rm PL}_{\mathcal I}}}
\newcommand{\fM}{{\sf m}_{\mathcal I}}
\newcommand{\fMM}{\bar {\sf{m}}_{\mathcal I}}
\newcommand{\fMR}{\tilde{\sf{m}}_{\mathcal I}}
\newcommand{\fMMR}{\tilde{\bar{\sf{m}}}_{\mathcal I}}
\newcommand{\StairOrd}{Q_O^{\mathcal I}}
\newcommand{\StairRel}{Q_R^{\mathcal I}}
\newcommand{\StairExtM}{Q_E^{\mathcal I}}
\newcommand{\StairExtMM}{Q_{E^-}^{\mathcal I}}

\newtheorem{theorem}{Theorem}[section]
\newtheorem{corollary}[theorem]{Corollary}
\newtheorem{proposition}[theorem]{Proposition}
\newtheorem{lemma}[theorem]{Lemma}
\newtheorem{definition}[theorem]{Definition}
\newtheorem{rmq}[theorem]{Remark}

\title{Structure and Stability of the 1-Dimensional Mapper}
\author{Mathieu Carri\`ere, Steve Oudot}
\date{}

\begin{document}

\maketitle

\begin{abstract}
Given a continuous function~$f:X\to\R$ and a cover~$\I$ of its image
by intervals, the Mapper is the nerve of a refinement of the pullback
cover $f^{-1}(\I)$.  Despite its success in applications, little is
known about the structure and stability of this construction from a
theoretical point of view. As a pixelized version of the Reeb graph
of~$f$, it is expected to capture a subset of its features (branches,
holes), depending on how the interval cover is positioned with respect
to the critical values of the function. Its stability should also
depend on this positioning. We propose a theoretical framework that
relates the structure of the Mapper to the one of the Reeb graph,
making it possible to predict which features will be present and which
will be absent in the Mapper given the function and the cover, and for
each feature, to quantify its degree of (in-)stability.  Using this
framework, we can derive guarantees on the structure of the Mapper, on
its stability, and on its convergence to the Reeb graph as the
granularity of the cover~$\I$ goes to zero.
\end{abstract}

\section{Introduction}
\label{sec:Intro}

Many data sets nowadays come in the form of point clouds with function
values attached to the points. Such data may come either from direct
measurements (e.g. think of a sensor field measuring some physical
quantity like temperature or humidity), or as a byproduct of some data
analysis pipeline (e.g. think of a word function in the quantization
phase of the bag-of-words model). There is a need for summarizing such
data and for uncovering their inherent structure, to enhance further
processing steps and to ease interpretation.

One way of characterizing the structure of a scalar field $f:X\to\R$
is to look at the evolution of the topology of its {\em level
  sets}---i.e. sets of the form $f^{-1}(\{\alpha\})$, for $\alpha$
ranging over~$\R$. This information is summarized in a mathematical
object called the {\em Reeb graph} of the pair $(X,f)$, denoted by
$\Reeb_f(X)$ and defined as the quotient space obtained by identifying
the points of $X$ that lie in the same connected component of the same
level set of~$f$~\cite{Reeb46}. The Reeb graph is known to be
a graph (technically, a multi-graph) when $X$ is a smooth
manifold and $f$ is a Morse function, or more generally when $f$ is of
{\em Morse type} (see Definition~\ref{def:Morse-type}). Moreover,
since the map~$f$ is constant over equivalence classes, there is a
well-defined induced map on $\Reeb_f(X)$.

The connection between the topology of the Reeb graph and the one of
its originating pair $(X, f)$ has been the object of much study in the
past and is now well understood. It has gained increasing interest in
the recent years, with the introduction of persistent
homology~\cite{Edelsbrunner02b} and its extended
version~\cite{Cohen09}. Indeed, the {\em extended
  persistence diagram} of the induced map~$\tilde f$ (see
Section~\ref{sec:Def} for a formal definition) describes the structure
of the Reeb graph in the following sense. It takes the form of a
multiset of points in the plane, where each point is matched with a
{\em feature} (branch or hole) of the Reeb graph in a one-to-one
manner. Furthermore, the coordinates of the point characterize the
{\em span} of the feature, that is, the interval of $\R$ spanned by
its image under~$\tilde f$. The vertical distance of the point to
the diagonal $\Delta=\{(x,x) : x\in\R\}$ measures the length of
that interval and thereby quantifies the prominence of the
feature. Thus, the extended persistence diagram plays the role of a
``bag-of-features'' type signature, summarizing the Reeb graph
through its list of features together with their spans, and forgetting
about the actual layout of those features.

One issue with the Reeb graph is its computation. Indeed, when the
pair $(X,f)$ is known only through a finite set of measurements, the
graph can only be approximated within a certain error. Quantifying
this error, in particular finding the right metric in which to measure
it, has been the object of intense investigation in the recent
years. A simple choice is to use the {\em bottleneck distance} between
the extended persistence diagram associated to the Reeb graph and the one
associated to its
approximation~\cite{Cohen07}. This pseudometric
treats the Reeb graph and its approximation as bags of features, and
it measures the differences between their respective sets of features. In
particular, it is oblivious to the layouts of the features in each of
the graphs. Other distances have been proposed recently, to capture a
greater part of this layout~\cite{Bauer14,Bauer15,deSilva16}.

Building approximations from finite point samples with scalar values
is a problem in its own right. A natural approach is to build a
simplicial complex (for instance the {\em Rips complex}) on top of
the point samples, to serve as a proxy for the underlying continuous
space; then, to extend the scalar values at the vertices to a
piecewise-linear (PL) function over the simplicial complex by linear
interpolation; finally, to apply some exact computation algorithm for
PL functions. This is the approach advocated by Dey and
Wang~\cite{Dey13a}, who rely on the $O(n\log n)$ expected time
algorithm of Harvey, Wenger and Wang~\cite{Harvey10} for
the last step. The drawbacks of this approach are:
\begin{itemize}
\item Its relative complexity: the Reeb graph computation from the PL
  function is based on collapses of its simplicial domain that may
  break the complex structure temporarily and therefore require some repairs.
\item Its overall computational cost: here, $n$ is not the number of
  data points, but the number of vertices, edges and triangles of the
  Rips complex, which, in principle, can be up to cubic in the number
  of data points. Indeed, triangles are needed to compute an
  approximation of the Reeb graph, in the same way as they are to
  compute 1-dimensional homology.
\end{itemize}

The {\em Mapper}\footnote{In this article we call {\em Mapper} the
  mathematical object, not the algorithm used to build it. Moreover,
  we focus on the case where the codomain of the function is~$\R$.} was
introduced by Singh, M\'emoli and Carlsson~\cite{Singh07} as a new
mathematical object to summarize the topological structure of a pair
$(X,f:X\to\R^d)$. Its construction depends on the choice of a cover~$\I$ of the
image of $f$ by open sets. Pulling back~$\I$ through $f$ gives an
open cover of the domain $X$. This cover may have some elements that
are disconnected, so it is refined into a connected cover by splitting
each element into its various connected components. Then, the Mapper
is defined as the nerve of the connected cover, having one vertex per
element, one edge per pair of intersecting elements, and more
generally, one $k$-simplex per non-empty ($k+1$)-fold
intersection. From a philosophical point of view, the Mapper can be
thought of as a {\em pixelized version} of the Reeb space, where the
resolution is prescribed by the cover~$\I$. From a practical point of
view, its construction from point cloud data is very easy to describe
and to implement, using standard graph traversals to detect connected
components. Furthermore, it only requires to build the 1-skeleton
graph of the Rips complex, whose size scales up at worst quadratically
(and not cubically) with the size of the input point cloud.

As a simple alternative to the Reeb space, the Mapper has been the
object of much interest by practitioners in the data sciences. It has
played a key role in several success stories, such as the
identification of a new subgroup of breast
cancers~\cite{Nicolau11}, or the elaboration of a new
classification of player positions in the NBA~\cite{NBA}, due to its ability to
deal with very general functions and datasets. Meanwhile,
it has become the flagship component in the software suite developed
by Ayasdi, a data analytics company founded in the late 2000's whose
interest is to promote the use of topological methods in the data
sciences. Somewhat surprisingly, despite this success, very little is
known to date about the structure of the Mapper and its stability with
respect to perturbations of the pair $(X,f)$ or of the cover $\I$. Intuitively, 
when $f$ is scalar, as a pixelized version of the Reeb graph, the Mapper should capture
some of its features (branches, holes) and miss others, depending on
how the cover $\I$ is positioned with respect to the critical values
of~$f$. How can we formalize this phenomenon?  The stability of the
structure of the Mapper should also depend on this positioning. How
can we quantify it?  These are the questions addressed in this
article.

\paragraph{Contributions.} We draw an explicit connection between the Mapper and the Reeb graph, 
from which we derive guarantees on the structure  of the Mapper and quantities to measure its stability. Specifically:

\begin{itemize}
\item The connection happens through an intermediate object, called
  the {\em MultiNerve Mapper}, which we define as the {\em
    multinerve} of the connected pullback cover in the sense
  of~\cite{Verdiere12}. The Mapper and its MultiNerve variant are related
  through the usual Nerve-vs-MultiNerve connection
  (see Lemma~\ref{lem:nerve-vs-multinerve}).

\item Given a pair $(X,f)$ and an interval cover $\I$, we relate the topological structure of the
  (MultiNerve) Mapper to the one of the Reeb graph.
  More precisely, we characterize the topology of these objects with particular 
  {\em zigzag persistence modules} that we relate to each other with the so-called
  Mayer-Vietoris half-pyramid (Theorem~\ref{th:ExDgMultiNerve}). 
  This correspondence is oblivious to the
  actual layouts of the topological features in the two graphs, which in principle
  could differ.

\item The previous connection allows us to derive a signature for the
  (MultiNerve) Mapper, which takes the form of an extended persistence
  diagram.
  The points in this diagram are in one-to-one correspondence
  with the features (branches, holes) in the (MultiNerve)
  Mapper. Thus, like the extended persistence diagram of the induced
  map~$\tilde f$ for the Reeb graph, our diagram for the (MultiNerve)
  Mapper serves as a bag-of-features type signature of its structure.

\item An interesting property of our signature is to be
  predictable\footnote{As a byproduct, we
    also clarify the relationship between the extended persistence diagrams of~$\tilde f$ and~$f$ (Theorem~\ref{thm:pdreeb}).}  given the extended persistence
  diagram of the induced map~$\tilde f$. Indeed, it is obtained from
  this diagram by removing the points lying in certain {\em
    staircases} that are defined solely from the cover~$\I$ and that
  encode the mutual positioning of the intervals of the cover. Thus,
  the signature for the (MultiNerve) Mapper is a subset of the one
  for the Reeb graph, which provides theoretical evidence to the
  intuitive claim that the Mapper is a pixelized version of the Reeb
  graph. Then, one can easily derive sufficient conditions under which
  the bag-of-features structure of the Reeb graph is preserved in the
  (MultiNerve) Mapper, and when it is not, one can easily predict
  which features are preserved and which ones disappear
  (Corollary~\ref{cor:Reeb_approx}).

\item The staircases also play a role in the stability of the
  (MultiNerve) Mapper, since they prescribe which features will
  (dis-)appear as the function~$f$ is perturbed. Stability is then
  naturally measured by a slightly modified version of the bottleneck
  distance, in which the staircases play the role of the diagonal. Our
  stability guarantees (Theorem~\ref{thm:MStab}) follow  easily
  from the general stability theorem for extended
  persistence~\cite{Cohen09}. Similar guarantees hold when the domain $X$ or
  the cover $\I$ is perturbed (Theorems~\ref{thm:DStab}
  and~\ref{thm:CStab}).
\item These stability guarantees can be exploited in practice to
  approximate the signatures of the Mapper and MultiNerve Mapper from
  point cloud data. The approach boils down to applying known scalar
  field analysis techniques~\cite{Chazal09a} then pruning the
  obtained extended persistence diagrams using the staircases
  (Theorem~\ref{th:sig-approx-general}). The approach becomes more
  involved if one wants to further guarantee that the approximate signature
  does correspond to some perturbed Mapper or MultiNerve Mapper
  (Theorem~\ref{th:sig-approx}).
\item We also refine the analysis by showing that
  the MultiNerve Mapper itself is a Reeb graph, for a perturbed pair
  $(X', f')$ (Theorem~\ref{th:Dg}). Furthermore, we are able to track
  the changes that occur in the structure of the Reeb graph as we go
  from the initial pair $(X,f)$ to its perturbed version $(X',
  f')$. This allows us to compute the {\em functional distortion distance}
  between the (MultiNerve) Mapper and the Reeb graph (Theorem~\ref{th:dfd}). 
  Our main proof technique consists in progressively perturbing
  the so-called {\em telescope}~\cite{Carlsson09b} associated with the
  pair $(X,f)$.
  To be more specific, we decompose the
  perturbation into a sequence of elementary perturbations with a
  predictable effect on the functional distortion distance. We believe the introduction of
  these elementary perturbations and the analysis of their effects
  on persistence diagrams and on the functional distortion distance between telescopes
  are of an independent interest (see Section~\ref{sec:stabdfd}). In particular, these
  elementary perturbations have already been used in other works about Reeb graphs and Mappers~\cite{Carriere17c,Carriere17a}.
\end{itemize}

\paragraph{Related work.}

Reeb graphs can be seen as particular types of {\em skeletonization},
when one tries to recover the geometric structure of data with graphs.
Several kinds of graph-like geometric structures have been studied
within the past few years, such as {\em persistent skeletons}~\cite{Kurlin15} or {\em graph-induced simplicial complexes}~\cite{Dey13b}.
See~\cite{Kurlin15} for a list of references. 
As mentioned previously, Reeb graphs are now well understood and have been used in a wide range of applications.  
Algorithms for their computation have been proposed,
as well as metrics for their comparison.
We refer the interested reader to the survey~\cite{Biasotti08} and to the
introductions of~\cite{Bauer14} and~\cite{Bauer15} for a comprehensive
list of references.  In a recent study, even more structure has been
given to the Reeb graphs by categorifying them~\cite{deSilva16}.  

A lot of variants of these graphs have also been studied in the last
decade to face the common issues that come with the Reeb graphs
(complexity and computational cost among others).  The
Mapper~\cite{Singh07} is one of them.  Chazal et al.~\cite{Chazal14a}
introduced the {\em $\lambda$-Reeb graph}, which is another type of Reeb
graph pixelization with intervals. It is the quotient space
obtained by identifying the points with the transitive closure of the
following relation: $x\sim y\ \Leftrightarrow$ $x,y$ belong to the
same level set and $x,y$ belong to the same element of a given
family of intervals. The computation is easier than for the Reeb
graph, and the authors can derive upper bounds on the Gromov-Hausdorff
distance between the space and its Reeb or $\lambda$-Reeb graph.  However,
this is too much asking in general; as a result, the hypothesis made on the space
are very strong (it has to be close to a metric graph in the Gromov-Hausdorff distance already).
Moreover, taking the transitive closure makes the structure of the output more
difficult to interpret.

Joint Contour Nets~\cite{Carr14,Chattopadhyay16} and Extended Reeb
graphs~\cite{Barra14} are Mapper-like objects.  The former is the
Mapper computed with the cover of the codomain given by rounding the function values,
while the latter is the Mapper computed from a partition of the domain with
no overlap. Both structures are used for scientific
purposes (visualization, shape descriptors, to name a few) and
algorithms are proposed for their computation.
Babu~\cite{Babu13} characterized the Mapper with coarsened levelset zigzag persistence modules and showed that, as
the lengths of the intervals in the cover~$\I$ go to zero uniformly,
the Mapper of a real-valued function converges to the continuous
Reeb graph in the bottleneck distance.
Similarly,
Munch and Wang~\cite{Munch16} recently characterized the Mapper
with constructible cosheaves and showed the same type of convergence for
both the Joint Contour Net and the Mapper 
in the so-called {\em interleaving
  distance}~\cite{deSilva16}. Their result holds in the general case
of vector-valued functions. 
Differently, we restrict the focus to
real-valued functions but are able to make non-asymptotic 
claims (Corollary~\ref{cor:Reeb_approx}).

On another front, Stovner~\cite{Stovner12} proposed a categorified
version of the Mapper, which is seen as a covariant functor from the
category of covered topological spaces to the category of simplicial
complexes.  Dey et al.~\cite{Dey16} pointed out the inherent
instability of the Mapper and introduced a multiscale variant, called
the {\em MultiScale Mapper}, which is built by taking the Mapper over
a hierarchy of covers of the codomain. 
They derived a stable signature by considering the persistence diagram
of this family.  Unfortunately, their construction is hard to relate
to the original Mapper. Here we work with the original Mapper directly
and answer two open questions from~\cite{Dey16}, introducing a
signature that gives a complete description of the set of features of
the Mapper together with a quantification of their stability and a
provable way of approximating them from point cloud data.

\section{Background}
\label{sec:Def}

Throughout the paper we work with singular homology with coefficients in the
field~$\Z_2$, which we omit in our notations for simplicity. We use
the term "connected" as a shorthand for "path-connected".
%
Given a real-valued function $f$ on a topological space $X$, and an
interval $I\subseteq\R$, we denote by $X_f^I$ the preimage
$f^{-1}(I)$. We omit the subscript $f$ in the notation when there is
no ambiguity in the function considered.

\subsection{Morse-Type Functions}
\label{sec:Morse-type}

We restrict our focus to the class of real-valued functions called
{\em Morse-type}. These are generalizations of the classical Morse
functions: 

\begin{definition}\label{def:Morse-type}
A continuous real-valued function $f$ on a topological space $X$ 
is \emph{of Morse type} if:
\begin{itemize}
\item[{\rm (i)}] There is a finite set $\emph{Crit}(f)=\{a_1<...<a_n\}$, 
called the set of \emph{critical values},
such that over every open interval $(a_0=-\infty,a_1),...,(a_i,a_{i+1}),...,(a_n,a_{n+1}=+\infty)$
there is a compact and locally connected space $Y_i$ 
and a homeomorphism $\mu_i:Y_i\times(a_i,a_{i+1})\rightarrow X^{(a_i,a_{i+1})}$
such that $\forall i=0,...,n, f|_{X^{(a_i,a_{i+1})}}=\pi_2\circ\mu_i^{-1}$, where
$\pi_2$ is the projection onto the second factor;

\item[{\rm (ii)}]$\forall i=1,...,n-1,\mu_i$ extends to a continuous function $\bar{\mu}_i:Y_i\times[a_i,a_{i+1}]\rightarrow X^{[a_i,a_{i+1}]}$
-- similarly $\mu_0$ extends to $\bar{\mu}_0:Y_0\times(-\infty,a_1]\rightarrow X^{(-\infty,a_1]}$
and $\mu_n$ extends to $\bar{\mu}_n:Y_n\times[a_n,+\infty)\rightarrow X^{[a_n,+\infty)}$;

\item[{\rm (iii)}]Each levelset $X^t$ has a finitely-generated homology.
\end{itemize}
\end{definition}

Morse functions are known to be of Morse type while the converse is
clearly not true.  In fact, Morse-type functions do not have to be
differentiable, and their domain does not have to be a smooth manifold
nor even a manifold at all.  Furthermore, it is possible to find
Morse-type functions that are not Morse even though they satisfy the
previous assumptions (think of the Gaussian curvature on a torus, for
instance).

\subsection{Extended Persistence}
\label{sec:Persistence}

Let $f$ be a real-valued function on a topological space $X$.
The family $\{X^{(-\infty, \alpha]}\}_{\alpha\in\R}$ of
  sublevel sets of $f$ defines a {\em filtration}, that is, it is
  nested w.r.t. inclusion: $X^{(-\infty, \alpha]}\subseteq
X^{(-\infty, \beta]}$ for all $\alpha\leq\beta\in\R$. The
  family $\{X^{[\alpha, +\infty)}\}_{\alpha\in\R}$ of superlevel sets
  of $f$ is also
  nested but in the opposite direction: $X^{[\alpha, +\infty)}\supseteq X^{[\beta, +\infty)}$
  for all $\alpha\leq\beta\in\R$. We can turn it into a filtration by
  reversing the real line. Specifically, let $\Rop=\{\tilde{x}:x\in\R\}$,
  ordered by $\tilde{x}\leq\tilde{y}\Leftrightarrow x\geq y$. We index
  the family of superlevel sets by $\Rop$, so now we have a
  filtration: $\{X^{[\tilde\alpha, +\infty)}\}_{\tilde\alpha\in\Rop}$, with
  $X^{[\tilde\alpha, +\infty)}\subseteq X^{[\tilde\beta, +\infty)}$ for all
  $\tilde\alpha\leq\tilde\beta\in\Rop$.

\begin{figure}[h]\centering
\includegraphics[width=13cm]{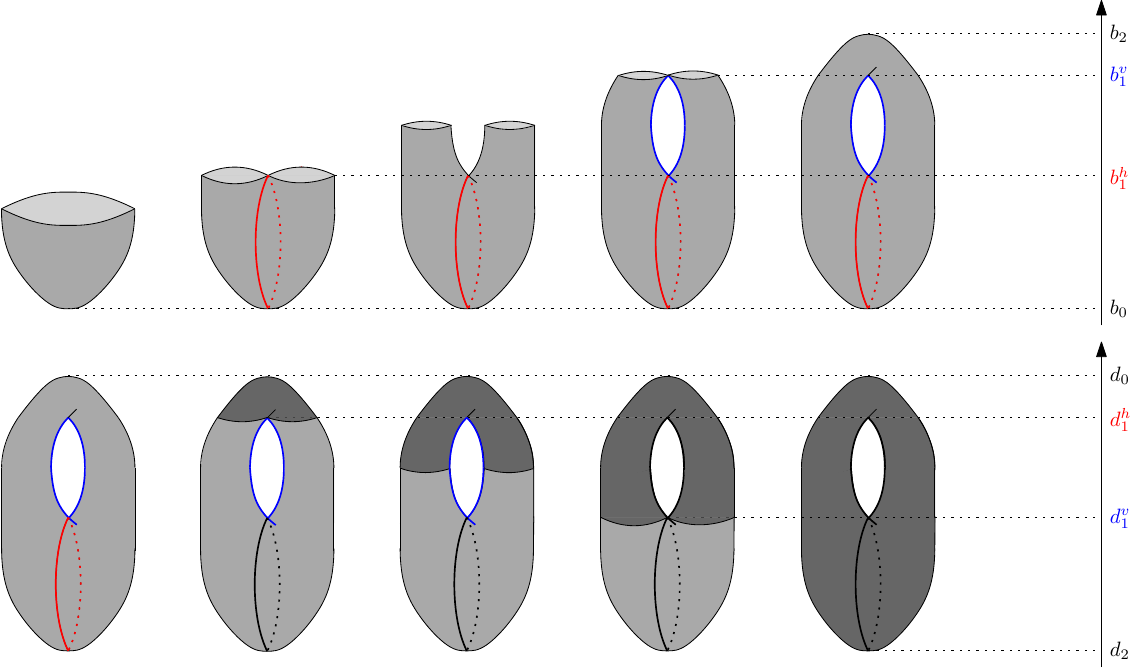}
\caption{The extended filtration of the
  height function on a torus.  The upper row displays the ordinary
  part of the filtration while the lower row displays the relative part.  The red and blue cycles both correspond to
  extended points in dimension 1.  The point corresponding to
  the red cycle is located above the diagonal ($d_1^h>b_1^h$), while
  the point corresponding to the blue cycle is located below the
  diagonal ($d_1^v>b_1^v$).}
\label{fig:ExFilt}
\end{figure}

Extended persistence connects the two filtrations at infinity as
follows. Replace each superlevel set $X^{[\tilde\alpha, +\infty)}$ by
  the pair of spaces $(X, X^{[\tilde\alpha, +\infty)})$ in the second
    filtration. This maintains the filtration property since we have
    $(X, X^{[\tilde\alpha, +\infty)})\subseteq (X, X^{[\tilde\beta,
          +\infty)})$ for all
        $\tilde\alpha\leq\tilde\beta\in\Rop$. Then, let
        $\Rext=\R\cup\{+\infty\}\cup\Rop$, where the order is
        completed by $\alpha<+\infty<\tilde{\beta}$ for all
        $\alpha\in\R$ and $\tilde\beta\in\Rop$. This poset is
        isomorphic to $(\R, \leq)$. Finally, define the {\em extended
          filtration} of $f$ over $\Rext$ by:
\[
\begin{array}{llll}
F_\alpha &=& X^{(-\infty, \alpha]} & \mbox{for $\alpha\in\R$}\\[0.5em]
F_{+\infty} &=& X \equiv (X,\emptyset)\\[0.5em]
F_{\tilde\alpha} &=& (X, X^{[\tilde\alpha, +\infty)}) & \mbox{for $\tilde\alpha\in\Rop$},
\end{array}
\]
where we have identified the space $X$ with the pair of spaces $(X,
\emptyset)$. This is a well-defined filtration since we have
$X^{(-\infty, \alpha]}\subseteq X \equiv (X, \emptyset) \subseteq (X,
X^{[\tilde\beta, +\infty)})$ for all $\alpha\in\R$ and
  $\tilde\beta\in\Rop$.  The subfamily $\{F_\alpha\}_{\alpha\in\R}$ is
  called the \textit{ordinary} part of the filtration, and the
  subfamily $\{F_{\tilde\alpha}\}_{\tilde\alpha\in\Rop}$ is called the
  \textit{relative} part. See Figure~\ref{fig:ExFilt} for an
  illustration.

Applying the homology functor $H_*$ to this filtration gives the so-called \textit{extended persistence module} ${\rm EP}(f)$:
%
\[
\begin{array}{llll}
{\rm EP}(f)_\alpha &=& H_*(F_\alpha)=H_*(X^{(-\infty, \alpha]}) & \mbox{for $\alpha\in\R$}\\[0.5em]
{\rm EP}(f)_{+\infty} &=& H_*(F_{+\infty})=H_*(X) \cong H_*(X,\emptyset)\\[0.5em]
{\rm EP}(f)_{\tilde\alpha} &=& H_*(F_{\tilde\alpha})=H_*(X, X^{[\tilde\alpha, +\infty)}) & \mbox{for $\tilde\alpha\in\Rop$},
\end{array}
\]
and where the linear maps between the spaces are induced by the
inclusions in the extended filtration.

For functions of Morse type, the extended persistence module can be
decomposed as a finite direct sum of closed-open {\em interval
  modules}---see e.g.~\cite{Chazal16}:
\[
{\rm EP}(f)\simeq\bigoplus_{k=1}^n \mathbb{I}[b_k, d_k),
\]
where each summand $\mathbb{I}[b_k, d_k)$ is made of copies of the
  field of coefficients at each index $\alpha\in [b_k, d_k)$, and of
    copies of the zero space elsewhere, the maps between copies of the
    field being identities.  Each summand represents the lifespan of a
    {\em homological feature} (connected component, hole, void, etc.) within the
    filtration. More precisely, the {\em birth time} $b_k$ and {\em
      death time} $d_k$ of the feature are given by the endpoints of
    the interval.  Then, a convenient way to represent the structure
    of the module is to plot each interval in the decomposition as a
    point in the extended plane, whose coordinates are given by the
    endpoints. Such a plot is called the \textit{extended persistence
      diagram} of $f$, denoted $\Dg(f)$.  The distinction
    between ordinary and relative parts of the filtration allows to
    classify the points in $\Dg(f)$ in the following way:
\begin{itemize}
\item points whose coordinates both belong to $\R$  are called \textit{ordinary} points; they correspond to homological features being born and then dying in the ordinary part of the filtration;
\item points whose coordinates both belong to $\Rop$ are called \textit{relative} points; they correspond to homological features being born and then dying in the relative part of the filtration;
\item points whose abscissa belongs to $\R$  and whose ordinate belongs to $\Rop$ are called \textit{extended} points; they correspond to homological features being born in the ordinary part and then dying in the relative part of the filtration. 
\end{itemize}
Note that ordinary points lie strictly above the diagonal
$\Delta=\{(x,x):x\in\R\}$ and relative points lie strictly below $\Delta$,
while extended points can be located anywhere, including on $\Delta$, e.g. connected components 
that lie inside a single critical level---see Section~\ref{sec:ReebGraphDef}.  It is common to
decompose $\Dg(f)$ according to this classification:
\[
\Dg(f)=\Ord(f)\sqcup\Rel(f)\sqcup\Ext^+(f)\sqcup\Ext^-(f),
\]
where by
convention $\Ext^+(f)$ includes the extended points located on the
diagonal~$\Delta$.

\paragraph{Persistence measure.}

From an extended persistence module ${\rm EP}(f)$ we derive a measure on the set of rectangles in the plane, called the \textit{persistence measure} and denoted $\mu_{\rm EP}$. Given a rectangle~$R=[a,b]\times[c,d]$ with $a<b\leq c<d$, we let 
\begin{equation}\label{eq:pers_meas}
\mu_{\rm EP}(R)=r_b^c-r_b^d+r_a^d-r_a^c,
\end{equation}
where $r_x^y$ denotes the rank of the linear map between the vector
spaces indexed by $x,y\in\R_{\Ext}$ in~${\rm EP}(f)$.
When ${\rm EP}(f)$ has a well-defined persistence diagram, $\mu_{\rm EP}(R)$ equals the total multiplicity of the diagram within the rectangle~$R$~\cite{Chazal16}.

\paragraph{Stability.}
An important property of extended persistence diagrams is to be stable
in the so-called {\em bottleneck distance} $d^\infty_{\rm b}$. The definition
of this distance is based on partial matchings between the
diagrams. Given two persistence diagrams $D,D'$, a \emph{partial
  matching} between $D$ and $D'$ is a subset $\Gamma$ of $D\times D'$
such that:
\[
\forall p\in D, \text{ there is at most one }p'\in D'\text{ s.t. }(p,p')\in\Gamma,
\]
\[
\forall p'\in D', \text{ there is at most one }p\in D\text{
  s.t. }(p,p')\in\Gamma.
\]
Furthermore, $\Gamma$ must match points of the same type (ordinary,
relative, extended) and of the same homological dimension only.  The
\emph{cost} of $\Gamma$ is:
\[
\text{cost}(\Gamma)=\max\left\{\max_{p\in D}\ \delta_D(p),\ \max_{p'\in D'}\ \delta_{D'}(p')\right\},
\]
where
\[
\delta_D(p)=\|p-p'\|_\infty \text{ if } \exists p'\in D'\text{ s.t. }(p,p')\in\Gamma \text{ and } d_\infty(p,\Delta)=\inf_{q\in\Delta} \|p-q\|_\infty \text{ otherwise},
\]
\[
\delta_{D'}(p')=\|p-p'\|_\infty \text{ if } \exists p\in D\text{ s.t. }(p,p')\in\Gamma \text{ and } d_\infty(p',\Delta)=\inf_{q\in\Delta} \|p'-q\|_\infty \text{ otherwise}.
\]
\begin{definition}
\label{def:bottleneck}
Let $D,D'$ be two persistence diagrams.
The \emph{bottleneck distance} between $D$ and $D'$ is:
\[
\distb(D,D')=\inf_{\Gamma}\ \emph{cost}(\Gamma),
\]
where $\Gamma$ ranges over all partial matchings between $D$ and $D'$.
\end{definition}
Note that $\distb$ is only a pseudometric, not a true metric, because points lying on $\Delta$ can be left unmatched at no cost.
\begin{theorem}[Stability Theorem in~\cite{Cohen09}]
\label{th:Stab}
For any  Morse-type functions $f,g:X\rightarrow\R$,
\[
d^\infty_{\emph{b}}(\Dg(f),\Dg(g))\leq \|f-g\|_\infty.
\]
\end{theorem}
Moreover, as pointed out in~\cite{Cohen09}, the theorem can be
strengthened to apply to each subdiagram
$\Ord$, $\Ext^+$, $\Ext^-$, $\Rel$ and to each homological dimension
individually.

\subsection{Zigzag persistence}
\label{sec:ZZ}

Let $f:X\rightarrow\R$ be a Morse-type function, and let $\Crit(f)=\{a_1,\cdots,a_n\}$ be its set of critical values.
Let $-\infty=a_0<s_0<a_1<s_1<a_2<\cdots<s_{n-1}<a_n<s_n<a_{n+1}=+\infty$. Then, for any $1\leq i \leq j \leq n$, we define $X_i^j=X^{[s_i,s_j]}$,
and the {\em levelset zigzag} as the following sequence of $2n+1$ nodes:
$$X_0^0\hookrightarrow X_0^1 \hookleftarrow X_1^1 \hookrightarrow X_1^2 \hookleftarrow \cdots \hookrightarrow X_{n-1}^n \hookleftarrow X_n^n,$$
where each arrow is the canonical inclusion. Applying the homology functor $H_*$ to the levelset zigzag gives the so-called
{\em levelset zigzag persistence module} $\LZZ(f)$, where the linear maps between the spaces are induced by the inclusions.
For functions of Morse type, the levelset zigzag persistence module decomposes as a finite direct sum of closed interval modules:
$$\LZZ(f)\simeq \bigoplus_{k=1}^{m}\mathbb{I}[X_i^{i'},X_j^{j'}],$$
where $i'$ is either $i$ or $i+1$, and similarly for $j'$. Hence, the classification given by Table~\ref{tab:classif}.

\begin{table}[h]\centering
\caption{Classification of intervals in a levelset zigzag persistence barcode.}
\label{tab:classif}
\begin{tabular}{|c|c|c|c|c|}
\hline
Type & I                                             & II                                            & III                                             & IV \\
\hline
     & \( \begin{array}{l} i_\alpha'=i_\alpha+1 \\ j_\alpha'=j_\alpha \end{array}\)  
     & \( \begin{array}{l} i_\alpha'=i_\alpha \\ j_\alpha'=j_\alpha+1 \end{array}\)   
     & \( \begin{array}{l} i_\alpha'=i_\alpha+1 \\ j_\alpha'=j_\alpha+1\end{array}\)  
     & \( \begin{array}{l} i_\alpha'=i_\alpha \\ j_\alpha'=j_\alpha\end{array}\) \\
\hline
\end{tabular}
\end{table} 

Moreover, each summand $\mathbb{I}[X_i^{i'},X_j^{j'}]$ is made of copies of the
field of coefficients for each space between $X_i^{i'}$ and $X_j^{j'}$ and of
copies of the zero space elsewhere, the maps between copies of the field
being identities. 
The disjoint union of all of these intervals is called the {\em levelset zigzag persistence barcode} $\DgZZ(f)$.

\begin{figure}\centering
\includegraphics[width=12.5cm]{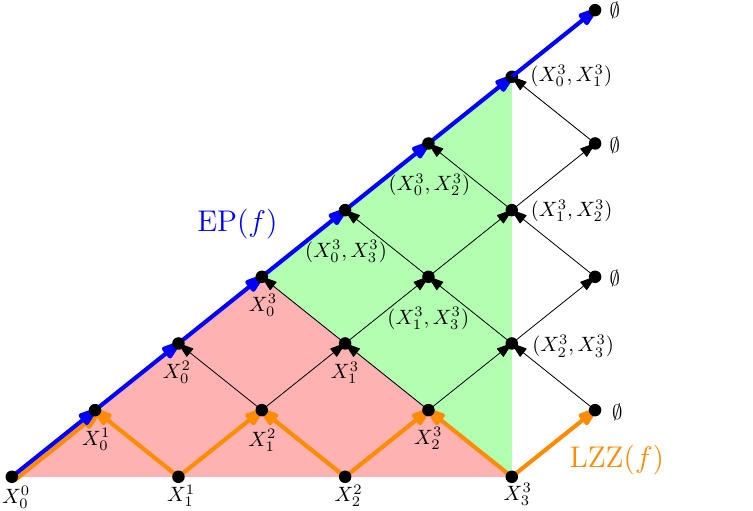}
\caption{Mayer-Vietoris half-pyramid when the Morse-type function has
three critical values. It is composed of two faces of the full Mayer-Vietoris pyramid: 
the south face (red) and the east face (green). The extended persistence module ${\rm EP}(f)$ is in blue and 
the levelset zigzag persistence module $\LZZ(f)$ is in orange.}
\label{fig:ZZpyramid}
\end{figure}

\paragraph{Mayer-Vietoris half-pyramid.} The Mayer-Vietoris half-pyramid is the diagram of topological spaces 
and inclusions displayed in Figure~\ref{fig:ZZpyramid}.
We refer the reader to~\cite{Carlsson09b} for more details.
Any zigzag within the Mayer-Vietoris half-pyramid that stretches from the left boundary (i.e. the node $X_0^0)$ 
to the right boundary without backtracking
is called {\em monotone}. Theorem~\ref{th:MY} 
below relates any two monotone zigzags. 

\begin{theorem}[Pyramid Theorem in~\cite{Carlsson09b}]\label{th:MY}
For any Morse-type function $f$, there exists a bijection between the barcodes of any pair of monotone zigzag persistence modules in the Mayer-Vietoris half-pyramid.
\end{theorem}

Since the extended persistence module of $f$ is a monotone zigzag persistence
module---more precisely the principal diagonal---of the Mayer-Vietoris half-pyramid, we have the following corollary.

\begin{corollary}[Table 1 in~\cite{Carlsson09b}]\label{cor:EPZZ}
For any Morse-type function $f$, there exists a bijection between $\Dg(f)$ and $\DgZZ(f)$,
which is described in Table~\ref{tab:EPZZ}.
\end{corollary}

\begin{table}[h]\centering
\caption{This table gives the correspondences between the points of $\Dg(f)$ and the intervals of $\DgZZ(f)$.
The minus sign on some intervals of $\DgZZ(f)$ means that the homological dimension of that interval is equal to the dimension of 
its corresponding point in $\Dg(f)$ minus 1. }
\label{tab:EPZZ}
\begin{tabular}{|l|l|l|l|l|}
\hline
Type                      & $\Ord$                         & $\Rel$                      & $\Ext^+$                          & $\Ext^-$                   \\
\hline
$\Dg(f)$                  & $[a_i,a_j)$                    & $[\tilde a_j, \tilde a_i)$  & $[a_i,\tilde a_j)$                & $[a_j, \tilde a_i)$     \\
\hline
\hline
$\DgZZ(f)$                & $[X_{i-1}^i,X_{j-1}^{j-1}]$    & $[X_i^i,X_{j-1}^j]^-$       & $[X_{i-1}^i,X_{j-1}^j]$           & $[X_i^i,X_{j-1}^{j-1}]^-$ \\
\hline
Type                      & I                              & II                          & III                               & IV \\
\hline
\end{tabular}
\end{table}
 
In particular, the bottleneck distance, as well as stability results, can be derived for levelset zigzag persistence barcodes using
this correspondence and Theorem~\ref{th:Stab}.

\subsection{Reeb Graphs}
\label{sec:ReebGraphDef}

\begin{definition}
Given a topological space $X$ and a continuous function 
$f:X\rightarrow\R$, we define the equivalence relation $\sim_f$ between points of $X$ by:
\[
x\sim_f y\ \Longleftrightarrow\ \left[ f(x)=f(y)\emph{ and }x,y\ \emph{belong to the same connected component of}\ f^{-1}(f(x))=f^{-1}(f(y))\right].
\]
The \emph{Reeb graph} $\Reeb_f(X)$ is the quotient space $X/\sim_f$.
\end{definition}

\begin{figure}[htb] 
\begin{center}
\includegraphics[height = 3.5cm]{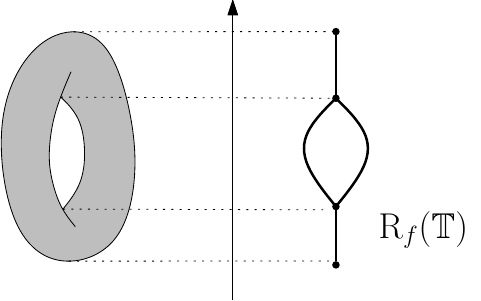}
\caption{ We consider the height function of the torus $\mathbb T$. Note how the critical points induce changes on the graph.}
\label{fig:ReebGraphTorus}
\end{center} 
\end{figure}

As $f$ is constant on equivalence classes, there is an induced map $\tilde{f}~:~\Reeb_f(X)~\rightarrow~\R$ such that $f~=~\tilde{f}~\circ~\pi$,
where $\pi$ is the quotient map $X\to \Reeb_f(X)$:
\begin{equation}\label{eq:triangleReeb}
\xymatrix@!C=5pt{
X \ar[rr]^-\pi \ar[dr]_-f && \Reeb_f(X) \ar[dl]^-{\tilde f}\\
& \R
}
\end{equation}
If $f$ is a function of Morse type, then
the pair $(X,f)$ is an {\em $\R$-constructible space} in the sense
of~\cite{deSilva16}.  This ensures that the Reeb graph is a
multigraph, whose nodes are in one-to-one correspondence with the connected components
of the critical level sets of $f$. 
In that case, computing the Reeb graph of a Reeb graph preserves all information,
as stated in the following remark.

\begin{rmq}\label{rem:idemreeb}
Let $f$ be a Morse-type function.
Then there is a bijection $b:\Reeb_{\tilde f}(\Reeb_f(X))\rightarrow\Reeb_f(X)$
such that $\tilde f\circ b = \tilde{\tilde f}$. In other words, computing the Reeb graph is an idempotent operation.
\end{rmq}

In the following, the combinatorial version of the Reeb graph (where each critical point is turned into a node) is
denoted by $\mathcal{C}\Reeb_f(X)$.

\begin{figure}[htb]
\centering
\includegraphics[width=10cm]{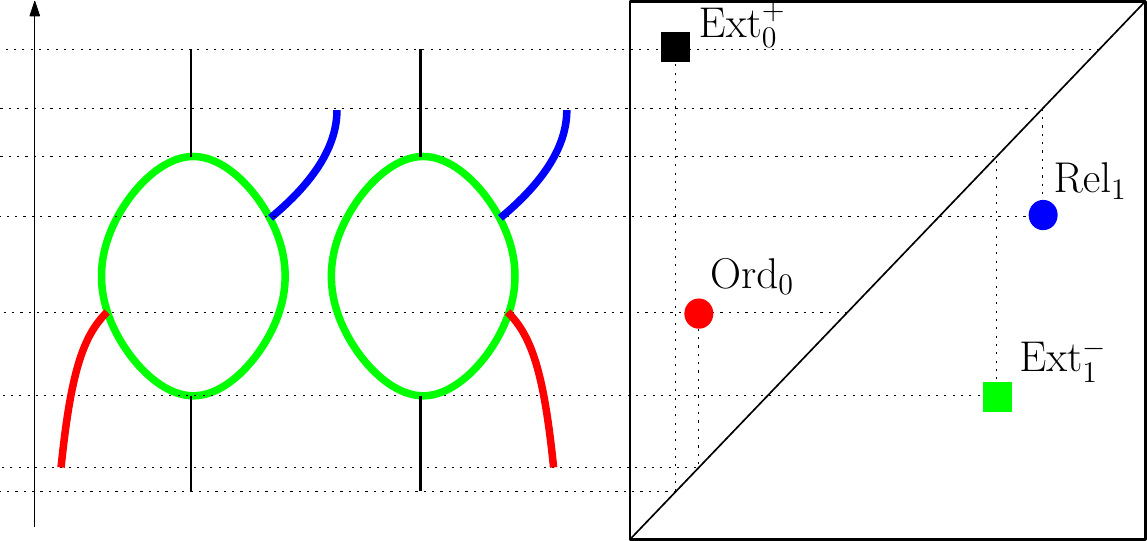}
\caption{Two Reeb graphs with the same set of features but not the same layout.}
\label{fig:Reeb_struct}
\end{figure}

\paragraph{Persistence-based bag-of-features signature.} 
There is a nice interpretation of $\Dg(\tilde f)$ in terms of the
structure of $\Reeb_f(X)$. We refer the reader to~\cite{Bauer14} and
the references therein for a full description as well as formal definitions
and statements.  Orienting the Reeb graph vertically so ${\tilde f}$ is the height function,
we can see each connected component of the graph as a trunk with
multiple branches (some oriented upwards, others oriented downwards)
and holes.  Then,
one has the following correspondences, where the
{\em vertical span} of a feature is the span of its image by~$\tilde
f$:
\begin{itemize}
\item The vertical spans of the trunks are given by the points in $\Ext_0^+(\tilde f)$;
\item The vertical spans of the branches that are oriented downwards are given by the points in $\Ord_0(\tilde f)$;
\item The vertical spans of the branches that are oriented upwards are given by the points in $\Rel_1(\tilde f)$;
\item The vertical spans of the holes are given by the points in $\Ext_1^-(\tilde f)$. 
\end{itemize}
The rest of the diagram of~$\tilde f$ is empty. These correspondences
provide a dictionary to read off the structure of the Reeb graph from
the extended persistence diagram of the induced map~$\tilde f$. Note that it
is a bag-of-features type signature, taking an inventory of all
the features (trunks, branches, holes) together with their vertical
spans, but leaving aside the actual layout of the features. As a
consequence, it is an incomplete signature: two Reeb graphs with the
same persistence diagram may not be isomorphic, as illustrated in
Figure~\ref{fig:Reeb_struct}.

\paragraph{Connection to the extended and zigzag persistence of~$f$.}
We now show that the topological structure of $\Reeb_f(X)$ is actually nothing but a simplification 
of the one of $f$. This can be phrased using the extended persistence diagrams of $f$ and $\tilde f$:

%
%
%
\begin{theorem}\label{thm:pdreeb}
Let $X$ be a topological space and
$f:X\rightarrow\R$ be a function of Morse type.
Then, the levelset zigzag persistence barcodes of $f$ and $\tilde f$ in dimension 0 are the same: $\DgZZ_0(f)=\DgZZ_0(\tilde f)$,
and the extended persistence diagram of $\tilde f$ is included in the one of $f$: $\Dg(\tilde{f})\subseteq \Dg(f)$.
More precisely:
\begin{align*}
\Dg_0(\tilde{f}) &= \Dg_0(f) \\
\Dg_1(\tilde{f}) &= \Dg_1(f)\setminus(\Ext^+_1(f)\cup\Ord_1(f)) \\
\Dg_p(\tilde{f}) &= \emptyset\text{  if }p\geq 2
\end{align*}
\end{theorem}

Note that  $\Ext_0^-(\tilde f)
= \emptyset$ because every
essential $0$-dimensional feature corresponds to some connected component of the
domain, and it is born at the minimum function value and killed at the
maximum function value over that connected component, hence it belongs to $\Ext_0^+$.
Similarly, $\Rel_0(\tilde f)=\emptyset$
because no 0-dimensional homology class (i.e. connected component) can be
created in the relative part of the extended filtration of~$f$. 
Hence, the structure of a Reeb graph can be read off from the levelset zigzag persistence module of $\tilde f$.
Indeed, since $\Ext_1^+(\tilde f)$, $\Ord_1(\tilde f),\Ext_0^-(\tilde f),\Rel_0(\tilde f)$ and $\Dg_p(\tilde f)$ for $p\geq 2$ are empty, 
it follows from Corollary~\ref{cor:EPZZ} that 
there is a bijection preserving types between $\Dg_0(\tilde f)\cup\Dg_1(\tilde f)$ and $\DgZZ_0(\tilde f)$. 
This is because all intervals in the 1-dimensional extended persistence module of $\tilde f$ 
are either of type $\Rel$ or $\Ext^-$, and thus their analogues in the levelset zigzag persistence module of $\tilde f$ have 
homological dimension 0 according to Table~\ref{tab:EPZZ}.

We provide a proof for completeness, as we have not seen this result  stated formally in the
literature. 
First, note that $\Crit(f)=\{a_1,\cdots,a_n\}=\Crit(\tilde f)$. 
Hence, given $i\leq j$ and $[s_i,s_j]$ as in Section~\ref{sec:ZZ},
we recall that $X_i^j$ denote $X^{[s_i,s_j]}=f^{-1}([s_i,s_j])$ and $\Reeb_f(X)_i^j$ denote $\Reeb_f(X)^{[s_i,s_j]}=\tilde f^{-1}([s_i,s_j])$.

\begin{lemma}\label{lem:ZZ0}
Let $\pi$ denote the quotient map $X\rightarrow\Reeb_f(X)$.
Let 
$i\leq j$, and 
$I=[s_i,s_j]$, as defined in Section~\ref{sec:ZZ}.
Then the morphism $\pi_*:H_0(X^I)\rightarrow H_0(\Reeb_f(X)^I)$ 
is an isomorphism.
\end{lemma}

The proof of Lemma~\ref{lem:ZZ0} is
simpler when $\pi$ admits {\em continuous sections}, i.e. when there exist
continuous maps $\sigma:\Reeb_f(X)\rightarrow X$ such that
$\pi\circ\sigma=\text{id}_{\Reeb_f(X)}$.
Below we give the proof under this hypothesis, deferring 
the general case of Morse-type functions to
Appendix~\ref{sec:proof_pdreeb}.  The hypothesis holds for instance
when $X$ is a compact smooth manifold and $f$ is a Morse function, or
when $X$ is a simplicial complex and $f$ is piecewise-linear.

\begin{proof}[Lemma~\ref{lem:ZZ0}]

Since $\pi$ is surjective, proving the result boils down to
showing that $x,y$ are connected in 
$X^j_i$ if and only if $\pi(x),\pi(y)$ are connected in $\Reeb_f(X)_i^j$.

\begin{itemize}
\item If $x,y$ are connected in $X^j_i$,
then $\pi(x),\pi(y)$ are connected in $\Reeb_f(X)^j_i$
by continuity of $\pi$ and
commutativity of~(\ref{eq:triangleReeb}).
\item If $\pi(x),\pi(y)$ are connected in $\Reeb_f(X)^j_i$,
then choose a path $\gamma$ connecting $\pi(x)$ to $\pi(y)$. By definition of~$\sigma$, we have
$\pi\circ\sigma\circ\pi(x)=\pi(x)$, thus $\sigma\circ\pi(x)$ and
$x$ lie in the same connected component of $f^{-1}(f(x))$.  Let $\gamma_x$ be a path
connecting $x$ to $\sigma\circ\pi(x)$. Similarly, let $\gamma_y$
be a path connecting $\sigma\circ\pi(y)$ to $y$. Then,
$\gamma_y\circ\sigma(\gamma)\circ\gamma_x$ is a path between $x$ and $y$ in $X^j_i$.
\end{itemize}
\qed\end{proof}

\begin{proof}[Theorem~\ref{thm:pdreeb}]

We first show that $\DgZZ_0(f)=\DgZZ_0(\tilde f)$. Let $\pi$ denote the quotient map $X\rightarrow\Reeb_f(X)$.
Since $\pi$ is continuous, it induces a morphism in homology $\pi_*$. We will show that $\pi_*$ 
induces an isomorphism between $\LZZ(f)$ and $\LZZ(\tilde f)$ in dimension~0.
First, note that $\Crit(f)=\{a_1,\cdots,a_n\}=\Crit(\tilde f)$. 
Hence both $\LZZ(f)$ and $\LZZ(\tilde f)$ have $2n+1$ nodes.
Now, let $1\leq i\leq n$.

\begin{itemize}
\item According to Lemma~\ref{lem:ZZ0}, $\pi_*:H_0(X_i^i)\rightarrow H_0(\Reeb_f(X)^i_i)$ is an isomorphism, and 
the same holds for $\pi_*:H_0(X_i^{i+1})\rightarrow H_0(\Reeb_f(X)^{i+1}_i)$. Hence $\pi_*$ induces a pointwise isomorphism 
in dimension~0 between $\LZZ(f)$ and $\LZZ(\tilde f)$.

\item Let $\iota:X^i_i\rightarrow X_i^{i+1}$ and
 $\iota^{\Reeb}:\Reeb_f(X)^{i}_i\rightarrow \Reeb_f(X)^{i+1}_i$
be canonical inclusions. Then, we have $\pi\circ\iota=\iota^{\Reeb}\circ\pi$ 
by definition of $\iota^\Reeb$.
Hence, the following diagram commutes: 
$$\xymatrix@!C{
H_0(X^i_i) \ar[r]^-{\iota_*} \ar[d]_-{\pi_*} & H_0(X^{i+1}_i) \ar[d]^-{\pi_*} \\
H_0(\Reeb_f(X)^i_i) \ar[r]_-{\iota^{\Reeb}_*} & H_0(\Reeb_f(X)^{i+1}_i)
}$$
and the same is true for the canonical inclusions $X^i_{i-1}\hookleftarrow X_i^{i}$ and
 $\Reeb_f(X)^{i}_{i-1}\hookleftarrow \Reeb_f(X)^{i}_i$.
\end{itemize}

Hence, the induced pointwise isomorphism is an isomorphism between $\LZZ_0(f)$ and $\LZZ_0(\tilde f)$.

Now, recall that there is a bijection $b_1$ preserving types between $\Dg(\tilde f)$ and $\DgZZ_0(\tilde f)$.
Since there is also a bijection $b_2$ preserving types between $\DgZZ_0(\tilde f)$ and $\DgZZ_0(f)$
and a bijection $b_3$ preserving types between $\DgZZ_0(f)$ and $\Ord_0(f)\cup\Ext_0^+(f)\cup\Rel_1(f)\cup\Ext_1^-(f)$ from Corollary~\ref{cor:EPZZ},
the result follows by considering the bijection $b_3\circ b_2 \circ b_1$. 

\qed\end{proof}

\paragraph{A distance between Reeb graphs.}
We now give the definition of the {\em functional distortion distance}~\cite{Bauer14} between Reeb graphs.
Note that any Reeb graph $\Reeb_f(X)$ can be equipped with a canonical metric: 
$d_f(x,x')={\rm min}_{\pi:x\rightarrow x'}\ \{{\max}_{t\in[0,1]}\ \tilde f\circ\pi(t)-{\min}_{t\in[0,1]}\ \tilde f\circ\pi(t)\}$, 
where $\pi:[0,1]\rightarrow\Reeb_f(X)$ ranges over the continuous paths from $x$ to $x'$ ($\pi(0)=x$ and $\pi(1)=x'$).
Then, given a pair of Reeb graphs, the functional distortion distance measures the {\em distortion} of their corresponding metrics.
Hence, it is very similar to the Gromov-Hausdorff distance. 
We use this distance in Section~\ref{sec:dfdconv} to provide a convergence result of the (MultiNerve) Mapper to the Reeb graph.

\begin{definition}
Let $X,Y$ be topological spaces and $f:X\rightarrow\R$ and $g:Y\rightarrow\R$ be continuous scalar functions.
The {\em functional distortion distance} between $\Reeb_f(X)$ and $\Reeb_g(Y)$ is:
\begin{equation}\label{eq:dfd}
\distfd(\Reeb_f(X),\Reeb_g(Y))=\underset{\phi,\psi}{\rm inf}\ {\rm max}\left\{\frac{1}{2}D(\phi,\psi),\ \|\tilde f-\tilde g\circ\phi\|_\infty,\ \|\tilde f\circ\psi-\tilde g\|_\infty\right\},
\end{equation}
where: 
\begin{itemize}
\item $\phi:\Reeb_f(X)\rightarrow\Reeb_g(Y)$ and $\psi:\Reeb_g(Y)\rightarrow\Reeb_f(X)$ are continuous maps, 
\item $D(\phi,\psi)=\ {\rm sup}\ \left\{|d_f(x,x')-d_g(y,y')| : (x,y),(x',y')\in C(\phi,\psi)\right\},$
\item $C(\phi,\psi)=\{(x,\phi(x)):x\in\Reeb_f(X)\}\cup\{(\psi(y),y):y\in\Reeb_g(Y)\}$.
\end{itemize}
\end{definition}

The functional distortion distance enjoys the following stability theorem:

\begin{theorem}[Theorem 4.1 in~\cite{Bauer14}]\label{th:dfdstab}
Let $X$ be a topological space and let $f,g:X\rightarrow\R$ be two Morse-type functions with continuous sections.
Then:
$$\distfd(\Reeb_f(X),\Reeb_g(X))\leq \|f-g\|_\infty.$$
\end{theorem}

Since $\distfd$ can be quite hard to compute and to interpret,
we also study the bottleneck distance between the extended persistence diagrams $\distb(\Dg(\tilde f),\Dg(\tilde g))$ in Section~\ref{sec:structure}.
Recall that $\distb$ is only a pseudometric---see Figure~\ref{fig:Reeb_struct}. However, it can be computed efficiently, it allows for interpretation
(recall that extended persistence diagrams act as bag-of-feature signatures) and it has been proven~\cite{Carriere17a} that $\distb$ and 
$\distfd$ are actually equivalent for close Reeb graphs. 

\subsection{Covers and Nerves}
\label{sec:basicdef}

Let $Z$ be a topological space. A \emph{cover} of $Z$ is a 
family $\U$ of subsets of $Z$, $\U=\{U_\alpha\}_{\alpha\in A}$, such
that $Z=\bigcup_{\alpha\in A}\ U_\alpha$. It is {\em open} if all its
elements are open subspaces of $Z$. It is {\em connected} if all its
elements are connected subspaces of $Z$. Its \emph{nerve} is the 
abstract simplicial complex $\mathcal{N}(\U)$ that has one $k$-simplex per $(k+1)$-fold intersection of elements of $\U$:
\[
\{\alpha_0,...,\alpha_k\}\in\mathcal{N}(\U)\Longleftrightarrow\bigcap_{i=0,...,k}U_{\alpha_i}\neq\emptyset.
\]
%
When $\V$ itself is a cover of $Z$, it is called
a \emph{subcover} of $\U$. It is {\em proper} if it is not equal to $\U$. 
Finally, $\U$ is called \emph{minimal} if it admits no proper
subcover or, equivalently, if it has no element included in the union
of the other elements.
Given a minimal cover $\U=\{U_\alpha\}_{\alpha\in A}$, for every
$\alpha\in A$ we let
\[
\tilde{U}_\alpha = U_\alpha \setminus \bigcup_{\alpha'\neq\alpha\in A} U_{\alpha'},
\]
be the \textit{proper subset} of $U_\alpha$,
that is the maximal subset of $U_\alpha$ that has an empty
intersection with the other elements of $\U$.
$\U$ is called \textit{generic} if no connected component of the 
proper subsets of its elements
is a singleton. \\

Consider now the special case where $Z$ is a subset of $\R$, equipped
with the subspace topology. A subset $U\subseteq Z$ is an {\em
  interval of $Z$} if there is an interval $I$ of $\R$ such that
$U=I\cap Z$. Note that $U$ is open in $Z$ if and only if $I$ can be
chosen open in $\R$.  A cover $\U$ of $Z$ is an {\em interval cover}
if all its elements are intervals. In this case, ${\rm End}(\U)$
denotes the set of all of the interval endpoints. Finally, the
\textit{granularity} of $\U$ is the
supremum of the lengths of its elements, i.e. it is the quantity $\text{sup}_{U\in\U}\ |U|$ where
$|U|=\sup(U)-\inf(U)\in\R\cup\{+\infty\}$.
%
%
%
\begin{lemma}\label{lem:minimal-cover_1-dim}
No more than
two elements of a minimal open interval cover can intersect at a time.
\end{lemma}
\begin{proof}
Assume for a contradiction that there are $k\geq 3$ elements
of $\U$: $U_1, \cdots, U_k$, that have a non-empty common intersection. For
every $i$, fix an open interval $I_i$ of $\R$ such that $U_i=I_i\cap Z$. Up to
a reordering of the indices, we can assume without loss of generality
that $I_1$ has the smallest lower bound and $I_2$ has the largest
upper bound. Since $I_1\cap I_2 \supseteq U_1\cap U_2 \neq\emptyset$,
the remaining intervals satisfy $I_i\subseteq I_1\cup I_2$. In
particular, we have $U_3= I_3\cap Z\subseteq (I_1\cup I_2)\cap Z =
(I_1\cap Z)\cup (I_2\cap Z) = U_1\cup U_2$, so the cover $\U$ is not
minimal.
\qed\end{proof}
\begin{lemma}
If $Z$ is $\R$ itself or a compact subset thereof, then any cover $\U$
of $Z$ has a minimal subcover.
\end{lemma}
\begin{proof}
When $Z$ is compact, there exists a subcover $\V$ of $\U$ that has
finitely many elements.  Any subcover of $\V$ with the minimum number of
elements is then a minimal cover of $Z$.

When $Z=\R$, the same argument applies to any subset of the form $[-n, n]$, $n\in \N$. Then, a simple induction on $n$ allows us to build a minimal subcover of $\U$. 
\qed\end{proof}

From now on, unless otherwise stated, all covers of $Z\subseteq \R$
will be generic, open, minimal, interval covers ({\em gomic} for
short). Given such a cover~$\U$, the proper subset $\tilde U$ of any
interval $U\in\U$ is itself an interval of~$Z$ since $\U$ is generic, therefore we call it
the \textit{proper subinterval} of~$U$.
Moreover, Lemma~\ref{lem:minimal-cover_1-dim} yields a total order on the intervals of~$\U$, 
so each one of them partitions into subintervals as follows:
\begin{equation}\label{eq:U_partition}
U=U_\cap^-\sqcup \tilde U\sqcup U_\cap^+,
\end{equation}
where $U_\cap^-$ is the intersection of $U$ with the element right below it in the cover ($U_\cap^-=\emptyset$ if that element does not exist), 
and where $U_\cap^+$ is the intersection of $U$ with the element right above it ($U_\cap^+=\emptyset$ if that element does not exist).

\subsection{Mapper}
\label{sec:ClassicalMapperDef}

Let $f:X\rightarrow Z$ be a continuous function. 
Consider a cover $\U$ of $\im(f)$, and pull
it back to $X$ via $f^{-1}$. Then, decompose every
$V_\alpha=f^{-1}(U_\alpha)\subseteq X$ into its connected components:
$V_\alpha=\bigsqcup_{i\in\{1,\cdots,c(\alpha)\}}V_\alpha^i$, where
$c(\alpha)$ is the number of connected components of $V_\alpha$.  Then,
$\V=\{V_\alpha^i\}_{\alpha\in A,i\in\{1,\cdots,c(\alpha)\}}$ is a connected
cover of $X$.  It is called the \textit{connected pullback cover}, and
its nerve $\mathcal{N}(\V)$ is the Mapper.

\begin{definition}
Let $X,Z$ be topological spaces, $f:X\rightarrow Z$ be a continuous
function, $\U$ be a cover of $\im(f)$ and $\V$ be the associated
connected pullback cover.  The \emph{Mapper} of $X$ is
$\Mapper_f(X,\U)=\mathcal{N}(\V)$.
\end{definition}

\begin{figure}[h]\centering
\includegraphics[height = 5cm]{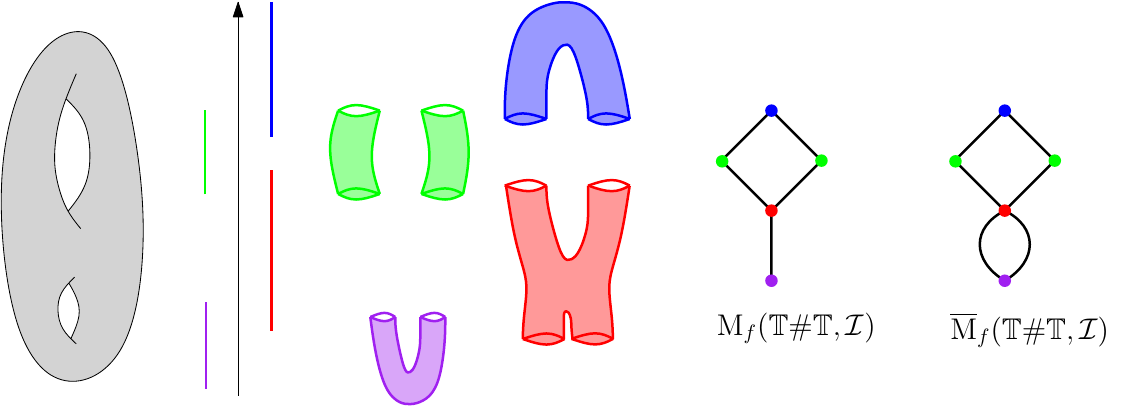}
\caption{Example of the Mapper and the MultiNerve Mapper computed on the double torus $\mathbb T\# \mathbb T$ with the height function $f$. 
The cover $\I$ of $\im(f)$ has four intervals (red, green, blue and purple), and the cover of the double torus 
has five connected components (one is blue, one is red, one is purple and the other two are green). The Mapper and the MultiNerve Mapper are displayed on the right.}
\label{fig:mmappervsmappertorus}
\end{figure}

See Figure \ref{fig:mmappervsmappertorus} for an illustration. Note
that, when $Z=\R$ and $\U$ is a gomic, the Mapper has a
natural 1-dimensional stratification since no more than two intervals
can intersect at a time by Lemma~\ref{lem:minimal-cover_1-dim}. Hence,
in this case, it has the structure of a (possibly infinite) simple graph
and therefore has trivial homology in dimension 2 and above.

\section{MultiNerve Mapper}
\label{MultiNerve}

In this section, we explain how to extend the construction of Mapper by using a slight modification of the nerve, called the {\em multinerve},
and whose definition relies on \textit{simplicial posets}~\cite{Verdiere12}.

\subsection{Simplicial Posets and MultiNerves}

\begin{definition}
A \emph{simplicial poset} is a partially ordered set $(P,\preceq)$, whose elements are called \emph{simplices}, and which satisfies the two following properties:

\begin{itemize}
\item[\rm (i)] $P$ has a least element called $0$ such that $\forall p\in P$, $0\preceq p$;
\item[\rm (ii)] $\forall p\in P$, $\exists d\in\mathbb{N}$ such that the lower segment 
$[0,p]=\{q\in P:q\preceq p\}$ is \emph{isomorphic} to the
set of simplices of the standard $d$-simplex with the inclusion as partial order,
where an isomorphism between posets is a bijective and order-preserving function.
\end{itemize}
\end{definition}

\begin{figure}[h]\centering
\includegraphics[height = 6cm]{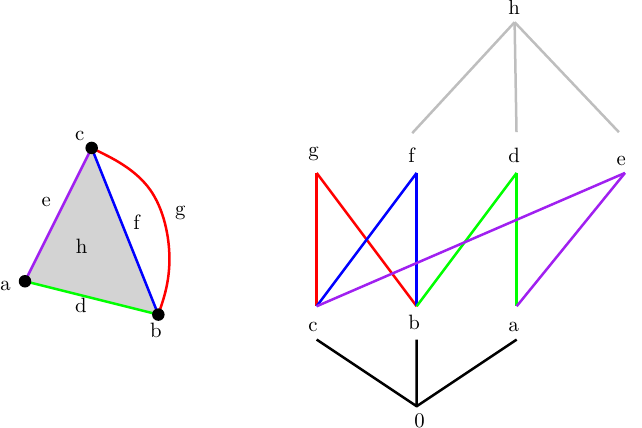}
\caption{Left: A simplicial poset that is not a 
simplicial complex. Indeed, edges $f$ and $g$ have the same vertices ($b$ and $c$). 
Right: The corresponding Hasse diagram showing the partial order on the simplices. Note that $f,g$
cannot be part of the same 2-cell.}
\label{fig:PosetExample}
\end{figure}

Simplicial posets are extensions of simplicial complexes: while every
simplicial complex is also a simplicial poset (with inclusion as
partial order and $\emptyset$ as least element), the converse is not
always true as different simplices may have the same set of vertices.
However, these simplices cannot be faces of the same
higher-dimensional simplex, otherwise $\rm (ii)$ would be false. See
Figure~\ref{fig:PosetExample} for an example of a simplicial poset
that is not a simplicial complex.

Given a cover $\U$ of a topological space $X$, the nerve is extended
to a simplicial poset as follows:

\begin{definition}
Let $\U=\{U_\alpha\}_{\alpha\in A}$ be a cover of a topological space $X$.
The \emph{multinerve $\M(\U)$}
is the simplicial poset defined by:
\[
\M(\U)=\left\{(\{\alpha_0, \cdots, \alpha_k\} ,C)\ :\ \bigcap_{i=0}^k U_{\alpha_i}\neq\emptyset\ \mbox{and}\ C\ \text{is a connected component of}\ \bigcap_{i=0}^k U_{\alpha_i}\right\}.
\]
\end{definition}

The proof that this set, together with the least element $(\emptyset,\bigcup_{\alpha\in A}U_\alpha)$
and the partial order $(F,C)\preceq(F',C')\Leftrightarrow F\subseteq F'\ \text{and}\ C'\subseteq C$, 
is a simplicial poset, can be found in \cite{Verdiere12}.
Given a simplex $(F,C)$ in the multinerve of a cover, its \textit{dimension}
is $|F|-1$. 
The dimension of the multinerve of a cover is the maximal dimension of
its simplices.  Given two simplices $(F,C),(F',C')$, we say that $(F,C)$ is a
\textit{face} of $(F',C')$ if $(F,C)\preceq (F',C')$.  



Given a connected pullback cover~$\V$, we extend the Mapper by using the
multinerve $\M(\V)$ instead of $\mathcal{N}(\V)$.  This variant will
be referred to as the MultiNerve Mapper in the following.

\begin{definition}
Let $X,Z$ be topological spaces, $f:X\rightarrow Z$ be a continuous function,
$\U$ be a cover of $\im(f)$ and $\V$ be the associated connected pullback cover.
The \emph{MultiNerve Mapper} of $X$ is $\MMapper_f(X,\U)=\mathcal{M}(\V)$.  
\end{definition}

See Figure~\ref{fig:mmappervsmappertorus} for an illustration.  For
the same reasons as Mapper, when $Z=\R$ and $\U$ is a gomic of
$\im(f)$, the MultiNerve Mapper is a (possibly infinite) multigraph
having trivial homology in dimension 2 and above.
Contrarily to the Mapper, the
MultiNerve Mapper also takes the connected components of the
intersections into account in its construction. As we shall see in
Section~\ref{sec:structure}, it is able to capture the same features
as the Mapper but with coarser gomics, and it is more naturally
related to the Reeb graph.

\subsection{Connection to Mapper}

The connection between the Mapper and the MultiNerve Mapper is
induced by the following connection between nerves and multinerves:

\begin{lemma}[\cite{Verdiere12}]\label{lem:nerve-vs-multinerve}
Let $X$ be a topological space and $\U$ be a cover of $X$.
Let $\pi_1:(F,C)\mapsto F$ be the projection of the simplices of
$\M(\U)$ onto the first coordinate. Then, $\pi_1(\M(\U))=\mathcal{N}(\U)$. 
\end{lemma}
%
%
\begin{corollary}
\label{cor:projMNonM}
Let $X,Z$ be topological spaces and $f:X\rightarrow Z$ continuous.
Let $\U$ be a cover of $\im(f)$. 
Then, $\Mapper_f(X,\U)=\pi_1(\MMapper_f(X,\U))$.
\end{corollary} 

\begin{figure}[h]\centering
\includegraphics[height=4cm]{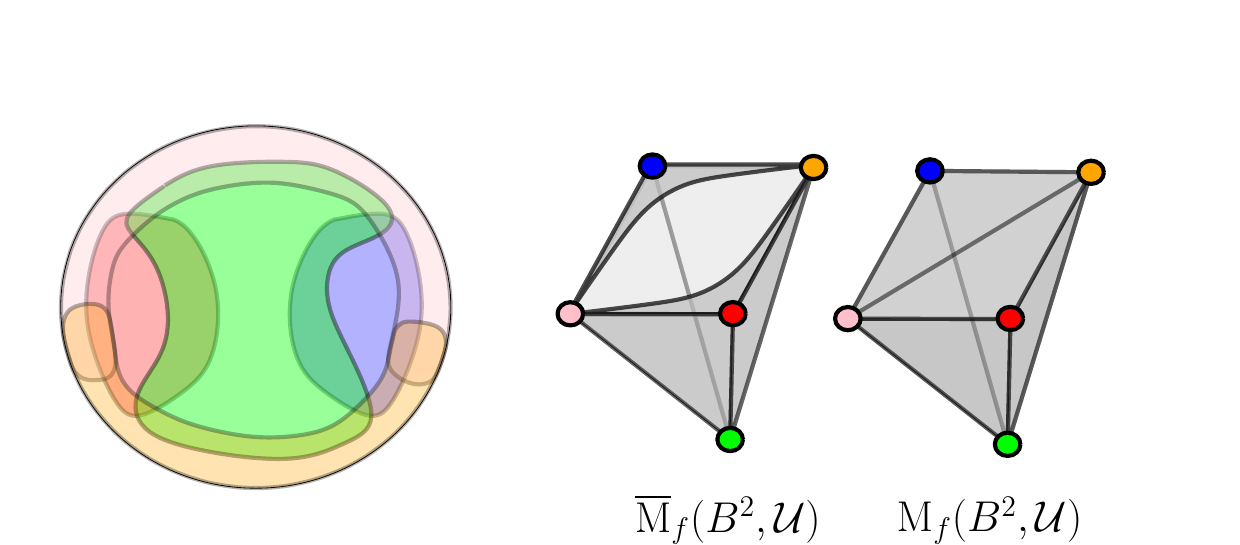}
\caption{The domain is the disk $B^2$, and we consider the identity function $f$,
as well as a generic open minimal cover $\U$ with five elements.
The MultiNerve Mapper is homeomorphic to the disk $B^2$ and the Mapper is homeomorphic to the sphere $S^2$.
Then, $H_2(\Mapper_f(B^2,\U))\neq 0$ while $H_2(\MMapper_f(B^2,\U))=0$.}
\label{fig:cover} 
\end{figure}

Thus, when $Z=\R$ and $\U$ is a gomic, the Mapper is the simple graph
obtained by gluing the edges that have the same endpoints in the
MultiNerve Mapper. In this special case it is even possible to
embed $\Mapper_f(X,\U)$ as a subcomplex of $\MMapper_f(X,\U)$.
Indeed, both objects are multigraphs over the same set of nodes
since they are built from the connected pullback cover. Then, it
is enough to map each edge of $\Mapper_f(X,\U)$ to one of its copies
in $\MMapper_f(X,\U)$, chosen arbitrarily, to get a subcomplex. This mapping serves as a simplicial section for the projection $\pi_1$, therefore:
\begin{lemma}\label{lem:pi1_surj}
When $Z=\R$ and $\U$ is a gomic, $\pi_1$  induces a surjective homomorphism in homology.
\end{lemma}
Note that this is not true in general when $\MMapper_f(X,\U)$ has a
higher dimension. See Figure~\ref{fig:cover} for an example.
%

\section{Structure of the MultiNerve Mapper}
\label{sec:structure}

In this section, we study and characterize the topological structure of the (MultiNerve) Mapper computed on a non discrete topological space.
More precisely, we show that this topological structure can be read off
from the extended persistence diagram of the Reeb graph. 
To prove this, we show that the MultiNerve Mapper $\MMapper_f(X,\I)$ is actually isomorphic (as a combinatorial multigraph)
to a specific Reeb graph, 
whose extended persistence diagram is related to 
the extended persistence diagram $\Dg(\tilde f)$ of $\Reeb_f(X)$.

\subsection{Topology of the MultiNerve Mapper}\label{sec:topoMultiNerve}

In order to show that the MultiNerve Mapper is a specific Reeb graph, 
we first show that (MultiNerve) Mappers can be equipped with functions. 

\begin{definition}\label{def:arbitfunc}
Let $\I=\{I_\alpha\}_{\alpha\in A}$ be a gomic of $\im(f)$ and 
$\V=\{V^i_\alpha\}_{\substack{1\leq i \leq c(\alpha), \alpha\in A}}$ be the associated
connected pullback cover. 
Then we define $\fMM:\MMapper_f(X,\I)\rightarrow\R$ 
as the piecewise-linear extension of the function defined
on the nodes of $\MMapper_f(X,\I)$ by
$V_\alpha^i\mapsto{\rm mid}(\tilde I_\alpha)$,
where ${\rm mid}(\tilde I_\alpha)$ is the midpoint of the proper subinterval $\tilde I_\alpha$ of $I_\alpha$.
The definition of $\fM:\Mapper_f(X,\I)\rightarrow\R$ is similar. 
\end{definition}

Hence, Reeb graphs can be computed from $\MMapper_f(X,\I)$ and $\Mapper_f(X,\I)$, once they are equipped with $\fMM$ and $\fM$ respectively.
Let us call them $\Reeb_{\fMM}(\MMapper_f(X,\I))$ and $\Reeb_{\fM}(\Mapper_f(X,\I))$,
with corresponding induced maps 
${\fMMR}:\Reeb_{\fMM}(\MMapper_f(X,\I))\rightarrow\R$ and $\fMR:\Reeb_{\fM}(\Mapper_f(X,\I))\rightarrow\R$.
The following lemma, which states that (MultiNerve) Mappers are isomorphic to their Reeb graphs, is a simple consequence of Remark~\ref{rem:idemreeb}.

\begin{lemma}\label{lem:idemMapper}
Let $X$ be a topological space and $f:X\rightarrow\R$ be a Morse-type function.
Let $\I$ be a gomic of $\im(f)$.
Then $\MMapper_f(X,\I)$ and $\mathcal C\Reeb_{\fMM}(\MMapper_f(X,\I))$ are isomorphic as combinatorial multigraphs.
The same is true for $\Mapper_f(X,\I)$ and $\mathcal C\Reeb_{\fM}(\Mapper_f(X,\I))$.
\end{lemma}  


Hence, by a slight abuse of notation, we rename $\tilde{\sf m}_{\I}$ 
and  $\tilde{\bar{\sf m}}_{\I}$ 
into $\fM$ and $\fMM$ for convenience.


%


We now state the main result of this section, which ensures that the extended persistence diagram $\Dg(\fMM)$, i.e. the bag-of-features signature of 
$\Reeb_{\fMM}(\MMapper_f(X,\I))$ and $\MMapper_f(X,\I)$,  
is nothing but a simplification of  $\Dg(\tilde f)$, i.e. the bag-of-features signature of $\Reeb_{f}(X)$.

\begin{figure}[h]\centering
\includegraphics[height = 4cm]{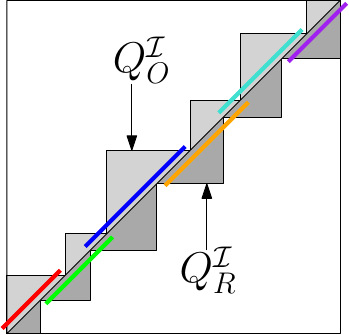}\ \ \ \ \includegraphics[height = 4cm]{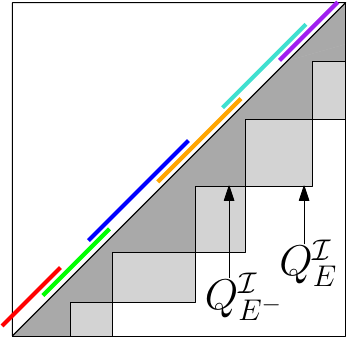}
\caption{Left: Staircases of ordinary (light grey) 
and relative (dark grey) types.
Right: Staircases of extended types---$Q_{E^-}^{\I}$ is in dark grey while $Q_{E}^{\I}$ 
is the union of $Q_{E^-}^{\I}$ with the light grey area.}
\label{fig:staircase}
\end{figure}

\begin{theorem}
\label{th:ExDgMultiNerve}
Let $X$ be a topological space and $f:X\rightarrow\R$ be a Morse-type function. Let $\Reeb_f(X)$ be the corresponding Reeb graph
and $\tilde f : \Reeb_f(X)\rightarrow\R$ be the induced map.
Let $\I$ be a gomic of ${\rm im}(f)$. 
There are bijections between:

\begin{center}
\begin{tabular}{ll}
{\rm (i)} 
$\Ord_0(\fMM)$ 
and $\Ord_0(\tilde f)\setminus Q_O^{\I}$ & 
{\rm (iii)} 
$\Ext_1^-(\fMM)$ 
and $\Ext_1^-(\tilde f)\setminus Q_{E^-}^{\I}$\\[0.5em]
{\rm (ii)} 
$\Rel_1(\fMM)$ 
and $\Rel_1(\tilde f)\setminus Q_R^{\I}$ &
{\rm (iv)} 
$\Ext_0^+(\fMM)$ 
and $\Ext_0^+(\tilde f)$
\end{tabular}
\end{center}
where $\StairOrd=\bigcup_{I\in\I}Q_{{\tilde I}\cup I_\cap^+}^+$,
$\StairRel=\bigcup_{I\in\I}Q_{{\tilde I}\cup I_\cap^-}^-$, and
$\StairExtMM=\bigcup_{I\in\I}Q_I^-$, and where, for any interval $I$ with endpoints $a\leq b$, we let
$Q_I^+ =\{(x,y)\in\R^2\,:\, a\leq x\leq y\leq b\}$ be the corresponding half-square
above the diagonal, and $Q_I^- =\{(x,y)\in\R^2\,:\, a\leq y< x\leq
b\}$ be  the half-square strictly below the diagonal. See Figure~\ref{fig:staircase} for an illustration.
\end{theorem}

The remaining of Section~\ref{sec:topoMultiNerve} is devoted to the proof of Theorem~\ref{th:ExDgMultiNerve}.
In order to state the proof, we first introduce {\em cover zigzag persistence modules}. 

\begin{definition}\label{def:coverZZ}
Let $X$ be a topological space and $f:X\rightarrow\R$ be a Morse-type function.
Let $\I=\{I_\alpha\}_{1\leq\alpha\leq m}$ be a gomic of $\im(f)$, sorted by the natural order defined in Section~\ref{sec:basicdef}. 

Let $\Crit(f)=\{-\infty=a_0,a_1,...,a_n,a_{n+1}=+\infty\}$. 
For any open interval $I$ with left endpoint $a$, 
we define the integers $l(I)$, $r(I)$ by $l(I)=\max\{i\,:\,a_i\leq a\}$
and $r(I)=\max\{l(I),\max\{i\,:\,a_i\in I\}\}$.
Then, we define the {\em cover zigzag persistence module} $\CZZ(f,\I)$ by 
$$\CZZ(f,\I)=H_*\left(X_{l(I_1)}^{r(I_1)}\hookleftarrow X_{l(I_1\cap I_2)}^{r(I_1\cap I_2)}\hookrightarrow X_{l(I_2)}^{r(I_2)}\hookleftarrow X_{l(I_2\cap I_3)}^{r(I_2\cap I_3)}
\hookrightarrow \cdots \hookleftarrow X_{l(I_{m-1}\cap I_m)}^{r(I_{m-1}\cap I_m)}\hookrightarrow X_{l(I_m)}^{r(I_m)}\right),$$
where the $X_i^j$ spaces are as in Section~\ref{sec:ZZ}.
We also let $\DgCZZ(f,\I)$ denote the barcode of this module.
\end{definition}

Note that cover zigzag persistence modules can be isometrically embedded (with the bottleneck distance) 
into the south face of the Mayer-Vietoris half-pyramid. Indeed, each node of $\CZZ(f,\I)$ belongs to this south face.
The only difficulty is that $\CZZ(f,\I)$ may include the same node several times consecutively when there is a sequence of 
consecutive intervals in the gomic that are all included between two consecutive critical values of $f$, i.e. for which $l(I)=r(I)$. 
However, in that case, the corresponding arrows in the module are isomorphisms.
Thus, composing these arrows leaves the resulting barcode unchanged. 


\begin{lemma}\label{lem:CoverZZ}
Let $X$ be a topological space and $f:X\rightarrow\R$ be a Morse-type function.
Let $\I$ 
be a gomic of ${\rm im}(f)$. 
Then, there is a bijection between $\Dg(\fMM)$ and $\DgCZZ_0(f,\I)$. 
\end{lemma}

\begin{proof}

Recall from Corollary~\ref{cor:EPZZ} that it suffices to show that 
$\LZZ_0(\fMM)$ and $\CZZ_0(f,\I)$ are isomorphic as zigzag persistence modules. 
Assume without loss of generality that $\I$ has $m$ elements, with $m\in\N^*$.
First, note 
that $\card(\Crit(\fMM))$ is equal to $m$. Hence, both $\LZZ(\fMM)$ and $\CZZ(f,\I)$ have exactly $2m+1$ nodes.
Moreover, since the MultiNerve Mapper tracks the connected components of the interval and intersection preimages of $f$,
each element of $\LZZ_0(\fMM)$ is of the form $H_0(f^{-1}(I))$, $I\in\I$, or $H_0(f^{-1}(I\cap J))$, $I,J$ consecutive in $\I$.

Let $I\in\I$.
Since $f$ is Morse-type, $X_{l(I)}^{r(I)}$ and $X^I=f^{-1}(I)$ have the same 
homotopy type. Indeed, recall from Section~\ref{sec:ZZ} that there exist $s_{l(I)}$ and $s_{r(I)}$ such that 
$X_{l(I)}^{r(I)}=f^{-1}\left(\left[s_{l(I)},s_{r(I)}\right]\right)$ 
and $s_{l(I)}$ (resp. $s_{r(I)}$) and the left (resp. right) endpoint of $I$
are located between the same consecutive critical values of $f$.  
In particular, $X_{l(I)}^{r(I)}$ and $X^I$ have the same number of connected components, meaning that 
$H_0(X^I)$ and $H_0(X_{l(I)}^{r(I)})$ are isomorphic groups.
The same is also true for any $I\cap J$, $I,J\in\I$.

Hence, we define a canonical pointwise isomorphism $\Psi$ in dimension 0 as follows:
for each node, 
send each connected component of one preimage, or equivalently each generator of one homology group, 
to the connected component of the other preimage which intersects it 
(there is only one since the preimages have the same number of connected components). 
By definition of the MultiNerve Mapper, $\Psi$ commutes with the canonical inclusion.
Hence, $\LZZ_0(\fMM)$ and $\CZZ_0(f,\I)$ are isomorphic. 
\qed\end{proof}

Finally, we  relate the cover zigzag persistence barcode to the extended persistence diagram of the Reeb graph.
Namely, we show that a specific simplification of this extended persistence diagram encodes the same information as the 
cover zigzag persistence barcode. 


\begin{proof}[Theorem~\ref{th:ExDgMultiNerve}]
Again, recall from Corollary~\ref{cor:EPZZ} 
that $\Dg(\tilde f)$ encodes the same information as $\DgZZ_0(\tilde f)$. Hence, since $\Dg(\fMM)$ and $\DgCZZ_0(f,\I)$ are equivalent
from Lemma~\ref{lem:CoverZZ}, we focus on the relation between
$\DgZZ_0(\tilde f)$ and $\DgCZZ_0(f,\I)$.
As mentioned after Definition~\ref{def:coverZZ},
the cover zigzag persistence module $\CZZ(f,\I)$ can be isometrically embedded in the south face of the Mayer-Vietoris half-pyramid.
Hence, we can assume without loss of generality that the set of nodes of $\CZZ(f,\I)$
is a subset of the nodes of a monotone zigzag module $\overline\CZZ(f,\I)$ 
that can be drawn along the south face of the Mayer-Vietoris half-pyramid
by interpolating the elements of $\CZZ(f,\I)$.
Thus, it suffices by Theorem~\ref{th:MY} to study which intervals disappear when going from $\DgZZ_0(\tilde f)$ to $\overline\DgCZZ_0(f,\I)$ and then to $\DgCZZ_0(f,\I)$ 
using the pyramid rules recalled in Figure~\ref{fig:shiftZZ}.

\begin{figure}[h]\centering
\includegraphics[width=12cm]{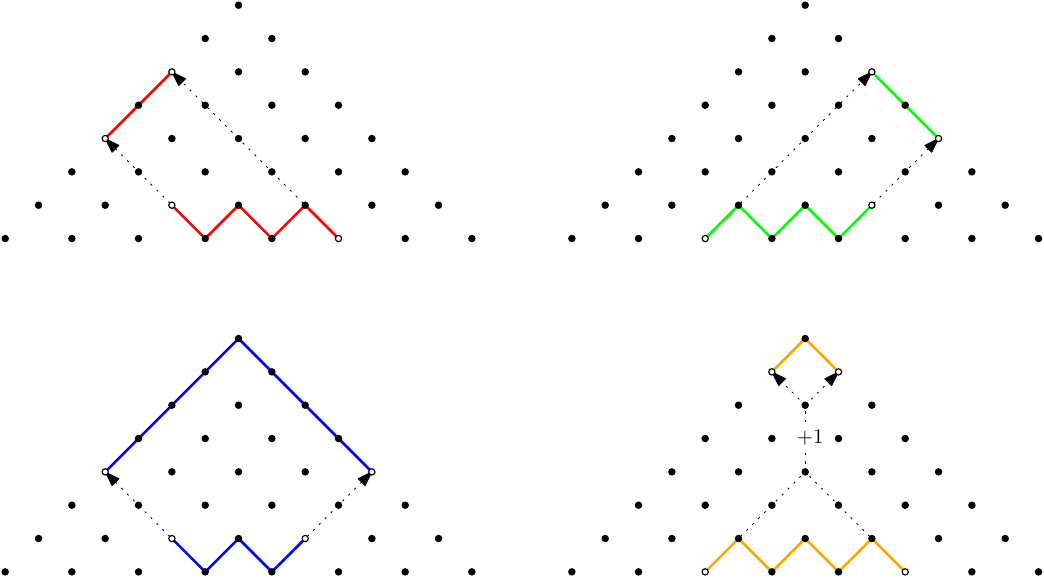}
\caption[Pyramid rules]{(From~\cite{Carlsson09b}) We show the
axis of travel of birth and death endpoints of intervals of
$\LZZ(f)$ to the up-down zigzag persistence module bounding the south face of
the Mayer-Vietoris half-pyramid for interval modules that correspond
to type I intervals (upper-left, red), type II intervals
(upper-right, green), type III intervals (down-left, blue), and type IV
intervals (down-right, orange). The +1 in the
down-right figure means that the homological dimension is increased by one.}
\label{fig:shiftZZ}
\end{figure}

We first give analogues of staircases for zigzag persistence. For any $I=I_\cap^-\sqcup \tilde I \sqcup I_\cap^+\in\I$, we define:
\begin{itemize}
\item ${\rm supp}_O(I)$ as the set of nodes of $\LZZ(f)$ that are located strictly between 
$X_{l(\tilde I\cup I_\cap^+)}^{l(\tilde I\cup I_\cap^+)}$ and $X_{r(\tilde I\cup I_\cap^+)-1}^{r(\tilde I\cup I_\cap^+)}$,
\item ${\rm supp}_R(I)$ as the set of nodes of $\LZZ(f)$ that are located strictly between 
$X_{l(I_\cap^-\cup \tilde I)}^{l(I_\cap^-\cup \tilde I)+1}$ and $X_{r(I_\cap^-\cup \tilde I)}^{r(I_\cap^-\cup \tilde I)}$, 
\item ${\rm supp}_{E^-}(I)$ as the set of nodes of $\LZZ(f)$ that are located strictly between $X_{l(I)}^{l(I)+1}$ and $X_{r(I)-1}^{r(I)}$. 
\end{itemize}

\begin{figure}[h]\centering
\includegraphics[width=15cm]{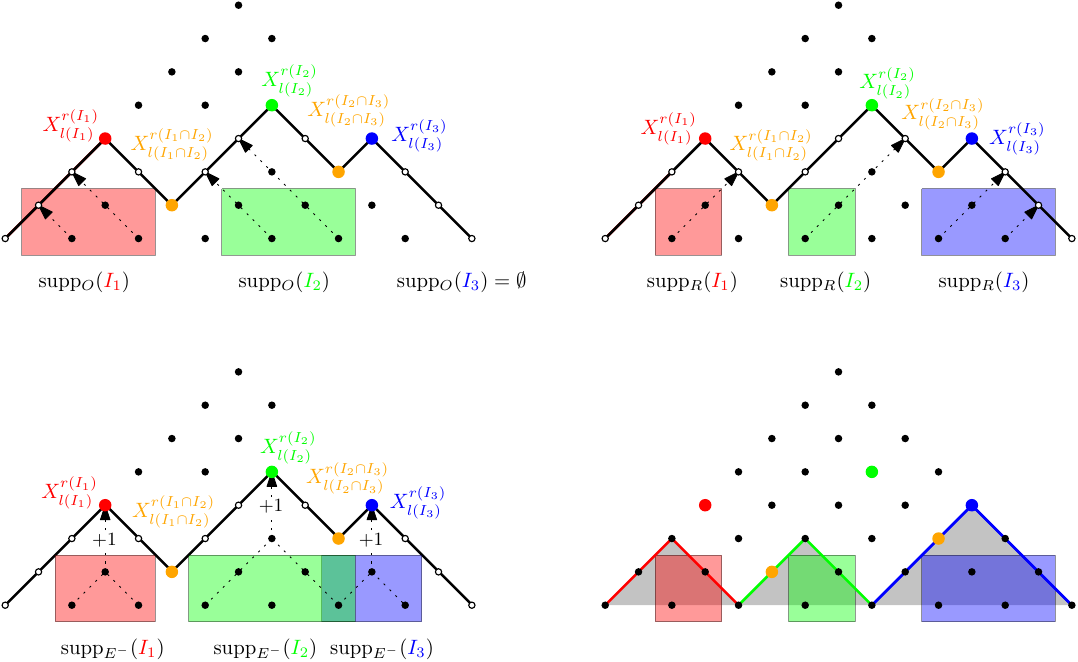}
\caption[Zigzag persistence modules in the half-pyramid]{The black path in the south face of the
  Mayer-Vietoris half-pyramid represents the monotone zigzag persistence module
  $\overline{\rm CZZ}(f,\mathcal I)$ for a gomic $\I$ with three intervals. The
  white disks on this path are the nodes that do not intersect the set
  of nodes of the cover zigzag persistence module $\CZZ(f,\I)$, which are colored
  according to the interval of $\I$ they represent (and are colored
  orange if they represent an intersection). The boxes outline the
  support of the intervals of $\DgZZ_0(f)$ that disappear in the
  MultiNerve Mapper depending on their types (upper-left for
  type I intervals, upper-right for type II intervals and
  down-left for type IV intervals). We also show (down-right) the
  analogue, drawn in grey color, of $Q^{\I}_R$ on the south face of the
  Mayer-Vietoris half-pyramid.}
\label{fig:suppZZ}
\end{figure}

There are two possible ways for an interval of $\DgZZ_0(f)$ to disappear in $\DgCZZ_0(f,\I)$: 
either its homological dimension is shifted by 1, or its intersection with the set of nodes of $\CZZ(f,\I)$ is empty
after being projected onto $\overline\DgCZZ_0(f,\I)$---see Figure~\ref{fig:suppZZ}.
According to the pyramid rules, we have that:
\begin{itemize}
\item Projections of type III intervals of $\DgZZ_0(f)$ onto $\overline\DgCZZ_0(f,\I)$ always intersect with the nodes of $\CZZ(f,\I)$
and their homological dimensions cannot be shifted.
Hence, none of them disappears. This proves (iv).
\item Projections of type IV intervals of $\DgZZ_0(f)$ onto $\overline\DgCZZ_0(f,\I)$ always intersect with the nodes of $\CZZ(f,\I)$.
However, their homological dimensions can be shifted by 1.
This happens when the endpoints collide in the south face of the Mayer-Vietoris half-pyramid. Hence, only those intervals whose support 
is included in ${\rm supp}_{E^-}(I)$ for some $I\in\I$
go through such a shift before getting to $\overline\DgCZZ_0(f,\I)$. This proves (iii).
\item Homological dimensions of type I intervals in $\DgZZ_0(f)$ cannot be shifted, but their projections onto $\overline\DgCZZ_0(f,\I)$ may 
not always intersect with the nodes of $\CZZ(f,\I)$.
This happens for those intervals whose support is included in ${\rm supp}_{O}(I)$ for some $I\in\I$, thus proving (i). 
\item Homological dimensions of type II intervals in $\DgZZ_0(f)$ cannot be shifted, but their projections onto $\overline\DgCZZ_0(f,\I)$ may
not always intersect with the nodes of $\CZZ(f,\I)$. 
This happens for those intervals whose support is included in ${\rm supp}_{R}(I)$ for some $I\in\I$, thus proving (ii). 
\end{itemize}
\qed\end{proof}

\subsection{A signature for MultiNerve Mapper} 
\label{sec:MultiNerveSign}

\begin{figure}[h]\centering
\includegraphics[height=10cm]{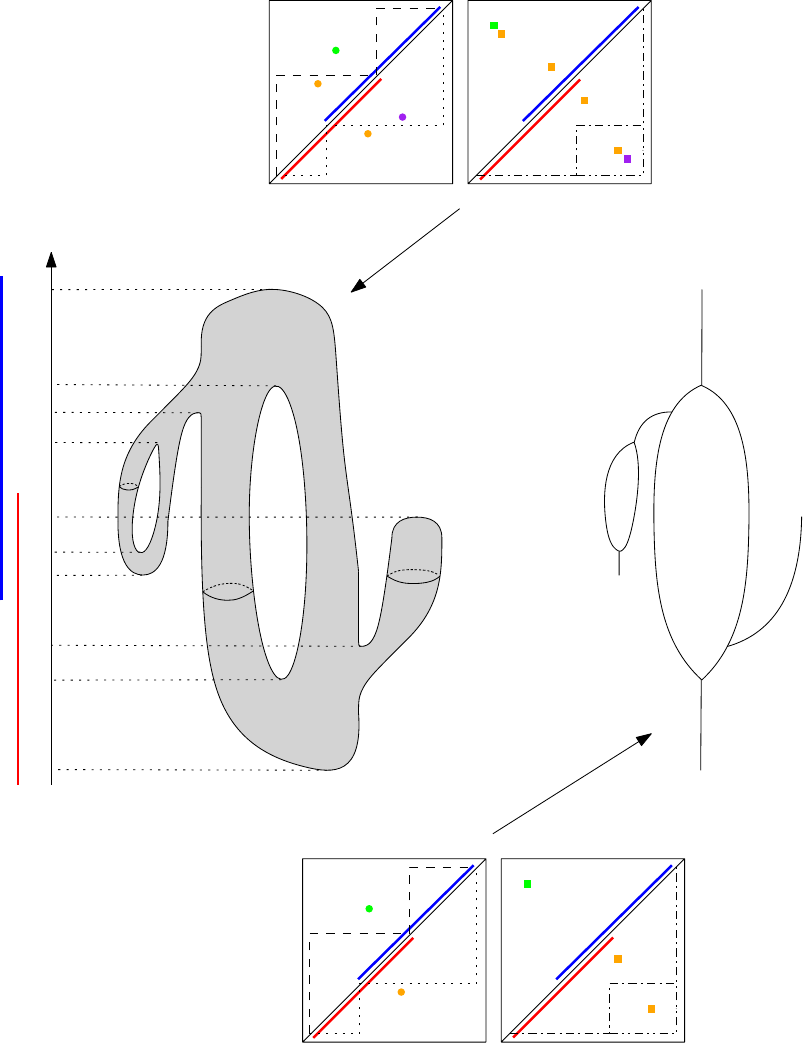}\ \ \ \ \ \ \ \ \includegraphics[height=10cm]{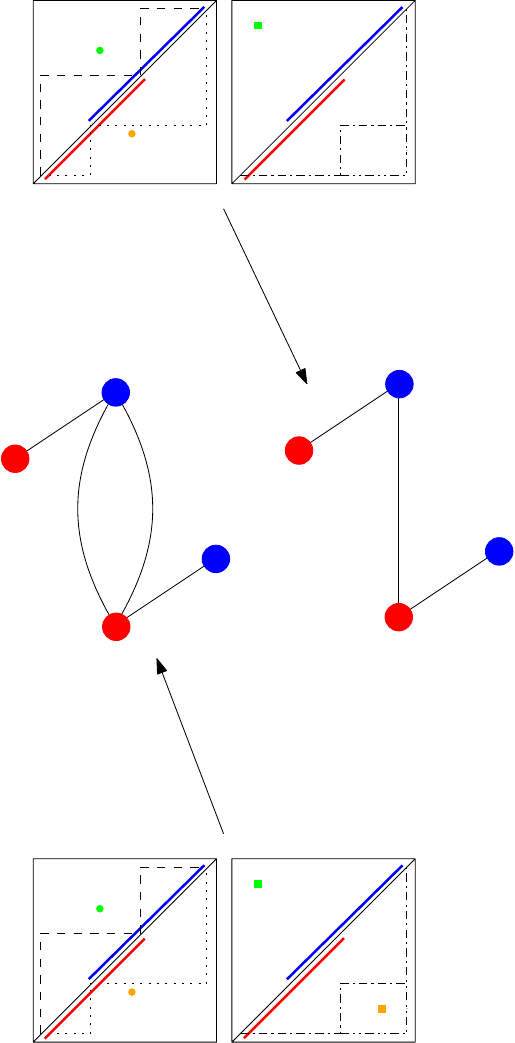}
\caption[Mapper as a pixelization of the Reeb graph]{From left to right: a 2-manifold equipped with the height
  function; the corresponding Reeb graph, MultiNerve Mapper, and
  Mapper. For each object, we display the extended persistence diagrams of
  dimension 0 (green points), 1 (orange points) and 2 (purple
  points). Extended points are squares while ordinary and relative
  points are disks (above and below the diagonal respectively).  The
  staircases are represented with dashed ($Q_O^{\I}$), dotted ($Q_{E^-}^{\I}$),
  dash-dotted ($Q_R^{\I}$), and dash-dot-dotted ($Q_{E}^{\I}$) lines.  
  One can see how to go from the extended
  persistence diagram of the height function to the one of the induced map
  (remove the points in dimension 2 and the points in dimension 1
  above the diagonal), then to the one of the MultiNerve Mapper
  (remove the points inside the staircases corresponding to their
  type), and finally, to the one of the Mapper (remove the extended
  points in $Q_{E}^{\I}$).}
\label{fig:ExDict}
\end{figure}

Theorem~\ref{th:ExDgMultiNerve} means that the dictionary
introduced in Section~\ref{sec:ReebGraphDef} can be used to describe
the structure of the MultiNerve Mapper from the extended persistence
diagram of the induced function~$\tilde f$. 
Indeed, the topological features of $\MMapper_f(X,\I)$
are in bijection with the points of
$\Dg(\tilde f)$ minus the ones that fall into the various
staircases ($Q_O^{\I}$, $Q_{E^-}^{\I}$, $Q_R^{\I}$) corresponding to their type. 
Moreover, by Theorem~\ref{thm:pdreeb}, $\Dg(\tilde f)$ itself is
obtained from $\Dg_0(f)$ and $\Dg_1(f)$ by removing the points of
$\Ext_1^+(f)$ and $\Ord_1(f)$.
Hence, we use the off-staircase part of
$\Dg(\tilde f)$ as a signature for the structure of the MultiNerve
Mapper\footnote{Recall that $\Ext_0^-(f)=\Rel_0(f)=\emptyset$.}:
\begin{align}\label{eq:sign1}
\begin{split}
\Dg(\MMapper_f(X,\I))&=
(\Ord(\tilde{f})\setminus Q_O^{\I})\cup
(\Ext(\tilde{f})\setminus Q_{E^-}^{\I})\cup
(\Rel(\tilde{f})\setminus Q_R^{\I})\\
&=
(\Ord_0(f)\setminus Q_O^{\I})\cup
((\Ext_0^+(f)\cup \Ext_1^-(f))\setminus Q_{E^-}^{\I})\cup
(\Rel_1(f)\setminus Q_R^{\I}).
\end{split}
\end{align}
We call this signature the {\em extended persistence diagram} of the MultiNerve Mapper.  
Note that this signature is not computed by applying persistence
to some function defined on the multinerve, but it is rather
a pruned version of the extended persistence diagram of ${\tilde f}$. 
As for Reeb graphs,
it serves as a bag-of-features type signature of the structure of
$\MMapper_f(X, \I)$. 
Moreover, the fact that $\Dg(\MMapper_f(X, \I)) \subseteq \Dg(\tilde
f)$ formalizes the intuition that the MultiNerve Mapper should be
viewed as a {\em pixelized version} of the Reeb graph, in which some
of the features disappear due to the staircases (prescribed by the cover). 
For instance, in Figure~\ref{fig:ExDict} we show a double torus
equipped with the height function, together with its associated Reeb
graph, MultiNerve Mapper, and Mapper. We also show the corresponding
extended persistence diagrams.  In each case, the points in the diagram
represent the features of the object: the extended points represent
the holes (dimension 1 and above) and the trunks (dimension 0) while
the ordinary and relative points represent the branches.

\paragraph{Convergence of the signature.} The following convergence
result (which is in fact non-asymptotic) is a direct consequence of
our previous results:
\begin{corollary}\label{cor:Reeb_approx}
Suppose the granularity of the gomic~$\I$ is at most $\e$. Then, 
\[
 \Dg(\tilde f)\setminus \{(x,y) \,:\, |y-x|\leq \e\} \subseteq \Dg(\MMapper_f(X, \I)) \subseteq \Dg(\tilde f).
\]
Thus, the features (branches, holes) of the Reeb graph that are
missing in the MultiNerve Mapper have spans at most~$\e$.  In
particular, we have $\distb(\Dg(\MMapper_f(X, \I)),
\Dg(\tilde f))\leq \e/2$.  Moreover, the two signatures become equal
when $\e$ becomes smaller than the smallest vertical distance of the
points of $\Dg(\tilde f)$ to the diagonal. Finally, $\MMapper_f(X,
\I)$ and $\Reeb_f(X)$ themselves become isomorphic as combinatorial graphs up to
one-step vertex splits and edge subdivisions (which are
topologically trivial modifications)
when $\e$ becomes smaller than the smallest absolute difference
between distinct critical values of~$f$.
\end{corollary}
We show a similar convergence
result in the functional distortion distance in Section~\ref{sec:stabdfd}.
Note that building the signature $\Dg(\MMapper_f(X, \I))$
requires computing the critical values of~$f$ exactly, which may not
always be possible. However, as for Reeb graphs, the signature can be
approximated efficiently and with theoretical guarantees under mild
sampling conditions using existing work on scalar fields analysis, as
we will see in Section~\ref{sec:discrete}.

\subsection{Induced signature for Mapper}
\label{sec:MapperSign} 

Recall from Lemma~\ref{lem:pi1_surj} that the projection
$\pi_1:\MMapper_f(X,\I)\to \Mapper_f(X,\I)$ induces a surjective
homomorphism in homology. Thus, the Mapper has a simpler structure
than the MultiNerve Mapper. To be more specific, $\pi_1$ identifies all the 
edges connecting the same pair of vertices. This eliminates the corresponding 
holes in $\MMapper_f(X,\I)$. Since the two vertices
lie in successive intervals of the cover, the corresponding diagram points lie 
in the following extended staircase 
(see the staircase $Q^{\I}_E$ displayed on the right in Figure~\ref{fig:staircase}):
\[
\StairExtM=\bigcup_{I\cup J \text{ such that } I\cap J\neq \emptyset}Q_{I\cup J}^-.
\]
The other staircases remain unchanged. Hence the following signature:
\begin{align}\label{eq:sign_Mapper}
\begin{split}
\Dg(\Mapper_f(X,\I))&=
(\Ord(\tilde{f})\setminus Q_O^{\I})\cup
(\Ext(\tilde{f})\setminus Q_{E}^{\I})\cup
(\Rel(\tilde{f})\setminus Q_R^{\I})\\
&=
(\Ord_0(f)\setminus Q_O^{\I})\cup
((\Ext_0^+(f)\cup \Ext_1^-(f))\setminus Q_{E}^{\I})\cup
(\Rel_1(f)\setminus Q_R^{\I}).
\end{split}
\end{align}
The interpretation of this signature in terms of the structure of the
Mapper follows the same rules as for the MultiNerve Mapper and Reeb
graph---see again Figure~\ref{fig:ExDict}.
Moreover, the convergence result stated in
Corollary~\ref{cor:Reeb_approx} holds for the Mapper as well.

\section{Stability in the bottleneck distance}
\label{sec:Stab_func}

Intuitively, for a point in the signature $\Dg(\MMapper_f(X,\I))$,
the $\ell^\infty$-distance to its corresponding
staircase\footnote{$Q_O^{\I}$, $Q_{E^-}^{\I}$ or $Q_R^{\I}$, depending on the type of the point.}  
measures the amount by which the
function~$f$ or the cover~$\I$ must be perturbed in order to eliminate
the corresponding feature (branch, hole) in the MultiNerve
Mapper. Conversely, for a point in the Reeb graph's signature
$\Dg(\tilde f)$ that is not in the MultiNerve Mapper's signature
(i.e. that lies inside its corresponding staircase), the
$\ell^\infty$-distance to the boundary of the staircase measures the
amount by which~$f$ or~$\I$ must be perturbed in order to create a
corresponding feature in the MultiNerve Mapper. Our goal here is to
formalize this intuition. For this we adapt the bottleneck distance so
that it takes the staircases into account. Our results are stated for
the MultiNerve Mapper, they hold the same for the Mapper with the
staircase $Q_{E^-}^{\I}$ replaced by its extension $Q_{E}^{\I}$.

\paragraph{An extension of the bottleneck distance.} Let $\Theta$ be a subset of
$\R^2$. Given a partial matching $\Gamma$ between two extended persistence
diagrams $\Dg, \Dg'$, the \emph{$\Theta$-cost} of $\Gamma$ is:
\[
\mathrm{cost}_\Theta(\Gamma)=\max\left\{\max_{p\in \Dg}\ \delta_{\Dg}(p),\ \max_{p'\in {\Dg'}}\ \delta_{\Dg'}(p')\right\},
\]
where:
\[
\delta_{\Dg}(p)=\|p-p'\|_\infty \text{ if } \exists p'\in \Dg'\text{ such that }(p,p')\in\Gamma \text{ and } d_\infty(p,\Theta) \text{ otherwise},
\]
\[
\delta_{\Dg'}(p')=\|p-p'\|_\infty \text{ if } \exists p\in \Dg\text{ such that }(p,p')\in\Gamma \text{ and } d_\infty(p',\Theta) \text{ otherwise}.
\]
The bottleneck distance becomes:
\[
\diststair(\Dg,\Dg')=\inf_{\Gamma}\ \mathrm{cost}_\Theta(\Gamma),
\]
where $\Gamma$ ranges over all partial matchings between $\Dg$ and
$\Dg'$. This is again a pseudometric and not a metric. 
Note that the usual bottleneck distance
is obtained by taking $\Theta$ to be the diagonal $\Diag$.  
Given a gomic $\I$, we choose different sets $\Theta$ depending on the
types of the points in the two diagrams. More precisely, we define the
distance between signatures as follows:

\begin{definition}
Given a gomic $\I$, we define the distance $d_{\mathcal I}$ between extended persistence diagrams $\Dg,\Dg'$ as:
\begin{equation}\label{eq:def-bottleneck-MMapper} 
\distcov(\Dg,\Dg')=\max\left\{
d_{\rm{b},Q_O^{\I}}(\Ord, \Ord'),\; 
d_{\rm{b},Q_{E^-}^{\I}}(\Ext, \Ext'),\; 
d_{\rm{b},Q_R^{\I}}(\Rel, \Rel') \right\}. 
\end{equation}
\end{definition}

\subsection{Stability with respect to perturbations of the function}

The distance $d_{\I}$ {\em stabilizes} the (MultiNerve) Mappers, as stated in the following theorem:

\begin{theorem}\label{thm:MStab}
Given a topological space $X$, Morse-type functions $f,g:X\to\R$ and a gomic $\I$ of granularity at most $\epsilon > 0$, 
the following stability inequality holds:
\begin{align}
&d_{\I}(\Dg(\Mapper_f(X,\I)),\Dg(\Mapper_g(X,\I))) \leq d_{\I}(\Dg(\MMapper_f(X,\I)),\Dg(\MMapper_g(X,\I))) \leq \|f-g\|_\infty. \label{ineq:stabilitydi}
\end{align}

Moreover, $d_{\I}$ and $\distb$ are related as follows:
\begin{align}
&\distb(\Dg(\MMapper_f(X,\I)),\Dg(\MMapper_g(X,\I))) \leq \frac \e2 + d_{\I}(\Dg(\MMapper_f(X,\I)),\Dg(\MMapper_g(X,\I))).\label{ineq:didbMM}\\
&\distb(\Dg(\Mapper_f(X,\I)),\Dg(\Mapper_g(X,\I))) \leq \e + d_{\I}(\Dg(\Mapper_f(X,\I)),\Dg(\Mapper_g(X,\I))).\label{ineq:didbM}
\end{align}
\end{theorem}
\begin{figure}[h]\centering
\includegraphics[height=4cm]{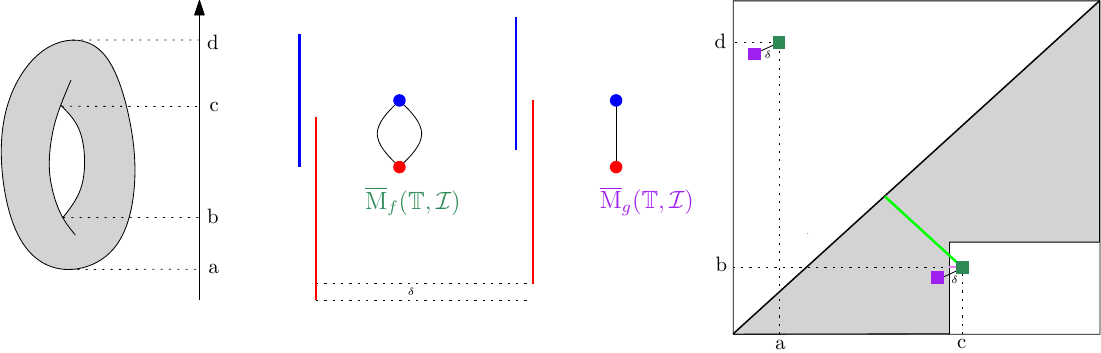}
\caption[Stability of the Mapper]{We compute the MultiNerve Mapper of the height function $f$ on the torus $\T$,
given a gomic $\I$ with two intervals.
We also compute the MultiNerve Mapper of a perturbed function $g$ such that $\|f-g\|_\infty\leq\delta$.
We plot the extended persistence diagrams of $\tilde f$ (dark green) and $\tilde g$ (purple). Note that the signature of $\MMapper_g(\T,\I)$
is obtained by removing the purple point beneath the diagonal since it belongs to a staircase, while the signature of $\MMapper_f(\T,\I)$
is equal to $\Dg({\tilde f})$. If we used the bottleneck distance to compare the two signatures, their distance would be equal to
the distance to the diagonal of the dark green point beneath $\Diag$ (green segment), which can be arbitrarily large, while, using 
$d_{\I}$, their distance becomes the distance of the
same point to the staircase (tiny pink segment), which is bounded by $\delta$.}
\label{fig:counterexample} 
\end{figure}

The proof of Theorem~\ref{thm:MStab} relies on the following monotonicity
property, which is immediate:
\begin{lemma}\label{lem:ineqsub}
Let $\Theta\subseteq\R^2$ be in the closure of $\Theta'\subseteq\R^2$. 
Then, $$d_{\Theta'}(\Dg,\Dg')\leq d_\Theta(\Dg,\Dg')\leq d_{\Theta'}(\Dg,\Dg')+\disth(\Theta,\Theta'),$$
where $\disth$ denotes the Hausdorff distance in the $\ell_\infty$-norm.
\end{lemma}

\begin{proof}[Theorem~\ref{thm:MStab}]
Equation~(\ref{ineq:didbMM}) and~(\ref{ineq:didbM}) are direct applications of Lemma~\ref{lem:ineqsub}.
Equation~(\ref{ineq:stabilitydi}) is proven by the following sequence of (in)equalities:
\begin{align*}
d_{\I}(\Dg(\MMapper_f(X,\I)),\Dg(\MMapper_g(X,\I))) &=d_{\I}(\Dg(\tilde f), \Dg(\tilde g))\\
&\leq d_{\rm{b},\Diag}(\Dg(\tilde f), \Dg(\tilde g))=\distb(\Dg(\tilde f), \Dg(\tilde g))\\
&\leq \distb(\Dg(f),\Dg(g))\\
&\leq \|f-g\|_\infty.
\end{align*}
The first equality comes from the observation that the points
of $\Dg(\tilde f)\sqcup\Dg(\tilde g)$ that lie inside
their corresponding staircase can be left unmatched and have a zero
cost in the matching, so removing them as in~(\ref{eq:sign1}) does not change the bottleneck cost. 
The first inequality follows from Lemma~\ref{lem:ineqsub} since the diagonal $\Diag$ is included in the closure of each of the staircases. 
The second inequality follows from Theorem~\ref{thm:pdreeb} and the fact that the matchings only match points of 
the same type (ordinary, extended, relative) and of the same homological dimension. The last inequality comes from  Theorem~\ref{th:Stab}.
\qed\end{proof}

\paragraph{Interpretation of the stability.} Note that the bottleneck distance $\distb$ is unstable in this context---see Figure~\ref{fig:counterexample}. 
The theorem allows us to make some interesting claims. 
For instance, denoting by $Q_p^{\I}$ the staircase corresponding to the type of a diagram point~$p$,  the quantity 
\[
d_{\I}(\Dg, \emptyset) = \max_{p\in \Dg} d_\infty(p, Q_p^{\I}) 
\]
measures the amount by which the diagram~$\Dg$ must be perturbed in the
metric~$d_{\I}$ in order to bring all its points to the staircase. Hence, by Theorem~\ref{thm:MStab}, given a pair~$(X,f)$, the quantity 
\[
d_{\I}(\Dg(\MMapper_f(X,\I)), \emptyset) = \max_{p\in \Dg(\MMapper_f(X,\I))} d_\infty(p, Q_p^{\I}) 
\]
is a lower bound on the amount by which $f$ must be perturbed in
the supremum norm in order to remove all the features (branches and
holes) from the MultiNerve Mapper. Conversely, 
\[
\min_{p\in \Dg(\MMapper_f(X,\I))} d_\infty(p, Q_p^{\I}) 
\]
is a lower bound on the maximum amount of perturbation allowed for~$f$
if one wants to preserve all the features in the MultiNerve
Mapper no matter what. Note that this does not prevent other features from
appearing. The quantity that controls those is related to the
points of~$\Dg(\tilde f)$ (including 
diagonal points) that lie in the staircases. More precisely, the quantity
\[
\min_{p\in \Dg(\tilde f) \cup \Delta} d_\infty(p, \partial Q_p^{\I}\setminus\Delta)
\]

is a lower bound on the maximum amount by which $f$ can be perturbed
if one wants to preserve the structure (set of features) of the
MultiNerve Mapper no matter what. Note that this lower bound is
in fact zero since $\partial Q_O^{\I}\setminus\Diag$ and $\partial
Q_R^{\I}\setminus\Diag$ come arbitrarily close to the diagonal~$\Delta$
(recall Figure~\ref{fig:staircase}). This means that, as small as the
perturbation of~$f$ may be, it can always make new branches appear in
the MultiNerve Mapper. However, it will not impact the set of holes if
its amplitude is less than
\[
\min_{p\in \Ext(\tilde f) \cup \Delta} d_\infty(p, \partial Q_{E^-}^{\I}\setminus\Delta).
\]
From this discussion we derive the following rule of thumb: having
small overlaps between the intervals of the gomic helps capture more
features (branches and holes) of the Reeb graph in the (MultiNerve)
Mapper; conversely, having large overlaps helps prevent new holes from
appearing in the (MultiNerve) Mapper under small perturbations of the
function. This is an important trade-off to consider in applications.

\subsection{Stability with respect to perturbations of the domain}

More generally, we can derive a stability result for perturbations of
the pair $(X,f)$, provided we make some extra assumptions on the
regularity of the domain and function. Typically, we will assume $X$
to be a compact Riemannian manifold (or, more generally, a compact
length space with curvature bounded above) and $f$ to be
Lipschitz-continuous. To measure the amount of perturbation of the
domain we use the concept of {\em correspondence} from metric
geometry: given another pair $(Y,g)$, a correspondence is a subset~$C$
of the product space $X\times Y$ such that the canonical projections
$C\to X$ and $C\to Y$ are surjective.  We consider the {\em functional
  distortion} associated with~$C$, which is the quantity:
\[
\e_{\mathfrak{f}}(C) = \sup_{(x,y)\in C} |f(x)-g(y)|.
\]
Similarly, writing respectively $d_X$ and $d_Y$ for the intrinsic
metrics of $X$ and $Y$, we consider the {\em metric distortion} of~$C$:
\[
\e_{\mathfrak{m}}(C) = \sup_{(x,y)\in C, (x',y')\in C} |d_X(x,x')-d_Y(y,y')|.
\]
The {\em Gromov-Hausdorff distance} between $X$ and $Y$ is then:
\[
d_{\text{GH}}(X,Y) = \frac{1}{2}\inf_C \e_{\mathfrak{m}}(C),
\]
where $C$ ranges over all correspondences between $X$ and $Y$. Now we
can derive a stability guarantee for the signatures of MultiNerve
Mappers in this context, using a variant of Theorem~\ref{th:Stab}
proven in~\cite{Carriere15b}:
\begin{theorem}\label{thm:DStab}
Fix a gomic~$\I$. Let $X$ and $Y$ be two compact Riemannian manifolds
or length spaces with curvature bounded above.  Denote by $\rho(X)$
and $\rho(Y)$ their respective convexity radii (i.e. the smallest radius for which any geodesic ball is convex).  Let
$f:X\rightarrow\R$ and $g:Y\to\R$ be Lipschitz-continuous Morse-type
functions, with Lipschitz constants $c_f$ and $c_g$ respectively.
Assume $d_\emph{GH}(X,Y)\leq\frac{1}{20}\,\min(\rho(X),\rho(Y))$.
Then, for any correspondence $C\in\mathcal{C}(X,Y)$ such that
$\epsilon_{\mathfrak{m}}(C)<\frac{1}{10}\,\min(\rho(X),\rho(Y))$,
\[
d_\I(\Dg(\MMapper_f(X,\I)),\Dg(\MMapper_g(Y,\I)))\leq (9(c_f+c_g)+\min\{c_f,c_g\})\epsilon_{\mathfrak{m}}(C)+\epsilon_{\mathfrak{f}}(C).
\]
\end{theorem}
\begin{proof}
The proof is the same sequence of (in)equalities as for
Theorem~\ref{thm:MStab}, except the last inequality is replaced by
$d_\Delta(\Dg(f),\Dg(g))\leq
(9(c_f+c_g)+\min\{c_f,c_g\})\epsilon_{\mathfrak{m}}(C)+\epsilon_{\mathfrak{f}}(C)$, which comes\footnote{Note
  that Theorem~3.4 in~\cite{Carriere15b} is stated only for the ordinary part of
  the persistence diagrams, however its analysis extends to the
  full extended filtrations at no extra cost.} from Theorem~3.4
in~\cite{Carriere15b}.
\qed\end{proof}

This result brings about the same discussion as in
Section~\ref{sec:Stab_func}, with $f$ replaced by the pair~$(X,f)$.

\section{Stability with respect to perturbations of the cover}
\label{sec:Stab_cover}

Let us now fix the pair $(X,f)$ and consider varying gomics.  For each
choice of gomic, 
Eqs.~(\ref{eq:sign1})-(\ref{eq:sign_Mapper}) tell which points
of the diagram~$\Dg(f)$ end up in the diagram of the (MultiNerve)
Mapper and thus participate in its structure. We aim for a
quantification of the extent to which this structure may change as the
gomic is perturbed.  For this we adopt the dual point of view: for any
two choices of gomics, we want to use the points of the
diagram~$\Dg(f)$ to assess the degree by which the gomics differ.
This is a reversed situation compared to Section~\ref{sec:Stab_func},
where the gomic was fixed and was used to assess the degree by which
the persistence diagrams of two functions differed.

\paragraph{A distance between gomics.} The diagram points that discriminate between the two gomics are the
ones located in the symmetric difference of the staircases, since they
witness that the symmetric difference is
non-empty. Moreover, their $\ell^\infty$-distances to the staircase of
the other gomic provide a lower bound on the Hausdorff distance
between the two staircases and thus quantify the extent to which the
two covers differ. We formalize this intuition as follows: given a
persistence diagram~$\Dg$ and two gomics $\I, \J$, we consider the
quantity:
\begin{equation}\label{eq:def-dist-covers}
d_{\Dg}(\I, \J) = \max_{*\in \{O, {E^-}, R\}}
\left\{ \sup_{p\in \Dg^*\cap (Q_*^{\I} \symdiff Q_*^{\J})} \max\left\{ d_\infty(p, Q_*^{\I}),\, d_\infty(p, Q_*^{\J})\right\}
\right\},
\end{equation}
where $\symdiff$ denotes the symmetric difference, where $\Dg^*$ stands
for the subdiagram of $\Dg$ of the right type ($\Ord$, $\Ext$ or
$\Rel$), and where we adopt the convention that $\sup_{p\in
  \emptyset} ... $ is zero instead of infinite.  Note that there is
always one of the two terms in (\ref{eq:def-dist-covers}) that is zero
since the supremum is taken over all points that lie in the symmetric
difference of the staircases.  Deriving an upper bound on $d_{\Dg}(\I,\J)$ 
in terms of the Hausdorff distances between the staircases is
straightforward, since the supremum in~(\ref{eq:def-dist-covers}) is
taken over points that lie in the symmetric difference between the
staircases:

\[
d_{\Dg}(\I, \J) \leq \max_{*\in \{O, {E^-}, R\}} \disth(Q_*^{\I},\, Q_*^{\J}),
\]
where $\disth$ stands for the Hausdorff distance in the $\ell^\infty$-norm.
%
The connection to the MultiNerve Mapper appears when we take $\Dg$ to be the
persistence diagram of the induced map~$\tilde f$ defined on the Reeb
graph $\Reeb_f(X)$. Indeed, we have
\[
\Ord(\tilde f)\cap (Q_O^{\I} \symdiff Q_O^{\J}) = (\Ord(\tilde f)\cap Q_O^{\I}) \symdiff
(\Ord(\tilde f)\cap Q_O^{\J}) = \Ord(\MMapper_f(X,\I)) \symdiff
\Ord(\MMapper_f(X,\J)),
\]
where the second equality follows from the definition of the
signature of the MultiNerve Mapper given in~(\ref{eq:sign1}). 
Similar equalities can be derived with $\Ext$ and $\Rel$. Thus,
$d_{\Dg(\tilde f)}(\I, \J)$ quantifies the proximity of each
signature to the other staircase.  In particular, having
$d_{\Dg(\tilde f)}(\I, \J) = 0$ means that there are no diagram points
in the symmetric difference, so the two gomics are equivalent from the
viewpoint of the structure of the MultiNerve Mapper. Differently,
having $d_{\Dg(\tilde f)}(\I, \J) > 0$ means that the structures of
the two MultiNerve Mappers differ, and the value of $d_{\Dg(\tilde
  f)}(\I, \J)$ quantifies by how much the covers should be perturbed
to make the two multigraphs isomorphic. Furthermore,
we have the following upper bound on this quantity:
\begin{theorem}\label{thm:CStab} 
Given a Morse-type function $f:X\to\R$, for any gomics $\I, \J$,
\[
d_{\Dg(\tilde f)}(\I, \J) \leq  \max_{*\in \{O, {E^-}, R\}} \disth(Q_*^{\I},\, Q_*^{\J}).
\]
%
\end{theorem}
\paragraph{Tightness.} It is easy to build examples where the upper bound is tight, for
instance by placing a diagram point at a corner of one of the
staircases\footnote{Which is easily done by choosing suitable critical
  values as coordinates for this point.}. On the other hand, there are
obvious cases where the bound is not tight, for instance we have
$d_{\Dg(\tilde f)}(\I, \J) = 0$ as soon as there are no diagram points
in the symmetric difference, whereas the symmetric difference itself
may not be empty.  What the upper bound measures depends on the
subdiagram. For instance, for $*=E^-$, we defined $Q_{E^-}^{\I}$ to be
the set $\bigcup_{(a,b)\in\I} \{(x,y)\in\R^2\,:\, a\leq y< x\leq b\}$,
so $\disth(Q_{E^-}^{\I},\, Q_{E^-}^{\J})$ measures the supremum of
the differences between the intervals in one cover to their closest
interval in the other cover:
\[
\disth(Q_{E^-}^{\I},\, Q_{E^-}^{\J}) = \max\left\{ \sup_{(a,b)\in\I} \inf_{(c,d)\in\J} 
\max\{|a-c|,\, |b-d|\},\; \sup_{(c,d)\in\J} \inf_{(a,b)\in\I} \max\{|a-c|,\, |b-d|\}\right\}.
\]
Similar formulas can be derived for the other subdiagrams.

\section{Convergence in the functional distortion distance}
\label{sec:stabdfd}

Since $\distb$ is merely a pseudometric, the relationship between the
(MultiNerve) Mapper and the Reeb graph is only partially explained by
Theorem~\ref{th:ExDgMultiNerve}.  In this section, we bound the {\em
  functional distortion distance} $\distfd$ (a true distance between
metric graphs equipped with continuous functions) between the
(MultiNerve) Mapper and the Reeb graph, and we provide an alternative
proof of Theorem~\ref{th:ExDgMultiNerve} as a byproduct.  To this end,
we connect the (MultiNerve) Mapper and the Reeb graph through a
sequence of metric spaces on which we can control the functional
distortion distance.  This connection has an interest in its own
right, as it was leveraged in other contributions on Mappers and Reeb
graphs recently---see e.g.~\cite{Carriere17c,Carriere17a}.

\subsection{Telescopes and Operators}\label{sec:telescope}

In this section we introduce the telescopes, which are our main
objects of study when we relate the MultiNerve Mapper to
the Reeb graph.


Recall that, given topological spaces $X$ and $A\subseteq Y$ together
with a continuous map $f:A\rightarrow X$, the {\em adjunction space}
$X\cup_f Y$ (also denoted $Y\cup_f X$) is the quotient of the disjoint
union $X\amalg Y$ by the equivalence relation induced by the
identifications $\{f(a)\sim a\}_{a\in A}$.


\begin{definition}[Telescope~\cite{Carlsson09b}]\label{def:telescope}
A {\em telescope} is an adjunction space of the following form:
$$T=\left(Y_0\times (a_0,a_1]\right) \cup_{\psi_0} \left(X_1\times \{a_1\}\right) \cup_{\phi_1} \left(Y_1
\times [a_1,a_2]\right) \cup_{\psi_1} ...\ \cup_{\phi_n} \left(Y_n\times [a_n,a_{n+1})\right),$$
where $-\infty=a_0<a_1<\cdots< a_n<a_{n+1}=+\infty$,
and where the $\phi_i:Y_i\times\{a_i\}\rightarrow X_i\times\{a_i\}$
and $\psi_i:Y_i\times\{a_{i+1}\}\rightarrow X_{i+1}\times\{a_{i+1}\}$ are continuous maps. 
The $a_i$ are called the {\em critical values} of $T$ and their set is denoted by $\Crit(T)$,
the $\phi_i$ and $\psi_i$ are called {\em attaching maps}, 
the $Y_i$ are compact and locally connected spaces called the {\em cylinders}
and the $X_i$ are topological spaces called the {\em critical slices}.
Moreover, all $Y_i$ and $X_i$ have finitely-generated homology.
\end{definition} 

\paragraph{Extended persistence diagram.} A telescope comes equipped with functions $\pi_1$ and $\pi_2$,
which are the projections onto the first factor and second factor respectively.
From now on, given any interval $I$, we let $T^I$ denote $\pi_1\circ\pi_2^{-1}(I)$.
Then, the extended persistence diagram $\Dg(\pi_2)$ can be described using the following Lemma.

\begin{lemma}\label{lem:defret} 
Since $\phi_i$ and $\psi_i$ are continuous,
\begin{align*}
\forall \alpha\in [a_i,a_{i+1}),\ T^{(-\infty,\alpha]}&\text{ deform retracts onto }T^{(-\infty,a_i]} \\
\forall \alpha\in (a_{i-1},a_i],\ T^{[\alpha,+\infty)}&\text{ deform retracts onto }T^{[a_i,+\infty)},
\end{align*}
where a topological space $X$ is said to {\em deform retract} onto $Y\subseteq X$ if there exists a continuous function $F:X\times[0,1]\rightarrow X$
such that $F(\cdot,0)={\rm id}_X$, $F|_{Y\times\{\alpha\}}(\cdot,\alpha)={\rm id}_Y$ for any $\alpha\in[0,1]$, and $F(X,1)\subseteq Y$. In particular, 
this means that the inclusion $Y\hookrightarrow X$ is a homotopy equivalence. 
\end{lemma}


\begin{corollary} The following inclusion holds:
$\Dg(\pi_2)\subseteq\Crit(T)\times\Crit(T)$.
\end{corollary}


\paragraph{Construction from a Morse-type function.}
One can build telescopes from the domain of Morse-type functions---see Definition~\ref{def:Morse-type}.
Indeed, a function $f:X\rightarrow\R$ of Morse type 
naturally induces a telescope $T(X,f)$ with
\begin{itemize} 
\item $\Crit(T(X,f))=\Crit(f)$, 
\item $X_i=f^{-1}(a_i)$, 
\item $Y_i=\pi_1\circ\mu_i^{-1}\circ f^{-1}((a_i,a_{i+1}))$,
\item $\phi_i:(y,a_i)\mapsto(\bar{\mu}_i|_{Y_i\times\{a_i\}}(y,a_i),a_i)$, $\forall y\in Y_i$, $\forall i\in\{1,...,n\}$,
\item $\psi_i:(y,a_{i+1})\mapsto(\bar{\mu}_i|_{Y_i\times\{a_{i+1}\}}(y,a_{i+1}),a_{i+1})$, $\forall y\in Y_i$, $\forall i\in\{0,...,n-1\}$,
\end{itemize}

$T(X,f)$ is well-defined thanks to the following Lemma:

\begin{lemma}
\label{lem:restrictions}
$\im(\phi_i)\subseteq f^{-1}(a_i)\times\{a_i\}$ and 
$\im(\psi_i)\subseteq f^{-1}(a_{i+1})\times\{a_{i+1}\}$.
\end{lemma}

\begin{proof}
Let $(y,a_{i+1})\in Y_i\times\{a_{i+1}\}$.  Consider the sequence
$(y,v_n)_{n\in\mathbb{N}}$, for an arbitrary
$(v_n)_{n\in\mathbb{N}}\in(a_i,a_{i+1})^{\mathbb{N}}$ that converges
to $a_{i+1}$.  Then, $(f\circ\bar{\mu}_i(y,v_n))_{n\in\mathbb{N}}$
converges to $f\circ\bar{\mu}_i(y,a_{i+1})$ by continuity of
$f\circ\bar{\mu}$.  Moreover, for all $n\in\N$ we have
$f\circ\bar{\mu}_i(y,v_n)=f\circ\mu_i(y,v_n)=v_n$ since
$f|_{f^{-1}(a_i,a_{i+1})}=\pi_2\circ\mu_i^{-1}$.  Therefore,
$(f\circ\bar{\mu}_i(y,v_n))_{n\in\mathbb{N}}$ converges also to
$a_{i+1}$. By uniqueness of the limit, we have
$f\circ\bar{\mu}_i(y,a_{i+1})=a_{i+1}$, meaning that
$\bar{\mu}_i(y,a_{i+1})\in f^{-1}(a_{i+1})$. Thus,
${\rm im}(\psi_i)\subseteq f^{-1}(a_{i+1})\times\{a_{i+1}\}$.  
The same argument applies to show that 
${\rm im}(\phi_i)\subseteq f^{-1}(a_i)\times\{a_i\}$.
\qed\end{proof}

\paragraph{Correspondence between $X$ and $T(X,f)$.} We now exhibit a homeomorphism between $T(X,f)$ and $X$. 
Let $\mu:T(X,f)\rightarrow X$ be defined by:
\[
 \mu(y,z)= \left \{ \begin{array}{l} 
y\text{ if } (y,z)\in X_i\times\{a_i\} \mbox{ for some $i$};\\
\mu_i(y,z) \text{ if } (y,z)\in Y_i\times (a_i,a_{i+1}) \mbox{ for some $i$}.
\end{array} \right.
\]
The map $\mu$ is bijective as every $\mu_i$ is. It is also continuous as every $\bar{\mu}_i$ is.
Since every continuous bijection from a compact space to a Hausdorff space is a homeomorphism 
(see e.g. Proposition 13.26 in \cite{Sutherland09}),
$\mu$ defines a homeomorphism between $T(X,f)$ and $X$. Moreover, $\pi_2=f\circ\mu$ so $\Dg(f)=\Dg(\pi_2)$. 

\subsubsection*{Operators on telescopes}
\label{sec:operator}

The decomposition of telescopes into cylinders can be used to define simple operators that 
modify the telescope structures in a predictable way. Specifically, we detail three
types of operators, corresponding to the cases where one asks for either removal of critical
values ($\Merge$ operator), 
duplication of critical values ($\Split$ operator), 
or translation of critical values ($\Shift$ operator). 
To formalize this, we use \textit{generalized attaching maps}:

\begin{center}
$\begin{array}{llll}
\phi_i^a: & Y_i\times\{a\}\rightarrow & X_i\times\{a\}; & (y,a) \mapsto(\pi_1\circ\phi_i(y,a_i),a), \\
\psi_i^a: & Y_i\times\{a\}\rightarrow & X_{i+1}\times\{a\}; & (y,a) \mapsto(\pi_1\circ\psi_i(y,a_{i+1}),a).
\end{array}$
\end{center}


\paragraph{Merge.} Merge operators merge all critival values located in $[a,b]$ into
a single critical value $\bar a=\frac{a+b}{2}$.  
\begin{definition}[Merge]
Let $T$ be a telescope.
Let $a\leq b$.  If $[a,b]$ contains at least one critical value,
i.e. $\exists i,j\in\N$ such that $a_{i-1}<a\leq a_i\leq a_j\leq b < a_{j+1}$, then the {\em Merge} on
$T$ between $a,b$ is the telescope $T'=\Merge_{a,b}(T)$ given by:
\[\arraycolsep=1.4pt\def\arraystretch{1.2}
\begin{array}{c}
...(Y_{i-1}\times[a_{i-1},a_i])\cup_{\psi_{i-1}}(X_i\times\{a_i\})\cup_{\phi_i}...\cup_{\psi_{j-1}}(X_j\times\{a_j\})\cup_{\phi_j}(Y_j\times[a_j,a_{j+1}])...\\
\rotatebox[origin=c]{270}{$\mapsto$}\\
...(Y_{i-1}\times[a_{i-1},\bar{a}])\cup_{f_{i-1}}(T^{[a,b]}\times\{\bar{a}\})\cup_{g_j}(Y_j\times[\bar{a},a_{j+1}])...
\end{array}
\]
where $\bar{a}=\frac{a+b}{2}$, where
$f_{i-1}=\psi_{i-1}^{\bar{a}}$ if $a=a_i$ and
$f_{i-1}={\rm id}_{Y_{i-1}\times\{\bar{a}\}}$
otherwise, and where $g_j=\phi_j^{\bar{a}}$ if $b=a_j$ and
$g_{j}={\rm id}_{Y_j\times\{\bar{a}\}}$
otherwise.

If $[a,b]$ contains no critical value, i.e. $a_{i-1}<a\leq
b<a_i$, then $\Merge_{a,b}(T)$ is given by:
\[\hspace{-30pt}\arraycolsep=1.4pt\def\arraystretch{1.2}
\begin{array}{c}
...(X_{i-1}\times\{a_{i-1}\})\cup_{\phi_{i-1}}(Y_{i-1}\times[a_{i-1},a_i])\cup_{\psi_{i-1}}(X_i\times\{a_i\})...\\
\rotatebox[origin=c]{270}{$\mapsto$}\\
...
\cup_{\phi_{i-1}}
(Y_{i-1}\times[a_{i-1},\bar{a}])\cup_{f_{i-1}}(T^{[a,b]}\times\{\bar{a}\})\cup_{g_{i-1}}(Y_{i-1}\times[\bar{a},a_i])
\cup_{\psi_{i-1}}
...
\end{array}
\]
where $\bar{a}=\frac{a+b}{2}$, and where
$f_{i-1}=g_{i-1}={\rm id}_{Y_{i-1}\times\{\bar{a}\}}$.
\end{definition}
\begin{figure}[h]\centering
\includegraphics[height=4cm]{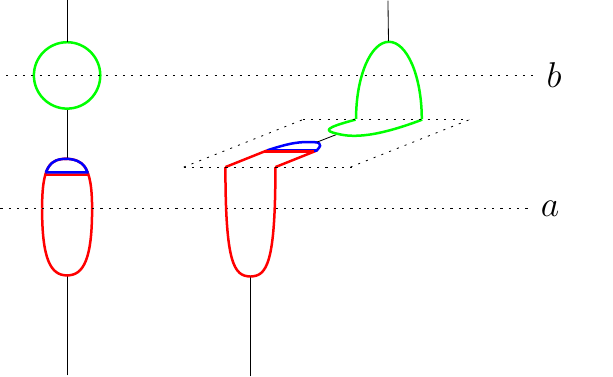}\ \ \ \includegraphics[width = 4cm]{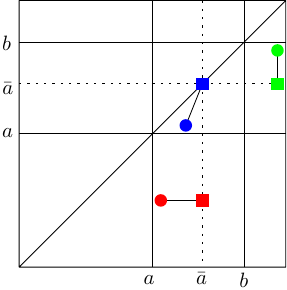}
\caption[Merge]{Left: Effect of a Merge on a telescope. 
Right: Effect on the corresponding extended persistence diagram. 
Points before the Merge are disks while points after the Merge are squares.}
\label{fig:mergespace} 
\end{figure}
See the left panel of Figure~\ref{fig:mergespace} for an illustration.
\paragraph{Merge for persistence diagrams.} 
Similarly, we define the Merge between $a,b$ on an extended persistence diagram~$\Dg$ as the diagram $\Merge_{a,b}(\Dg)$ given by
$\Merge_{a,b}(x,y)=({\bar x},{\bar y})$, where:
\[
{\bar x}=\left\{ \begin{array}{ll} x\text{ if }x\notin[a,b] \\ {\bar a}\text{ otherwise } \end{array} \right.
\text{ and }
{\bar y}=\left\{ \begin{array}{ll} y\text{ if }y\notin[a,b] \\ {\bar a}\text{ otherwise } \end{array} \right.
\]
Points in the strips $x\in[a,b]$, $y\in[a,b]$ are snapped to the lines
$x=\bar{a}$ and $y=\bar{a}$ respectively.  See the right panel of
Figure~\ref{fig:mergespace}.  See also the first intermediate points
along the trajectories of the red points in
Figure~\ref{fig:pdTRANS} for another illustration on extended
persistence diagrams.  

\paragraph{Commutativity of the operators.} We now prove 
that extended persistent homology commutes with this operator, i.e. $\Dg(\Merge)=\Merge(\Dg)$.

\begin{lemma}
\label{lem:mergeprop}
Let $a\leq b$ and $T'=\Merge_{a,b}(T)$. 
Let $\pi_2':T'\rightarrow\R$ be the projection onto the second factor. Then, $\Dg(\pi_2')=\Merge_{a,b}(\Dg(\pi_2))$.
\end{lemma}

\begin{proof}
We only study the sublevel sets of the functions, which means that we
only prove the result for the ordinary part of the diagrams.  The
proof is symmetric for superlevel sets, leading to the result for
the extended and the relative parts. 

Assume $a_{i-1}<a\leq a_i\leq a_j\leq b<a_{j+1}$.  Given $x\leq y$, we
let $\Pi_{x,y}:H_*(T^{(-\infty,x]})\rightarrow
      H_*(T^{(-\infty,y]})$ and
    $\Pi'_{x,y}:H_*((T')^{(-\infty,x]})\rightarrow
  H_*((T')^{(-\infty,y]})$ be the homomorphisms induced by
inclusions.  Since $f$ is of Morse type, Lemma~\ref{lem:defret}
relates $\Pi'$ to $\Pi$ as follows (see Figure~\ref{fig:areas}):

\begin{equation}\label{eq:rel}
\Pi'_{x,y}=\left\{ \begin{array}{ll} \Pi_{x,y}\text{ if }x,y\notin[a,b]\text{ (green)} & \Pi_{x,a_{i-1}}\text{ if }x<a,y\in[a,\bar{a})\text{ (pink)} \\
					   \Pi_{a_{i-1},y}\text{ if }x\in[a,\bar{a}),y>b\text{ (blue)} & \Pi_{a_{i-1},a_j}\text{ if }x\in[a,\bar{a}),y\in[\bar{a},b]\text{ (orange)}\\
					   \Pi_{a_j,y}\text{ if }x\in[\bar{a},b],y>b\text{ (grey)} & \text{id}^*_{Y_{i-1}}\text{ if }x,y\in[a,\bar{a})\text{ (brown)} \\
					   \Pi_{x,a_j}\text{ if }x<a,y\in[\bar{a},b]\text{ (turquoise)} & \text{id}^*_{Y_j}\text{ if }x,y\in[\bar{a},b]\text{ (purple)} 
		\end{array} \right.
\end{equation}

\begin{figure}[h]\centering
\includegraphics[width=7cm]{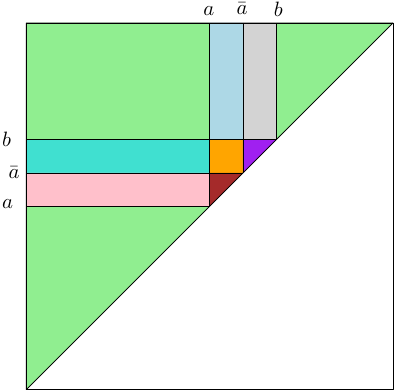}\ \ \includegraphics[width=8cm]{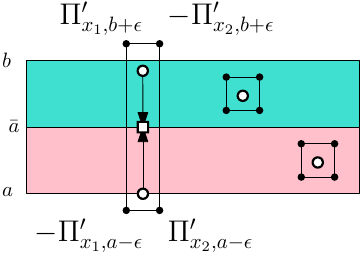}
\caption[Persistence measure for Merge]{
Left: Areas of the extended persistence diagram used in the proof.
Right: Examples of the boxes we use to prove the result (circles represent points before the Merge, squares represent points after the Merge).}
\label{fig:areas}
\end{figure}

The equality between the diagrams follows from these relations and the
inclusion-exclusion formula~(\ref{eq:pers_meas}).  Consider for
instance the case where the point $(x,y)\in\Dg(\pi_2)$ belongs to the
union $A$ of the pink and the turquoise areas.  One can select two
abscissae $x_1<x<x_2$ and an arbitrarily small $\epsilon>0$. Then, the
total multiplicity of the corresponding rectangle $R$ in $\Dg(\pi_2')$
(displayed in the right panel of Figure~\ref{fig:areas}) is given by:

\[ \text{mult}(R)=\text{rank}\ \Pi'_{x_2,a-\epsilon}-\text{rank}\ \Pi'_{x_2,b+\epsilon}+\text{rank}\ \Pi'_{x_1,b+\epsilon}-\text{rank}\ \Pi'_{x_1,a-\epsilon}. \]

The first relation in (\ref{eq:rel}) shows that $R$ has exactly the same multiplicity in $\Dg(\pi_2)$, since all its corners belong to the green area.
As this is true for arbitrarily small $\epsilon>0$, it means that $R'=R\cap A$  also has the same multiplicity in $\Dg(\pi_2)$ as in $\Dg(\pi'_2)$.
Now, if we pick a point inside $R'$ with an ordinate different than $\bar{a}$, we can compute its multiplicity in $\Dg(\pi'_2)$ 
by surrounding it with a box included in the turquoise area (if the ordinate is bigger than $\bar{a}$) or in the pink area (if it is smaller). 
Boxes in the turquoise area have multiplicity 
$\text{rank}\ \Pi'_{x_2,y_1}-\text{rank}\ \Pi'_{x_2,y_2}+\text{rank}\ \Pi'_{x_1,y_2}-\text{rank}\ \Pi'_{x_1,y_1}
=\text{rank}\ \Pi_{x_2,a_j}-\text{rank}\ \Pi_{x_2,a_j}+\text{rank}\ \Pi_{x_1,a_j}-\text{rank}\ \Pi_{x_1,a_j}=0$. 
Similarly, boxes in the pink area also have  multiplicity zero. Thus, all points of $R'$ in $\Dg(\pi'_2)$ have 
ordinate $\bar{a}$. Again, as it is true for $x_1,x_2$ as close to each other as we want, it means that $(x,y)$ 
is snapped to $(x,\bar{a})$ in $\Dg(\pi'_2)$. The treatment of the other areas in the plane is similar.

Now, if $[a,b]$ contains no critical values, then $\Pi'=\Pi$, so the result is clear.
\qed\end{proof}


\paragraph{Split.} Split operators split a critical value $a_i$ into two different ones $a_i-\e$ and $a_i+\e$. 

\begin{definition}[Split]\label{def:split}
Let $T$ be a telescope.
Let $a_i\in\Crit(T)$ and $\epsilon$ such that $$0\leq\epsilon<\displaystyle\min\{a_{i+1}-a_i,a_i-a_{i-1}\}.$$ 
The $\epsilon$-Split on $T$ at $a_i$ is the telescope $T'=\Split_{\epsilon,a_i}(T)$ given by:
\[~\hspace{-10mm}
\begin{array}{c}
...(Y_{i-1}\times[a_{i-1},a_i])\cup_{\psi_{i-1}}(X_i\times\{a_i\})\cup_{\phi_i}(Y_i\times[a_i,a_{i+1}])... \\
\rotatebox[origin=c]{270}{$\mapsto$}\\
...(Y_{i-1}\times[a_{i-1},a_i-\epsilon])\cup_{\psi_{i-1}^{a_i-\epsilon}}(X_i\times\{a_i-\epsilon\})
\cup_{{\rm id}}(X_i\times[a_i-\epsilon,a_i+\epsilon])\cup_{{\rm id}}(X_i\times\{a_i+\epsilon\})
\cup_{\phi_i^{a_i+\epsilon}}(Y_i\times[a_i+\epsilon,a_{i+1}])...
\end{array}
\]
\end{definition}

\begin{figure}[h]\centering
\includegraphics[height=4cm]{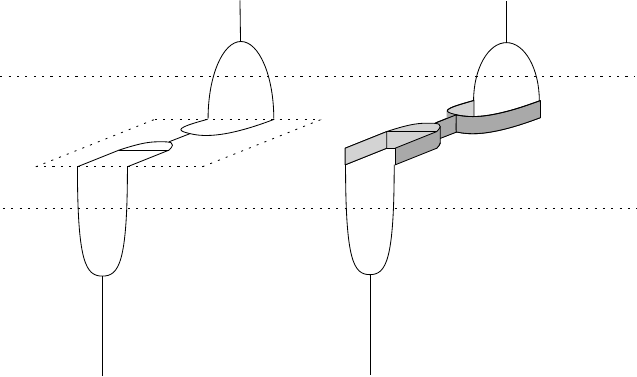}\ \ \ \includegraphics[width = 4cm]{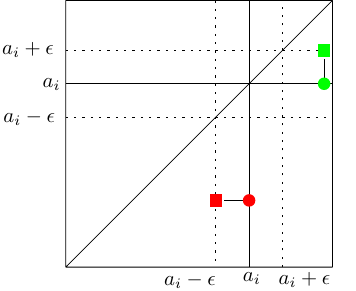}
\caption[Split]{Left: Effect of a Split on a telescope. 
Right: Effect on the corresponding extended persistence diagram. Points before the Split are disks while points after the Split are squares.}
\label{fig:splitspace} 
\end{figure}

See the left panel of Figure~\ref{fig:splitspace} for an
illustration.  

\begin{figure}[h]\centering
\includegraphics[width=7cm]{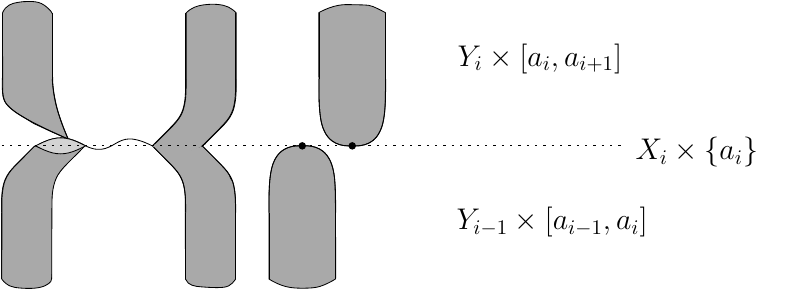}\includegraphics[width=9cm]{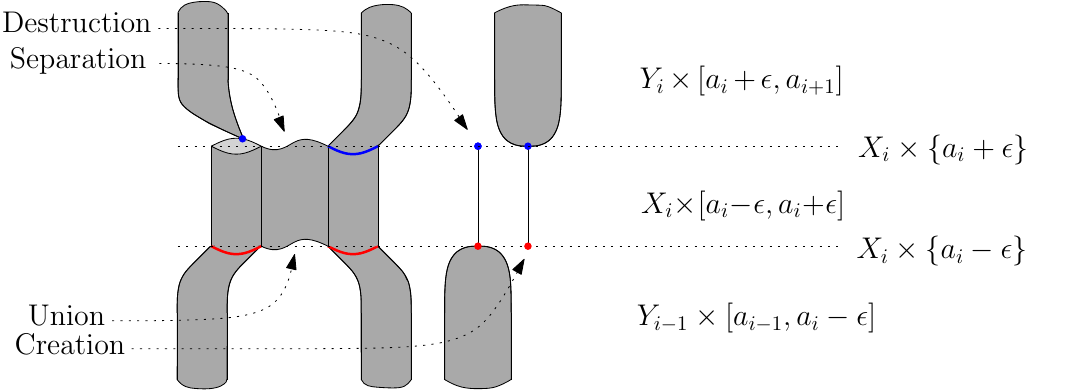}
\caption[Up- and down-forks]{
Left and right panels display the space before and after a Split respectively.
Subsets of $X_i$ that are colored in red and blue correspond to $\im(\pi_1\circ\psi_{i-1})$
and $\im(\pi_1\circ\phi_i)$ respectively.}
\label{fig:split}
\end{figure}

\paragraph{Down- and up-forks.} Splits create particular critical values called {\em down-} and {\em up-forks}.
Intuitively, Split operations allow to distinguish between all possible types of
changes in 0- and 1-dimensional homology of the sublevel and superlevel sets, namely: 
union of two connected components, creation of a connected component, destruction of a connected component,
and separation of a connected component. Unions and creations occur at down-forks while separations and destructions occur at up-forks.
See Figure~\ref{fig:split} for an illustration. We formalize and prove this intuition in Lemma~\ref{lem:dict}.

\begin{definition}
A critical value $a_i\in\Crit(T)$ is called an \emph{up-fork} if $\psi_{i-1}$ is an homeomorphism, and
it is called a \emph{down-fork} if $\phi_i$ is a homeomorphism.
\end{definition}

Since the attaching maps introduced by the Split are
identity maps, we have the following lemma:

\begin{lemma}\label{lem:classif}
The critical values $a_i-\epsilon$ and $a_i+\epsilon$ created with $\Split$ are down- and up-forks respectively.
\end{lemma}

The next lemma is a direct consequence of the existence and continuity
of $\phi_i^{-1}$ (resp. $\psi_{i-1}^{-1}$) when $a_i\in\Crit(T)$ is a
down-fork (resp. up-fork):

\begin{lemma}\label{lem:forks} Let $a_i\in\Crit(T)$. If $a_i$ is an up-fork, then
    $T^{(-\infty,a_i]}$ deform retracts onto
    $T^{(-\infty,\alpha]}$ for all $\alpha \in
  (a_{i-1},a_i]$. 
    If $a_i$ is a down-fork, then
    $T^{[a_i,+\infty)}$ deform retracts onto
      $T^{[\alpha,+\infty)}$ for all $\alpha \in [a_i,a_{i+1})$.
\end{lemma}

Now we can prove the previous intuition concerning down- and up-forks correct:

\begin{lemma}\label{lem:dict}
Let $a_i\in\Crit(T)$.
If $a_i$ is an up-fork, then it can only be the birth time of relative cycles and the death time of relative and extended cycles in $\Dg(\pi_2)$.
If $a_i$ is a down-fork, then it can only be the birth time of ordinary and extended cycles and the death time of ordinary cycles in $\Dg(\pi_2)$.
\end{lemma}

\begin{proof}
Let $0\leq \epsilon,\epsilon' < \min\{a_{i+1}-a_i,a_i-a_{i-1}\}$. Consider the extended persistence module of $\pi_2$:
\[
\begin{array}{lllllllll}
... & \longrightarrow & H_*(T^{(-\infty,a_i-\epsilon]}) & \longrightarrow & H_*(T^{(-\infty,a_i]}) & \longrightarrow & H_*(T^{(-\infty,a_i+\epsilon']}) &\longrightarrow & ... \\
... & \longrightarrow & H_*(T,T^{[a_i+\epsilon',+\infty)}) & \longrightarrow & H_*(T,T^{[a_i,+\infty)}) & \longrightarrow & H_*(T,T^{[a_i-\epsilon,+\infty)}) & \longrightarrow & ...
\end{array}
\]
If $a_i$ is an up-fork, then the composition
$H_*(T^{(-\infty,a_i-\epsilon]}) \to
        H_*(T^{(-\infty,a_i+\epsilon']})$ is an
      isomorphism since $T^{(-\infty,a_i+\epsilon']}$ deform
    retracts onto $T^{(-\infty,a_i-\epsilon]}$ by
      Lemmas~\ref{lem:defret} and \ref{lem:forks}. As
      $\epsilon,\epsilon'$ can be chosen arbitrarily small, there cannot be
      any creation of ordinary or extended cycle at $a_i$. There also
      cannot be any destruction of ordinary cycle.

Similarly, if $a_i$ is a down-fork, then the composition
$H_*(T,T^{[a_i+\epsilon',+\infty)}) \to
  H_*(T,T^{[a_i-\epsilon,+\infty)})$ is an isomorphism
    since $T^{[a_i-\epsilon,+\infty)}$ deform retracts onto
      $T^{[a_i+\epsilon',+\infty)}$. Again, there cannot be any
        destruction of extended or relative cycle at $a_i$. There also
        cannot be any creation of relative cycle.
\qed\end{proof}

\paragraph{Split for persistence diagrams.}
Similarly, we define the $\epsilon$-Split at $a_i$ on a
diagram~$\Dg$ as the diagram $\Split_{\epsilon,a_i}(\Dg)$ given by
$\Split_{\epsilon,a_i}(x,y)= (\bar x, \bar y)$, where:
\[
\bar x = \left\{\begin{array}{l}
x\ \mbox{if}\ x\neq a_i\\
a_i+\epsilon\ \mbox{if}\ 
x=a_i\ \mbox{and}\ 
(x,y)\in\Rel\\
a_i-\epsilon\ \mbox{if}\ 
x=a_i\ \mbox{and}\ 
(x,y)\notin\Rel
\end{array}\right.
\ \mbox{and}\ 
\bar y = \left\{\begin{array}{l}
y\ \mbox{if}\ y\neq a_i\\
a_i-\epsilon\ \mbox{if}\ 
y=a_i\ \mbox{and}\ 
(x,y)\in\Ord\\
a_i+\epsilon\ \mbox{if}\ 
y=a_i\ \mbox{and}\ 
(x,y)\notin\Ord
\end{array}\right.
\]
Points located on the lines $x,y=a_i$ are snapped to the lines
$x,y=a_i\pm\epsilon$ according to their type. Note that the definition
of $\Split_{\epsilon,a_i}(\Dg)$ assumes implicitly that $\Dg$ contains no
point within the horizontal and vertical bands $[a_i-\epsilon, a_i)\times\R$, 
$(a_i, a_i+\epsilon]\times\R$, $\R\times [a_i-\epsilon, a_i)$ and 
$\R\times (a_i, a_i+\epsilon]$, which is the
case under the assumptions of Definition~\ref{def:split}.
See the right panel of Figure~\ref{fig:splitspace} for an
illustration. See also the second intermediate points along the
trajectories of the red points in Figure~\ref{fig:pdTRANS} for another
illustration on extended persistence diagrams. 

\paragraph{Commutativity of the operators.} We now prove 
that extended persistent homology commutes with this operator, i.e.
$\Dg(\Split)=\Split(\Dg)$.

\begin{lemma}
\label{lem:splitprop}
Let $a_i\in\Crit(T)$.
Let $0<\epsilon<\displaystyle\min\{a_{i+1}-a_i,a_i-a_{i-1}\}$, 
$T'=\Split_{\epsilon,a_i}(T)$ and $\pi_2':T'\rightarrow\R$ the projection 
onto the second factor. Then, $\Dg(\pi_2')=\Split_{\epsilon,a_i}(\Dg(\pi_2))$.
\end{lemma}

\begin{proof}
Note that $T=\Merge_{a_i-\epsilon,a_i+\epsilon}(T')$. Hence, by
Lemma~\ref{lem:mergeprop},  $\Dg(\pi_2)$ can be obtained from $\Dg(\pi_2')$ 
with $\Dg(\pi_2)=\Merge_{a_i-\epsilon,a_i+\epsilon}(\Dg(\pi_2'))$.  Note also that
$\pi'_2$ has no critical value within the open interval
$(a_i-\epsilon,a_i+\epsilon)$, so $\Dg(\pi'_2)$ has no point within the
horizontal and vertical bands $\R\times (a_i-\e, a_i+\e)$ and
$(a_i-\e, a_i+\e)\times\R$. Finally, Lemma~\ref{lem:classif} ensures
that $a_i+\epsilon,a_i-\epsilon$ are up- and down-forks respectively,
so Lemma~\ref{lem:dict} tells us exactly where the preimages of the
points of $\Dg(\pi_2)$ through the Merge are located depending on
their type.
\qed\end{proof}


\paragraph{Shift.} Shift operators translate critical values.

\begin{definition}[Shift]\label{def:shift}
Let $T$ be a telescope.
Let $a_i\in\Crit(T)$ and $\epsilon$ such that 
$$0\leq|\epsilon|<\displaystyle\min\{a_{i+1}-a_i,a_i-a_{i-1}\}.$$ 
The $\epsilon$-Shift on $T$ at $a_i$ is the telescope $T'=\Shift_{\epsilon,a_i}(T)$ given by:
\[\begin{array}{c}
...(Y_{i-1}\times[a_{i-1},a_i])\cup_{\psi_{i-1}}(X_i\times\{a_i\})\cup_{\phi_i}(Y_i\times[a_i,a_{i+1}])...\\
\rotatebox[origin=c]{270}{$\mapsto$}\\
...(Y_{i-1}\times[a_{i-1},a_i+\epsilon])\cup_{\psi_{i-1}^{a_i+\epsilon}}(X_i\times\{a_i+\epsilon\})\cup_{\phi_i^{a_i+\epsilon}}(Y_i\times[a_i+\epsilon,a_{i+1}])...
\end{array}\]
\end{definition}

\begin{figure}[h]\centering
\includegraphics[height=4cm]{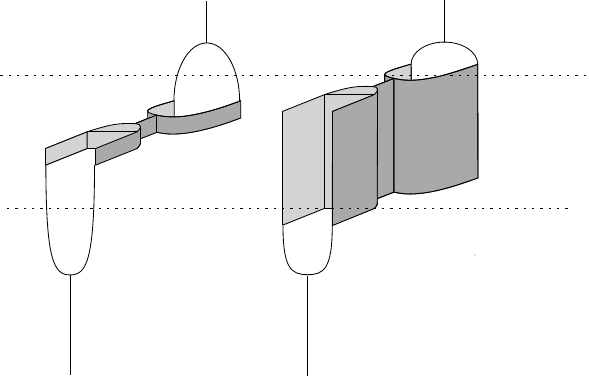}\ \ \ \includegraphics[width = 4cm]{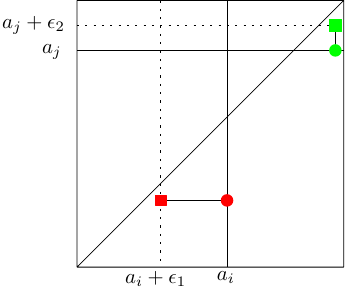}
\caption[Shift]{Left: Effect of a double Shift with
  amplitudes $\epsilon_1<0<\epsilon_2$. 
  Right: Effect on the corresponding
  extended persistence diagram. Points before the Shift are disks while
  points after the Shift are squares.}
\label{fig:shiftspace} 
\end{figure}

See the left panel of Figure~\ref{fig:shiftspace} for an illustration. 

\paragraph{Shift for persistence diagrams.} Similarly, we
define the $\epsilon$-Shift at $a_i$ on a diagram~$\Dg$
as the diagram
$\Shift_{\epsilon,a_i}(\Dg)$ given by $\Shift_{\epsilon,a_i}(x,y)= (\bar x, \bar y)$ where:
\[
\bar x = \left\{\begin{array}{l}
x\ \mbox{if}\ x\neq a_i\\
a_i+\epsilon\ \mbox{otherwise}
\end{array}\right.
\text{ and }
\bar y = \left\{\begin{array}{l}
y\ \mbox{if}\ y\neq a_i\\
a_i+\epsilon\ \mbox{otherwise}
\end{array}\right.
\]
Points located on the lines $x,y=a_i$ are snapped to the lines
$x,y=a_i+\epsilon$.  Note that the definition of
$\Shift_{\epsilon,a_i}(\Dg)$ assumes implicitly that $\Dg$ contains no
point within the horizontal and vertical bands delimited by $a_i$ and
$a_i+\e$, which is the case under the assumptions of
Definition~\ref{def:shift}.  See the right panel of
Figure~\ref{fig:shiftspace} for an illustration. See also the third
intermediate points along the trajectories of the red points in
Figure~\ref{fig:pdTRANS} for another illustration on extended
persistence diagrams. 

\paragraph{Commutativity of the operators.} We now prove 
that extended persistent homology commutes with this operator, i.e.
$\Dg(\Shift)=\Shift(\Dg)$.

\begin{lemma}
\label{lem:shiftprop}
Let $a_i\in\Crit(T)$,  $\epsilon\in \left(a_{i-1}-
a_i,\ a_{i+1}-a_i\right)$, $T'=\Shift_{\epsilon,a_i}(T)$
and $\pi_2':T'\rightarrow\R$ the projection onto the second factor.
Then, $\Dg(\pi_2')=\Shift_{\epsilon,a_i}(\Dg(\pi_2))$.
\end{lemma}

\begin{proof}
Again, the following relations coming from Lemma~\ref{lem:defret}:
\[\hspace{-1cm} \Pi'_{x,y}=\left\{ 
\begin{array}{ll}    
\Pi_{x,y}\text{ if }x,y\notin(a_{i-1},a_{i+1})\text{ (green)} & \Pi_{a_i,y}\text{ if }x\in[a_i+\epsilon,a_{i+1}),y\geq a_{i+1}\text{ (grey)} \\
\Pi_{x,a_{i-1}}\text{ if }x\leq a_{i-1},y\in(a_{i-1},a_i+\epsilon)\text{ (pink)} &  \Pi_{a_{i-1},y}\text{ if }x\in(a_{i-1},a_i+\epsilon), y\geq a_{i+1}\text{ (blue)} \\
\Pi_{x,a_i}\text{ if }x\leq a_{i-1},y\in[a_i+\epsilon,a_{i+1})\text{ (turquoise)} & \text{id}^*_{Y_{i-1}}\text{ if }x,y\in(a_{i-1},a_i+\epsilon)\text{ (brown)} \\
\Pi_{a_{i-1},a_i}\text{ if }x\in(a_{i-1},a_i+\epsilon),y\in[a_i+\epsilon,a_{i+1})\text{ (orange)} & \text{id}^*_{Y_i}\text{ if }x,y\in[a_i+\epsilon,a_{i+1})\text{ (purple)}				   
\end{array} \right. \]
allow us to prove the result similarly to
Lemma~\ref{lem:mergeprop}---see Figure~\ref{fig:areashift}.  For
instance, one can choose a box that intersects the lines
$y=a_i+\epsilon$ and $y=a_i$, show that the total multiplicity is
preserved, then choose another small box that does not intersect
$y=a_i+\epsilon$ inside the first box, and show that its multiplicity
is zero.

\begin{figure}[h]\centering
\includegraphics[width=7cm]{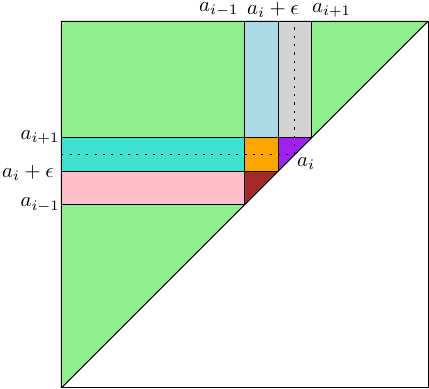}
\caption[Persistence measure for Shift]{
Areas of the extended persistence diagram used in the proof, with $\e <0$.}
\label{fig:areashift}
\end{figure}

\qed\end{proof}

\subsection{Operators on MultiNerve Mapper}
\label{sec:opeMulti}

We first provide invariance results for MultiNerve Mappers computed on telescopes as defined in Section~\ref{sec:telescope}. 
The result is stated in a way that is adapted to its use in the following
sections. The conclusion would still hold under somewhat weaker
assumptions.
\begin{proposition}\label{prop:MNInv}
Let $T$ be a telescope, $\pi_2$ be the projection onto the second coordinate, and $\I$ be a gomic of $\im(\pi_2)$.
Let $\End(\I)$ denote the set of endpoints of intervals of $\I$, sorted in ascending order.
All isomorphisms mentioned in the following items are in the category of combinatorial multigraphs.

\begin{itemize}
\setlength{\itemsep}{0.2pt}
\item[\rm (i)] Let $a\leq b$ such that there exists an interval $I\in\I$ for which $a,b$ belong to 
either $I_\cap^-$, $\tilde I$ or $I_\cap^+$.
Then, $\MMapper_{\pi_2}(\Merge_{a,b}(T),\I)$ is isomorphic to $\MMapper_{\pi_2}(T,\I)$.
 
\item[{\rm (ii)}] Let $a_i\in\Crit(T)\setminus \End(\I)$, and
$a<a_i<b$ with $a,b$ consecutive in $\End(\I)$.  If
$a_{i-1}<a<b<a_{i+1}$ and $0<\e<\min\{a_i-a, b-a_i\}$, then
$\MMapper_{\pi_2}(\Split_{\e,a_i}(T),\I)$ is isomorphic to $\MMapper_{\pi_2}(T,\I)$.

\item[{\rm (iii)}] Let $a_i\in\Crit(T)\setminus \End(\I)$, and $b<a_i<c<d$ with $b,c,d$ consecutive in $\End(\I)$.
If $a_i$ is an up-fork, $(b,c)=I \cap J$ is an intersection, and 
$c-a_i< \e< \min\{d, a_{i+1}\} - a_i$, then $\MMapper_{\pi_2}(\Shift_{\e,a_i}(T),\I)$ is isomorphic to $\MMapper_{\pi_2}(T,\I)$. 

\item[{\rm (iv)}] Let $a_i\in\Crit(T)\setminus \End(\I)$, and $a<b<a_i<c$ with $a,b,c$ consecutive in $\End(\I)$.
If $a_i$ is a down-fork, $(b,c)=I\cap J$ is an intersection, 
and $\max\{a, a_{i-1}\} - a_i<\e<b-a_i$, then
$\MMapper_{\pi_2}(\Shift_{\e,a_i}(T),\I)$ is isomorphic to
$\MMapper_{\pi_2}(T,\I)$.

\end{itemize}
\end{proposition}

\begin{proof}
Under the assumptions given by each item,
the connected components in every intersection $I\cap J$, $I,J\in\I$ and 
in every element $I\in\I$ remain the same after each operation.
Given any intersection $K=I\cap J$, $I,J\in\I$,
or interval $K=I\in\I$, 
we recall that $T^K$ denotes $\pi_1\circ\pi_2^{-1}(K)$. Then, we have:  

\begin{itemize}

\item[(i)]- (ii) $T^K$ deform retracts onto $(\Merge_{a,b}(T))^K$ and $(\Split_{\e,a_i}(T))^K$ deform retracts onto $T^K$; 

\item[(iii)]- (iv) The Shifts move the up-fork to the upper proper subinterval, and the down-fork 
to the lower proper subinterval, which preserves the connected components in each of the
two intervals as well as in their intersection.

\end{itemize}
Thus, the MultiNerve Mapper is not changed by any of the aforementioned operations.   
\qed\end{proof}

\subsection{Connection between the (MultiNerve) Mapper and the Reeb graph}
\label{sec:connection}

In this section, we describe a sequence of metric spaces linking the MultiNerve Mapper and the Reeb graph.
Let $f:X\to\R$ be of Morse type, and let $\I$ be a gomic of $\im(f)$.
Let $T(X,f)$ be the corresponding telescope.
The idea is to move all critical values out 
of the intersection preimages $f^{-1}(I\cap J)$, so that
the MultiNerve Mapper and the Reeb graph become isomorphic.  
For any interval $I\in\I$, we let $a_{\tilde I}<b_{\tilde I}$ be the
endpoints of its proper subinterval~$\tilde I$, so we have $\tilde
I=[a_{\tilde I}, b_{\tilde I}]$. For any non-empty intersection $I\cap J$, 
we fix a subinterval $[a_{I\cap J}, b_{I\cap J}] \subset I\cap J$
such that every critical value within $I\cap J$ falls into $[a_{I\cap J}, b_{I\cap J}]$
(which is possible because $f$ is of Morse type hence has finitely many critical values).  
We then define three different operations individually as follows:

\begin{itemize}
\item 
$\Merge_{\I}$
is the composition of all the
$\Merge_{a_{\tilde I},b_{\tilde I}}$, $I\in\I$, and of all the	
$\Merge_{a_{I\cap J},b_{I\cap J}}$, $I, J \in \I$ and  $ I\cap J\neq \emptyset$.  
All these functions commute, so their composition is well-defined. 
The same holds for the following compositions.
\item $\Split_{\I}$ is the composition of all the 
$\Split_{\e,\bar a}$ with $\bar a$ a critical value after  $\Merge_{\I}$ 
(therefore not an interval endpoint) and
$\e>0$ such that the assumptions of
Proposition~\ref{prop:MNInv}~(ii) are satisfied.
\item $\Shift_{\I}$ is the composition of all the
$\Shift_{\e,\bar a_+}$ with $\bar a_+$ an up-fork critical
value after the $\Split_{\I}$ and $\e>0$ such that the
assumptions of Proposition~\ref{prop:MNInv}~(iii) are satisfied, 
and of all the $\Shift_{\e,\bar a_-}$
with $\bar a_-$ a down-fork critical value after the $\Split_{\I}$ and
$\e<0$ such that the assumptions of
Proposition~\ref{prop:MNInv}~(iv) are satisfied.
After $\Shift_{\I}$ there are no more critical values
located in the intersections of consecutive intervals of $\I$.
\item  $\Merge'_{\I}$ is the composition of all the $\Merge_{a_{\tilde I},b_{\tilde I}}$, $I\in\I$.  
\end{itemize}

We can now define our sequence of intermediate spaces:

\begin{definition}\label{eq:final_transf}
Let $X$ be a topological space, $f:X\rightarrow \R$ be a Morse-type function, and $\I$ be a gomic of $\im(f)$.
Let  $T(X,f)$ be the telescope associated to $f$. We define the telescope $\MMtel$ with:
$$\MMtel(X,f)=\Merge'_{\I}\circ\Shift_{\I}\circ\Split_{\I}\circ\Merge_{\I}(T(X,f)).$$
We also let $\fMMtel$ denote the projection of $\MMtel$ onto the second factor.
\end{definition}

\begin{figure}[h!]\centering
\includegraphics[width=12cm]{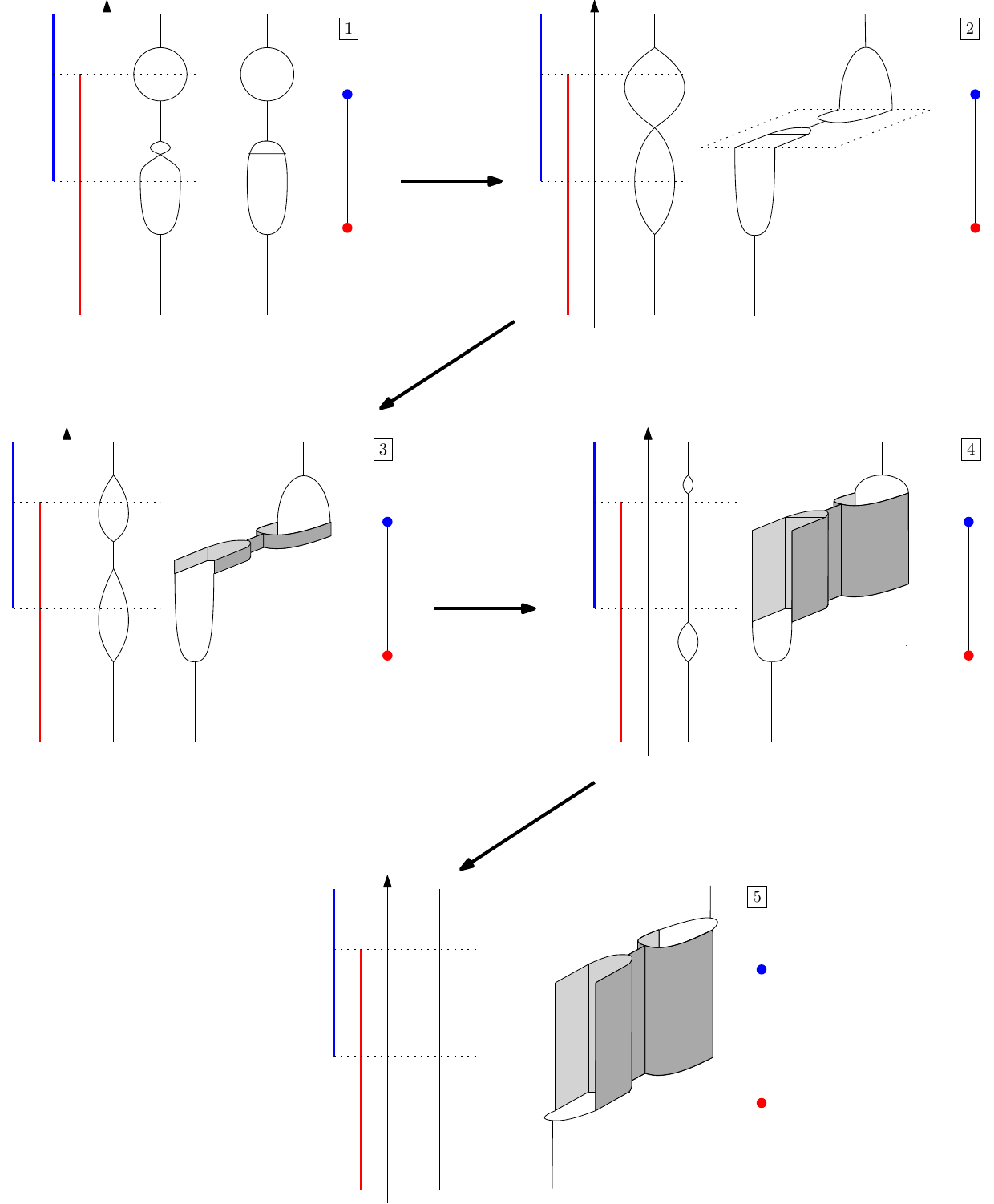}
\caption[Full transformation on spaces]{Illustration of the sequence of
  transformations in~(\ref{eq:final_transf}) on the features located in an interval
  intersection. For each figure, we display the original space
  (middle), its Reeb graph (left) and its MultiNerve Mapper (right).}
\label{fig:summary} 
\end{figure}

See Figure~\ref{fig:summary} for an
illustration of this sequence of transformations.
When often write $\MMtel$ instead of $\MMtel(X,f)$ when the pair $(X,f)$ is clear from the context.
%
In the following, we identify the pair $(T,\pi_2)$
with $(X,f)$ since they are isomorphic in the category of
$\R$-constructible spaces.   
We also let $\fRMMtel:\Reeb_{\fMMtel}(\MMtel)\rightarrow\R$ denote the induced map defined on the Reeb graph of $\MMtel$. 

Thanks to Proposition~\ref{prop:MNInv} and the
choice of the $a_{{\tilde I}},b_{{\tilde I}},a_{I\cap J},b_{I\cap J},\epsilon$ 
in the definitions of $\Merge_{\I}$, $\Split_{\I}, \Shift_{\I}$ and $\Merge'_{\I}$,
we provide Lemma~\ref{lem:inv_Mapper} below, which states that the MultiNerve Mapper
is not affected by this sequence of transformations. 

\begin{lemma}\label{lem:inv_Mapper}
For $(\MMtel,\fMMtel)$ defined as in Definition~\ref{eq:final_transf},
$\MMapper_{\fMMtel}(\MMtel,\I)$ and $\MMapper_f(X,\I)$ are isomorphic as combinatorial multigraphs.
\end{lemma}

This allows us to prove the following result, which states that the MultiNerve Mapper $\MMapper_f(X,\I)$
is actually the same object than the perturbed Reeb graph $\Reeb_{\fMMtel}(\MMtel)$. 

\begin{theorem}\label{th:Dg}
For $(\MMtel,\fMMtel)$ defined as in Definition~\ref{eq:final_transf}, $\MMapper_f(X,\I)$ and
$\mathcal{C}\Reeb_{\fMMtel}(\MMtel)$ are isomorphic as combinatorial multigraphs.
\end{theorem}
We know from Lemma~\ref{lem:inv_Mapper} that
$\MMapper_f(X,\I)$ and $\MMapper_{\fMMtel}(\MMtel,\I)$ are isomorphic as combinatorial multigraphs. Theorem~\ref{th:Dg} is 
then a consequence of the following result, whose hypothesis is satisfied by the~$\MMtel$ of
Definition~\ref{eq:final_transf}:
\begin{lemma}\label{lem:MMapper-Reeb-isom} 
Let $T$ be a telescope and let $\pi_2:T\to\R$ be the projection onto the
second factor.  Suppose that every proper subinterval~$\tilde I$ in the
cover~$\I$ contains exactly one critical value of~$\pi_2$, and that the
intersections $I\cap J$ contain none. Then, $\MMapper_{\pi_2}(T,\I)$ and
$\mathcal{C}\Reeb_{\pi_2}(T)$ are isomorphic as combinatorial multigraphs.
\end{lemma}
\begin{proof}
The nodes of $\mathcal{C}\Reeb_{\pi_2}(T)$ represent the connected components of the
preimages of all critical values of $\pi_2$, while the nodes of
$\MMapper_{\pi_2}(T,\I)$ represent the connected components of the preimages of all
$I\in\I$.  The hypothesis of the lemma implies that there is exactly
one critical value per interval $I\in\I$, hence the nodes of
$\MMapper_{\pi_2}(T,\I)$ and of $\mathcal{C}\Reeb_{\pi_2}(T)$ are in
bijection.  Meanwhile, the edges of $\mathcal{C}\Reeb_{\pi_2}(T)$
are given by the connected components of the $Y_i\times[a_i,a_{i+1}]$. Since the proper subintervals
contain one critical value each and the $I\cap J$ contain none, the
pullbacks of all intersections of consecutive intervals also span the
$Y_i\times[a_i,a_{i+1}]$. Hence, the edges of $\MMapper_{\pi_2}(T,\I)$
are in bijection with the ones of
$\mathcal{C}\Reeb_{\pi_2}(T)$. Moreover, their endpoints are defined
in both cases by the $\phi_i$ and $\psi_i$. Hence the multigraph
isomorphism.
\qed\end{proof}

In passing, it is interesting to study the behavior of the MultiNerve
Mapper as the hypothesis of the lemma is weakened. For instance:
\begin{lemma}\label{lem:MMapper-Reeb-quasi-isom} 
Let $T$ be a telescope and let $\pi_2:T\to\R$ be the projection onto
the second factor.  Suppose that every interval~$I$ in the cover~$\I$
contains at most one critical value of~$\pi_2$. Then, $\MMapper_{\pi_2}(T,\I)$
is obtained from $\mathcal{C}\Reeb_{\pi_2}(T)$ by splitting some
vertices into two and by subdividing some edges once.
\end{lemma}
Thus, the MultiNerve Mapper may non longer be `exactly' isomorphic to
the combinatorial Reeb graph (counter-examples are easy to build, by
making some of the critical values fall into intersections of
intervals in the cover), however it is still isomorphic to it up to
vertex splits and edge subdivisions, which are
topologically trivial modifications.  
\begin{proof}[Lemma~\ref{lem:MMapper-Reeb-quasi-isom}]
The proof is constructive and it proceeds in 3 steps:

\noindent 1. For every interval $I\in\I$ that does not contain a
critical value, add a dummy critical value (with identities as
connecting maps) in the proper subinterval $\tilde I$. The effect on the
Mapper is null, while the effect on the Reeb graph is to subdivide once
each edge crossing the dummy critical value. At this stage, every
interval of $\I$ contains exactly one critical value. For simplicity
we identify $T$ with the new telescope.

\noindent 2. For every interval $I\in\I$ whose corresponding critical
value does not lie in the proper subinterval $\tilde I$ but rather in some
intersection $I\cap J$ (defined uniquely since $\I$ is a gomic), merge
$I$ and $J$ into a single interval $I\cup J$. The coarser cover
$\J$ thus obtained is still a gomic and it has the extra property that
every proper subinterval contains exactly one critical
value and every intersection contains none. Then, by
Lemma~\ref{lem:MMapper-Reeb-isom}, the MultiNerve Mapper
$\MMapper_{\pi_2}(T,\J)$ is isomorphic to the combinatorial Reeb graph
$\mathcal{C}\Reeb_{\pi_2}(T)$.

\noindent 3. There remains to study the differences between
$\MMapper_{\pi_2}(T,\I)$ and $\MMapper_{\pi_2}(T,\J)$. The only
difference between the two covers is that some isolated pairs of
intervals $(I,J)$ have been merged because their intersection $I\cap
J$ contained a critical value $a_i$. For every such pair, there are as
many connected components in the preimage $\pi_2^{-1}(I)$ as in $\pi_2^{-1}(J)$ as in
$\pi_2^{-1}(I\cap J)$ as in $\pi_2^{-1}(I\cup J)$ because $I\cup J$ contains
no critical value other than $a_i$. Hence, every vertex of
$\MMapper_{\pi_2}(T,\J)$ corresponding to a connected component of $\pi_2^{-1}(I\cup J)$
is split into two in $\MMapper_{\pi_2}(T,\I)$. Moreover, the two
copies are connected by a single edge, given by the corresponding connected component
of $\pi_2^{-1}(I\cap J)$. Now, assuming without loss of generality that
$J$ lies above $I$, we have $(I\cup J)_\cap^+=J_\cap^+$, which by
assumption contains no critical value, so the connections between the
vertex copy corresponding to $\pi_2^{-1}(J)$ and the vertices lying above
it in $\MMapper_{\pi_2}(T,\I)$ are the same as the connections between
the original vertex and the vertices lying above it in
$\MMapper_{\pi_2}(T,\J)$. Similarly, $(I\cup J)_\cap^-=I_\cap^-$
contains no critical value by assumption, so the connections between
the vertex copy corresponding to $\pi_2^{-1}(I)$ and the vertices lying
below it in $\MMapper_{\pi_2}(T,\I)$ are the same as the connections
between the original vertex and the vertices lying below it in
$\MMapper_{\pi_2}(T,\J)$.
\qed\end{proof}

\paragraph{Extension to the Mapper.} Due to the simple relation between the Mapper and the MultiNerve Mapper
given by Corollary~\ref{cor:projMNonM}, Theorem~\ref{th:Dg} can be extended for Mappers. 

\begin{definition}\label{cor:finaltransfMapper}
Let $X$ be a topological space, and $f:X\rightarrow\R$ be a Morse-type function. 
Let $(\MMtel,\fMMtel)$ be defined as in Definition~\ref{eq:final_transf}.
Let ${\rm Cyl}(\MMtel)$ be the set of the connected components of the cylinders of $\MMtel$. 
We define the equivalence relation
$\sim$ between elements of ${\rm Cyl}(\MMtel)$ as: $$C\sim C'\Leftrightarrow \left\{\begin{array}{l} C,C'\text{ are connected components of the same cylinder} \\ 
\phi_i(C\times\{a_i\})\text{ and }\phi_i(C'\times\{a_i\})\text{ belong to the same connected component} \\
\psi_i(C\times\{a_{i+1}\})\text{ and }\psi_i(C'\times\{a_{i+1}\})\text{ belong to the same connected component} \end{array}\right.$$
Then, we define $\Mtel$ as $\MMtel/\sim$, equipped with the projection onto 
the second factor that we call $\fMtel$.
\end{definition}

Intuitively, we glue the pairs $C,C'$ of connected components of the same cylinder whose images under the attaching maps are in the same connected component of the critical slice, 
i.e. those that induce edges with the same endpoints in the multinerve. Hence, we obtain the following corollary using Corollary~\ref{cor:projMNonM}:

\begin{corollary}\label{cor:transfMapper}
$\mathcal C\Reeb_{\fMtel}(\Mtel)$ and $\Mapper_f(X,\I)$ are isomorphic as combinatorial multigraphs.   
\end{corollary}


\subsection{Convergence results}
\label{sec:dfdconv}

Recall that the $\distfd$ compares metric graphs, whereas the (MultiNerve) Mappers
are combinatorial graphs.
However, since $\MMapper_f(X,\I)$ and $\Reeb_{\fMMtel}(\MMtel)$ are essentially the same according to Theorem~\ref{th:Dg},
we can use $\Reeb_{\fMMtel}(\MMtel)$ as a metric graph representation of $\MMapper_f(X,\I)$,   
when computing the functional distortion distance.
Note that we could also use $\Reeb_{\fMM}(\MMapper_f(X,\I))$ since it is isomorphic to $\MMapper_f(X,\I)$ as well according to Lemma~\ref{lem:idemMapper},
but its connection to $\Reeb_f(X)$ is unclear. 
On the opposite, even though $\distfd$ is most of the time untractable, its computation is possible with $\Reeb_{\fMMtel}(\MMtel)$
thanks to the sequence of transformations of Definition~\ref{eq:final_transf}.
We will see at the end of the section that $\fMM$ and $\fRMMtel$ actually coincide on $\MMapper_f(X,\I)$.

Theorem~\ref{th:dfd} below shows that $\Reeb_{\fMMtel}(\MMtel)$ is close to $\Reeb_f(X)$ if $\I$ has a small granularity. To prove it,
we use the following lemma, whose proof is just a simple extension of the one of Proposition~3.1 in~\cite{Carriere17a},
and is deferred to Appendix~\ref{sec:proof_dfd}:

\begin{lemma}\label{lem:dfdmerge}
Let $S$ be a set of pairwise disjoint bounded open intervals, and let $\Merge_S$ be defined as
the composition of all $\Merge_{a,b}$, $(a,b)\in S$. 
Let $\Reeb_g$ be a Reeb graph such that $\Crit(\tilde g)\subset \bigcup_{I\in S} I$
and let $\Reeb_{g'}$ be the Reeb graph
of the telescope $\Merge_S(\Reeb_g)$.
Then $\distfd(\Reeb_g,\Reeb_{g'})\leq \sup\{{\rm length}(I)\,:\,I\in S\}$.
\end{lemma}

Given a gomic $\I$, we let $\epsilon_1(\I)=\sup\{{\rm length}(\tilde I)\,:\,I\in\I\}$ and
$\epsilon_2(\I)=\sup\{{\rm length}(I\cap J)\,:\,I,J\in\I\}$. 
Note that $\epsilon_1(\I)$ and $\epsilon_2(\I)$ can be thought of as different types of granularity measures of $\I$.
They are both bounded from above by the granularity of $\I$ as defined in Section~\ref{sec:basicdef}. 

\begin{theorem}\label{th:dfd}
Suppose the granularity of the gomic~$\I$ is at most $\e$.
For $(\MMtel,\fMMtel)$ defined as in Definition~\ref{eq:final_transf}, we have
$\distfd(\Reeb_{\fMMtel}(\MMtel),\Reeb_f(X))\leq \epsilon_1(\I)+\epsilon_2(\I)+\max\{\epsilon_1(\I),\epsilon_2(\I)\}\leq 3\epsilon$.

Moreover, for $(\Mtel,\fMtel)$ defined as in Definition~\ref{cor:finaltransfMapper}, we have
$\distfd(\Reeb_{\fMtel}(\Mtel),\Reeb_f(X))\leq 7\e/2$.
\end{theorem}

\begin{proof}
We start with the MultiNerve Mapper.
By the triangle inequality, it suffices to bound the functional distortion distance for each of the four operations $\Merge_{\I}$,
$\Shift_{\I}$, $\Split_{\I}$ and $\Merge'_{\I}$ individually.
Let $\Reeb_1$ be the Reeb graph of the telescope $\Merge_{\I}(\Reeb_f(X))$.
Similarly, let $\Reeb_2$ be the Reeb graph of $\Split_{\I}(\Reeb_1)$, $\Reeb_3$ be the Reeb graph of $\Shift_{\I}(\Reeb_2)$ 
and $\Reeb_4=\Reeb_{\fMMtel}(\MMtel)$ be the Reeb graph of $\Merge'_{\I}(\Reeb_3)$. Examples of such Reeb graphs can be seen in the left
parts of Figure~\ref{fig:summary}.

Then we have
$\distfd(\Reeb_{\fMMtel}(\MMtel),\Reeb_f(X))\leq \distfd(\Reeb_f(X),\Reeb_1)+\distfd(\Reeb_1,\Reeb_2)+\distfd(\Reeb_2,\Reeb_3)+\distfd(\Reeb_3,\Reeb_4)$.

\begin{itemize}

\item By Lemma~\ref{lem:dfdmerge}, we have $\distfd(\Reeb_f(X),\Reeb_1)\leq \max\{\epsilon_1(\I),\epsilon_2(\I)\}$
and $\distfd(\Reeb_3,\Reeb_4)\leq \epsilon_1(\I)$.

\item Assume without loss of generality that $\Split_{\I}$ is the composition of all $\Split_{\alpha,\bar a}$, 
where $\bar a$ is a critical value of $\Reeb_1$. 
Since $\Reeb_1$ is obtained from $\Reeb_2$ by taking the composition of all $\Merge_{\bar a - \alpha,\bar a + \alpha}$,
it follows from Lemma~\ref{lem:dfdmerge} that $\distfd(\Reeb_1,\Reeb_2)\leq 2\alpha$.

\item Since the assumptions of Prop.~\ref{prop:MNInv}~(iii) and Prop.~\ref{prop:MNInv}~(iv) are satisfied by $\Shift_{\I}$,
it follows that $\Reeb_2$ and $\Reeb_3$ are isomorphic, because the number, the types and the ordering of the critical values of $\Reeb_2$
are preserved when transformed into $\Reeb_3$. It is then straightforward that the functional distortion distance between
$\Reeb_2$ and $\Reeb_3$ is the maximal amplitude of the $\Shift$ operations involved. According to the assumptions 
of Proposition~\ref{prop:MNInv}~(iii) and Proposition~\ref{prop:MNInv}~(iv),
these amplitudes are all bounded by $\epsilon_2(\I)$.

\end{itemize}

The result follows by letting $\alpha \rightarrow 0$.

Concerning the Mapper, the result is obtained by adding an extra $\epsilon/2$ to the previous upper bound,
which corresponds to the functional distortion distance cost of gluing edges with the same endpoints. 
\qed\end{proof}

Note that a similar result can be obtained directly by using the convergence result of the so-called 
{\em interleaving distance}---see Theorem 4.1 in~\cite{Munch16},
and the strong equivalence between the functional distortion distance and this interleaving distance---see Theorem 14 in~\cite{Bauer15}.
However, the upper bound gets larger ($7\epsilon$) and there is no clear intuition on the Reeb graph used to represent the Mapper 
(also called geometric Mapper) in~\cite{Munch16}.

Finally, Theorems~\ref{th:dfd} and~\ref{th:dfdstab} allow us to derive the following result with the triangle inequality:

\begin{corollary} Let $X$ be a topological space, and
let $f,g:X\rightarrow\R$ be two Morse-type functions with continuous sections.
Let $(\MMtel(X,f),\fMMtel)$ (resp. $(\MMtel(X,g),\bar{g}_{\I})$) denote the pair computed with function $f$ (resp. $g$) as in Definition~\ref{eq:final_transf}.
Finally, let $\I$ be a gomic of granularity at most $\epsilon$. Then:
$$\distfd(\Reeb_{\fMMtel}(\MMtel(X,f)),\Reeb_{\bar{g}_{\I}}(\MMtel(X,g)))\leq \|f-g\|_\infty + 6\epsilon.$$
Moreover, for $(\Mtel(X,f),\fMtel)$ and $(\Mtel(X,g),g_{\I})$ computed as in Definition~\ref{cor:finaltransfMapper}, we have:
$$\distfd(\Reeb_{\fMtel}(\Mtel(X,f)),\Reeb_{g_{\I}}(\Mtel(X,g)))\leq \|f-g\|_\infty + 7\epsilon.$$
\end{corollary}

%
%
%

\subsection{An alternative proof of Theorem~\ref{th:ExDgMultiNerve}}
\label{sec:alternateproof}

\begin{figure}[h]\centering
\includegraphics[width=14cm]{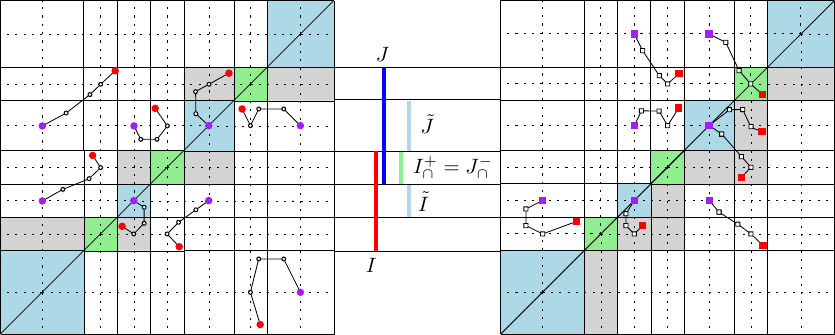}
\caption[Full transformation on persistence diagrams]{
The left panel displays
the trajectories of points in $\Ord$ (disks above the diagonal) and $\Rel$ (disks below the diagonal) while
the right panel displays the trajectories of points in $\Ext$.
For both diagrams, the original point is red, the final point is purple, and intersections and proper subintervals are colored
in green and light blue respectively. 
}
\label{fig:pdTRANS}
\end{figure}

We can prove Theorem~\ref{th:ExDgMultiNerve} again by studying the effect of 
the transformation defined in Definition~\ref{eq:final_transf} on the extended persistence
diagram of~$f$. These effect is illustrated in
Figure~\ref{fig:pdTRANS}.  There are two grids in this figure: the one
with solid lines is defined by the interval endpoints, while the one
with dotted lines is defined by midpoints of proper subintervals and intersections.  In the following, we
use the term \textit{cell} to designate a rectangle of the first grid.
Cells are closed if they correspond to proper subintervals for both
coordinates, they are open if they correspond to intersections for
both coordinates, and they are neither closed nor open otherwise.
Blue and green cells in Figure~\ref{fig:pdTRANS} correspond to squares
associated to a proper subinterval (blue) or intersection (green). We
can now interpret the effects of the transformations
in~(\ref{eq:final_transf}) on the persistence diagram visually:
\begin{itemize}
\item $\Merge_{\I}$ snaps all the points to the second grid by
  Lemma~\ref{lem:mergeprop}.
\item $\Split_{\I}$ moves the points 
to one of the four possible quarters inside their cell, 
depending on the point's type by Lemma~\ref{lem:splitprop}.
More precisely, ordinary points are moved to the down-left quarter,
extended points are moved to the up-left quarter, and relative points are moved to the up-right quarter.
\item Then, $\Shift_{\I}$ moves the points to a neighboring cell by Lemma~\ref{lem:shiftprop}.
This neighboring cell is given by the point's type (as in the case of
$\Split_{\I}$) and by the coordinates of the point.  For instance, an
extended point $(x,y)$ lying in the row of a green cell and in the
column of another green cell, has coordinates that both belong to
interval intersections.
Then, this point is moved to the upper-left neighboring
cell. Differently, an extended point whose abscissa (resp. ordinate)
is in an intersection and whose ordinate (resp. abscissa) is not, is
only moved to the left (resp. upper) cell.  The same can be said for
ordinary (resp. relative) points by changing up-left to down-left
(resp. up-right).  Points whose coordinates both belong to proper
subintervals
are not moved by $\Shift_{\I}$, regardless of their type.
\item Finally, $\Merge'_{\I}$ re-snaps the points to the second grid by Lemma~\ref{lem:mergeprop}.
\end{itemize}
Thus, each point of $\Dg(f)$ can be tracked through the successive
operations of~(\ref{eq:final_transf}), and this tracking leads to the
following elementary observations:
\begin{enumerate}
\item The points of $\Ord(f)$ or $\Rel(f)$ that end their course on
  the diagonal after the sequence of
  operations of Definition~\ref{eq:final_transf} disappear in $\Dg(\fMMtel)$. This is
  because ordinary and relative points cannot be located on the
  diagonal. The rest of the points of $\Ord(f)$ or $\Rel(f)$ are
  preserved with the same type in $\Dg(\fMMtel)$.
\item Differently, all the points of $\Ext(f)$ are preserved with the
  same type ($\Ext$) in $\Dg(\fMMtel)$. However, some of the points of
  $\Ext^-(f)$ may end their course on or across the diagonal after the
  sequence of operations~(\ref{eq:final_transf}), thus switching from
  $\Ext^-(f)$ to $\Ext^+(\fMMtel)$. 
\item All the points lie in the rows and columns of blue cells after
  $\Shift_{\I}$. Therefore, the points that end their course on the diagonal after
  the sequence of operations of Definition~\ref{eq:final_transf} are the ones
  located in blue cells after $\Shift_{\I}$.
\item Since transfers between cells occur only during $\Shift_{\I}$, a
  point $p$ that is not in a blue or green cell initially ends up in a
  blue cell~$B$ after $\Shift_{\I}$ if and only if:
\begin{itemize}
	\item $p$ is extended and it is in the down, right, or
          down-right neighboring cell of $B$ (grey cells in the right
          diagram of Figure~\ref{fig:pdTRANS}), or
	\item $p$ is ordinary and it is in the up neighboring cell of
          $B$ (grey cells above the diagonal in the left diagram of
          Figure~\ref{fig:pdTRANS}), or
	\item $p$ is relative and it is in the down neighboring cell
          of $B$ (grey cells below the diagonal in the left diagram
          of Figure~\ref{fig:pdTRANS}).
\end{itemize}
\item Points that belong to a blue or green cell initially are snapped
  to the diagonal by $\Merge_{\I}$. Among them, those that
  belong to a blue cell stay there until the end, whereas those that
  belong to a green cell eventually leave it---they end up in a blue
  cell after $\Shift_{\I}$ if they are ordinary or relative, while they end
  up in a white cell above the diagonal if they are extended.
\end{enumerate}
The outcome of these observations is the following one. Observations 1, 3,
4, 5 imply that the points of $\Ord(f)$ that disappear in $\Dg(\fMMtel)$
are the ones located initially in the staircase made of
the blue, green and grey areas above the diagonal in the left
panel of Figure~\ref{fig:pdTRANS}, which is nothing but~$Q_O^{\I}$. 
Similarly, the points of $\Rel(f)$
that disappear in $\Dg(\fMMtel)$ are the ones located initially in the
staircase made of the blue, green and grey areas below the
diagonal in the left panel of Figure~\ref{fig:pdTRANS}, which is nothing but~$Q_R^{\I}$. 
Finally, observations 2, 3, 4, 5 imply that the points of $\Ext^-(f)$ that
switch to $\Ext^+(\fMMtel)$ are the ones located initially in the
staircase made of the blue, green and grey areas below the
diagonal in the right panel of Figure~\ref{fig:pdTRANS}, which is nothing but~$Q_{E^-}^{\I}$. 
The rest of the points of $\Dg(f)$ are preserved (albeit shifted) with the same type ($\Ord$,
$\Rel$, $\Ext^+$, $\Ext^-$) in $\Dg(\fMMtel)$. Hence, 
there is a perfect matching between:
\begin{center}
\begin{tabular}{ll}
{\rm (i)} $\Ord(\fMMtel)$ and $\Ord(f)\setminus Q_O^{\I}$ & 
{\rm (iii)} $\Ext^-(\fMMtel)$ and $\Ext^-(f)\setminus Q_{E^-}^{\I}$\\
{\rm (ii)} $\Rel(\fMMtel)$ and $\Rel(f)\setminus Q_R^{\I}$ &
{\rm (iv)} $\Ext^+(\fMMtel)$ and 
$\Ext^+(f) \cup (\Ext^-(f)\cap Q_{E^-}^{\I})$
\end{tabular}
\end{center}

This, combined with Theorem~\ref{thm:pdreeb}, is equivalent to Theorem~\ref{th:ExDgMultiNerve}. 
This matching also shows that the critical points of $\fMMtel$ and $\fMtel$ are  located at the midpoints
of proper subintervals of the gomic's elements. Hence, $\fRMMtel$ and $\fRMtel$ actually coincide with $\fMM$ and $\fM$,
which allows us to state this final result:

\begin{theorem}
\label{th:MultiNerveCover}
Let $X$ be a topological space and $f:X\rightarrow\R$ be a Morse-type function. 
Let $\I$ be a gomic with granularity at most $\epsilon$.
Let $\fMM$ and $\fM$ be as in Definition~\ref{def:arbitfunc}. Then:
\begin{align}
\distb(\Dg(\fMM),\Dg(\tilde f))&\leq \epsilon /2,\label{eq:BottleIneqMMapperReeb} \\
\distb(\Dg(\fM),\Dg(\tilde f))&\leq \epsilon.\label{eq:BottleIneqMapperReeb}
\end{align}
Moreover, in both cases, the matching achieving the distance 
is actually a bijection preserving types. 
\end{theorem}

\section{Discrete Case}
\label{sec:discrete}

In this section we discuss the approximation of the (MultiNerve) Mapper
and of its signature~(\ref{eq:sign1}) when the pair $(X,f)$ is known
only through a finite set of sample points equipped with function
values. Throughout the section, $X$ is a compact Riemannian manifold of $\R^D$, 
$f:X\rightarrow\R$ a Morse-type function, $\I$ a gomic, and $\Xs_n$ a point cloud in~$X$
with $n$ points. 
\subsection{Approximation tools}
\label{sec:tools}

\paragraph{Rips complex.} All constructions take a neighborhood graph as
input, such as for instance the 1-skeleton graph of the Rips complex, defined as follows:
\begin{definition}
Let $\delta \geq 0$ be a scale parameter.
The \emph{Rips complex} of $\Xs_n$ of parameter~$\delta$ is
the simplicial complex $\Rips_\delta(\Xs_n)$ defined by:
\[ 
\{x_{i_0},...,x_{i_k}\}\in\Rips_\delta(\Xs_n)\Leftrightarrow \|x_{i_p}-x_{i_q}\|\leq\delta\text{ for any }0\leq p,q\leq k,. 
\]
Its 1-skeleton graph is called the {\em Rips graph} of parameter~$\delta$ and denoted by~$\Rips_\delta^1(\Xs_n)$.

Moreover, given a geometric realization $|\Rips_\delta(\Xs_n)|$,
we let $\frips{f}:|\Rips_\delta(X_n)|\rightarrow\R$ denote the piecewise-linear
interpolation of $f$ along the simplices of $\Rips_\delta(X_n)$, and
we let $\Reeb_{\frips{f}}(|\Rips_\delta(X_n)|)$ 
denote its Reeb graph, with induced function $\frips{\freeb{f}}$.
\end{definition}

\paragraph{Geometric quantities.} Two geometric quantities that assess the smoothness of topological spaces
will be used in the hypotheses of the results in this section. 
\begin{itemize}

\item {\em Reach.} The {\em medial axis} of $X\subset \R^D$ is the set of points in $\R^D$ with at least two nearest neighbors in $X$:
$${\rm med}(X)=\{y\in\R^D\,:\,{\rm card}\{x\in X\,:\,\|y-x\|=\|y,X\|\} \geq 2\},$$
where $\|y,X\|=\inf\{\|y-x\|\,:\,x\in X\}$.
The {\em reach} of $X$, denoted by ${\rm rch}(X)$, is the distance of $X$ to its medial axis:
$${\rm rch}(X)=\inf\{\|x-m\|\,:\,x\in X,m\in{\rm med}(X)\}.$$

\item {\em Convexity radius.} A set $Y\subseteq X$ is said to be {\em convex} whenever every
geodesic path in $X$ between two points of $Y$ stays in $Y$. The {\em convexity radius} of $X$
is the smallest radius $\rho$ for which every geodesic ball in $X$ of radius less than $\rho$
is convex. 

\end{itemize}

\paragraph{Regularity of the function.}
Intuitively, approximating  a Reeb graph computed with a function $f$ that has large variations is more difficult 
than for a smooth function, for some notion of regularity that we now specify. 
Our result is given in a general setting by considering the {\em modulus of continuity of $f$}. 

\begin{definition}\label{def:modulus}
Let $f:X\rightarrow\R$ be a Morse-type function. 
The {\em modulus of continuity} $\omega_f$ 
of $f$ is: 
\begin{equation*}
\omega_f: \left\{ \begin{array}{lll} \R_+ & \rightarrow & \R_+ \\ \delta & \mapsto & \sup\{|f(x) -  f(x')|\,:\,x,x'\in X,\ \|x - x'\| \leq \delta \} \end{array} \right.
\end{equation*} 
\end{definition}

It follows from the Definition~\ref{def:modulus} that $\omega_f$  satisfies :
\begin{enumerate}
\item $ \omega_f(\delta)  \rightarrow \omega(0) = 0$  as $ \delta \rightarrow 0$ ;
\item $\omega_f$ is nonnegative and non-decreasing on $\R^+$ ;
\item $\omega_f$ is subadditive: $\omega_f(\delta_1 + \delta_2) \leq \omega_f(\delta_1) +\omega_f(\delta_2) $ for any $\delta_1$, $\delta_2 >0$;
\item $\omega_f$  is continous on $\R^+$.
\end{enumerate}

\paragraph{Modulus of continuity.} More generally, we say that a function $\omega$ defined on $\R_+$ is 
\emph{a modulus of continuity} if it satisfies the four properties above, 
and we say that it is \emph{a modulus of continuity for $f$} if, in addition, we have
$\omega\geq\omega_f$.

\subsection{Reeb graph}

The following theorem states that the Rips complex of a point cloud can be used as a proxy for the
 original space $X$. Hence, it is possible to approximate the Reeb graph of $X$,
in the bottleneck distance, by computing the Reeb graph of the Rips complex built on top of the point cloud.

\begin{theorem}[Theorem~4.6 and Remark~2 in~\cite{Dey13a}]\label{th:Reebapprox1}
Assume $X$ has positive reach ${\rm rch}$ and convexity radius $\rho$.
Let $\delta\geq 0$ be a scale parameter,
and let $\omega$ be a modulus of continuity for $f$.

If $4\disth(\Xset,X_n)\leq\delta\leq\min\left\{\frac14 {\rm rch}, \frac14 \rho\right\}$, 
then:
$$\distb(\Ext_1^-(\freeb f),\Ext_1^-(\frips{\freeb f}))\leq 2\omega(\delta).$$
\end{theorem}

Note that the original version of this theorem is only proven for Lipschitz functions in~\cite{Dey13a}, but it extends at no cost,
i.e. with the same proof, to functions with modulus of continuity.

\begin{theorem}[Theorem~2 in~\cite{Chazal09b}]\label{th:Reebapprox0}
Assume $X$ has
convexity radius $\rho >0$. 
Let $\delta>0$ be a scale parameter, 
and let $\omega$ be a modulus of continuity for $f$.

If $4\disth(\Xset,X_n)\leq\delta\leq \rho$, 
then:
$$\max\{\distb(\Ord_0(\freeb f),\Ord_0(\frips{\freeb f})),
\distb(\Ext_0^+(\freeb f),\Ext_0^+(\frips{\freeb f})),
\distb(\Rel_1(\freeb f),\Rel_1(\frips{\freeb f}))\}\leq \omega(\delta).$$
\end{theorem}

Again, the original version of this theorem is only proven for Lipschitz functions 
in~\cite{Chazal09b}, but it extends at no cost
to functions with modulus of continuity. 
Moreover, three more remarks need to be made. 
Firstly, 
this theorem is originally stated only for the ordinary part of the
persistence diagrams but its proof extends to the full extended
filtrations at no extra cost. Secondly, it is stated for a
nested pair of Rips complexes, however, as pointed out
in Section~4.3 in~\cite{Chazal09b}, in 0-dimensional homology a single Rips graph
is sufficient for the theorem to hold. 
Thirdly, its approximation function is piecewise-constant and not piecewise-linear
as in this article.
However, the filtrations induced by the sublevel sets and the superlevel sets of the piecewise-constant function
are actually lower- and upper-star filtrations,
and it is known in that case that piecewise-linear and piecewise-constant functions induce the same persistence diagram.
See Section~2.5 of~\cite{Morozov08} for a proof of this statement.

Combining the two theorems gives the following complete approximation result: 
\begin{theorem}\label{th:Reebapprox}
Under the assumptions of Theorem~\ref{th:Reebapprox1}, we have $\distb(\Dg(\freeb f),\Dg(\frips{\freeb f}))\leq 2\omega(\delta).$
\end{theorem}

\subsection{(MultiNerve) Mapper}
\label{sec:Mapper_constr}

We report three possible constructions for the (MultiNerve) Mapper from
the pair $(X_n, f)$: the first is from the original Mapper
paper~\cite{Singh07}, the second is inspired from the graph-induced complex
paper~\cite{Dey13b}, and the third is simply the Rips complex approximation.  
Given a choice of neighborhood parameter~$\delta$ and the
corresponding Rips graph, the construction
from~\cite{Singh07} uses the vertices as witnesses for the connected components of the
pullback cover on $\Rips^1_\delta(X_n)$ and for their pairwise
intersections. Differently, the construction from~\cite{Dey13b} uses
the edges as witnesses for the pairwise intersections. Thus, both
constructions have the same vertex set but potentially different edge sets. 

\paragraph{Vertex-based connectivity.}

Given an arbitrary interval~$I$ in $\R$, the preimage of~$I$ in~$X_n$ is
defined to be $X_n\cap f^{-1}(I)$, and its connected components are defined to be the connected components of the
induced subgraph $\Rips^1_\delta(X_n\cap f^{-1}(I))$. Then, the
vertices in the (MultiNerve) Mapper are the connected components of the preimages of the
intervals $I\in \I$. Given two intersecting intervals $I,J$ of~$\I$,
given a connected component $C_I$ in the preimage of $I$ and a connected component $C_J$ in the preimage
of $J$, the corresponding vertices are connected by an edge in the Mapper if
there is a connected component in the preimage of $I\cap J$ that is contained in both
$C_I$ and $C_J$; in the MultiNerve Mapper, there are as many copies of
this edge as there are connected components in the preimage of $I\cap J$ that are
contained in $C_I\cap C_J$. We denote these two constructions by 
$\Mapperv_f(\Rips^1_\delta(X_n), \I)$ and $\MMapperv_f(\Rips^1_\delta(X_n), \I)$ respectively.
Moreover, functions $\fM^\bullet$ and $\fMM^\bullet$ can be defined on these constructions
exactly like in Definition~\ref{def:arbitfunc}.

\paragraph{Edge-based connectivity.}
The vertex set of the (MultiNerve) Mapper is the same as in the
previous construction. Now, for any intersecting intervals $I,J$
of~$\I$, we redefine the preimage of the intersection $I\cap J$ to be
the subset of $\Rips^1_\delta(X_n)$ spanned not only by the points of
$X_n\cap f^{-1}(I\cap J)$ and the graph edges connecting them, but also by
the relative interiors of the edges of $\Rips^1_\delta(X_n)$ that have
one vertex in $X_n\cap f^{-1}(I)$ and the other in $X_n\cap f^{-1}(J)$. Then, given a
connected component $C_I$ in the preimage of $I$ and a connected component $C_J$ in the preimage of $J$,
we connect the corresponding vertices by an edge in the Mapper if
there is a connected component of the redefined preimage of $I\cap J$ that
connects\footnote{By which we mean that the closure of the connected component in
  $\Rips^1_\delta(X_n)$ contains points from $C^I$ and from $C^J$.}
$C_I$ and $C_J$ in $\Rips^1_\delta(X_n)$; in the MultiNerve Mapper, we
add as many copies of this edge as there are connected components in the redefined
preimage of $I\cap J$ that connect $C_I$ and $C_J$. We denote these
two constructions by $\Mappere_f(\Rips^1_\delta(X_n), \I)$ and
$\MMappere_f(\Rips^1_\delta(X_n), \I)$ respectively.
Again, we also define functions $\fM^\triangle$ and $\fMM^\triangle$ on these
constructions. 

\begin{figure}\centering
\includegraphics[width=15cm]{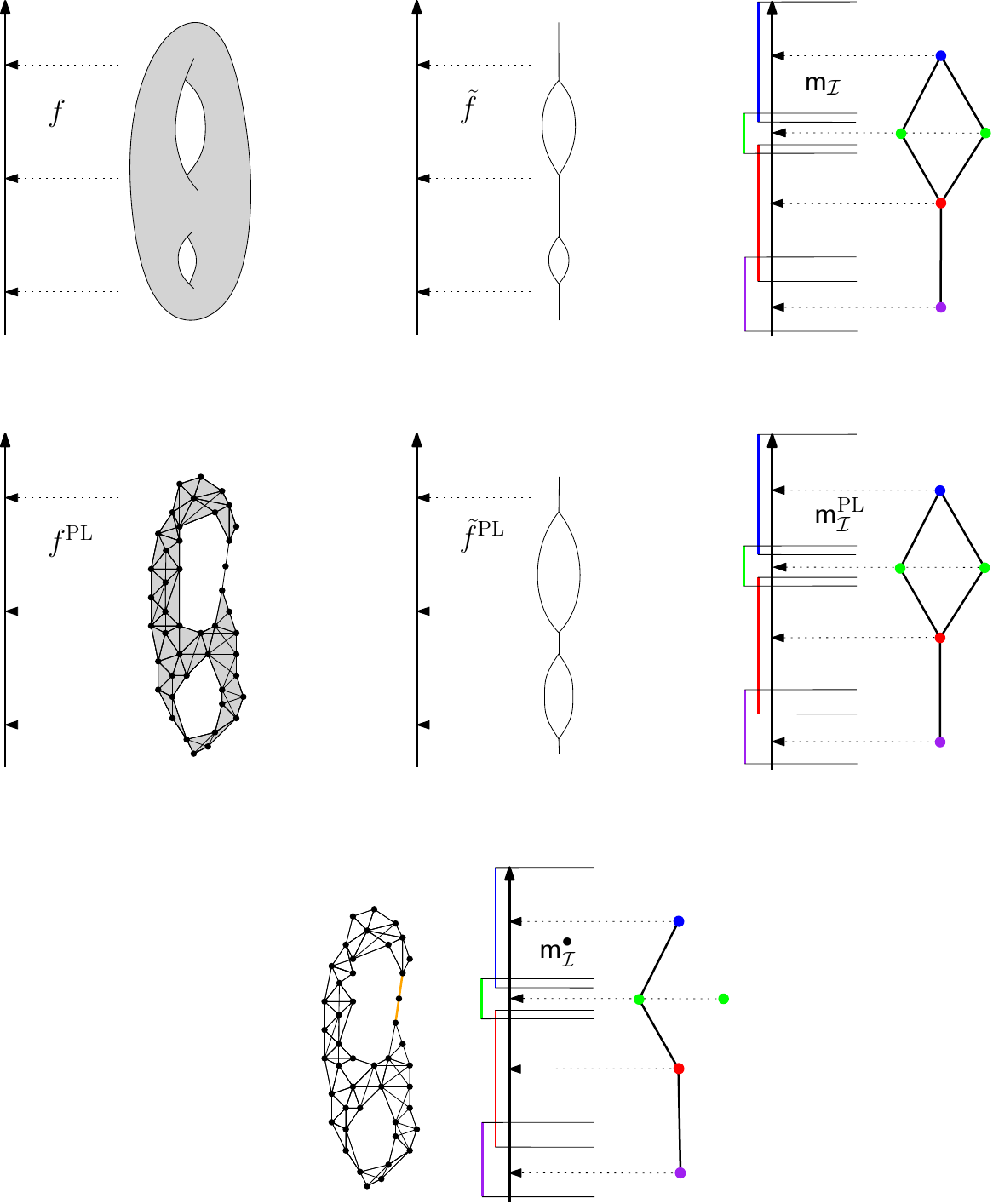}
\caption[Functions on Mapper]{
Examples of the function defined on the original space (left column), its induced function defined on the Reeb graph (middle column) and 
the function defined on the Mapper (right column).
Note that the Mapper computed from the geometric realization of the Rips complex (middle row, right)
is not isomorphic to the standard Mapper (last row), since there are two intersection-crossing edges
in the Rips complex (outlined in orange).}
\label{fig:recap}
\end{figure}

\paragraph{Rips complex.} As for Reeb graphs, one can compute (MultiNerve) Mappers
$\MMapper_{\funcPL}(|\Rips_\delta(X_n)|,\I)$ and $\Mapper_{\funcPL}(|\Rips_\delta(X_n)|,\I)$
from a geometric realization of a Rips complex built on top of $X_n$ with parameter $\delta$,
using the piecewise-linear extension $\funcPL$.  Seeing $|\Rips_\delta(X_n)|$ as a topological space, 
we also define the telescopes $\MMtel(|\Rips_\delta(X_n)|,\funcPL)$ 
and  $\Mtel(|\Rips_\delta(X_n)|,\funcPL)$ 
with Definition~\ref{eq:final_transf} and Definition~\ref{cor:finaltransfMapper},
with corresponding projections onto the second factor $\fPLMMtel$ and $\fPLMtel$.
We recall that $\mathcal{C}\Reeb_{\fPLMMtel}(\MMtel(|\Rips_\delta(X_n)|,\funcPL))$
and $\mathcal{C}\Reeb_{\fPLMtel}(\Mtel(|\Rips_\delta(X_n)|,\funcPL))$ are isomorphic
to $\MMapper_{\funcPL}(|\Rips_\delta(X_n)|,\I)$ and $\Mapper_{\funcPL}(|\Rips_\delta(X_n)|,\I)$
respectively, according to Theorem~\ref{th:Dg} and Corollary~\ref{cor:transfMapper}.
Moreover, the induced maps $\tilde{\bar f}^{\rm PL}_{\I}=\fRPLMMtel$ and $\tilde f^{\rm PL}_{\I}=\fRPLMtel$
are related to $\frips{\freeb{f}}$ according to Theorem~\ref{th:MultiNerveCover}.
See Figure~\ref{fig:recap} for an illustration of these functions. 

\subsection{Relationships between the constructions}
\label{sec:approx-mapper}

In each of the three constructions detailed above, the Mapper is included in the
MultiNerve Mapper by definition. Moreover, the preimages of the
intersections in the second construction are supersets of the preimages in the
first construction, and two different connected components in the same preimage in the first
construction cannot be connected in the second construction, therefore
the (MultiNerve) Mapper from the first construction is included in its
counterpart from the second construction. Hence the following diagram of
inclusions:
\begin{equation}\label{eq:cons_inc}
\xymatrix{
\Mapperv_f(\Rips^1_\delta(X_n), \I)\ar[r]\ar[d] &
\MMapperv_f(\Rips^1_\delta(X_n), \I)\ar[d] \\
\Mappere_f(\Rips^1_\delta(X_n), \I)\ar[r] &
\MMappere_f(\Rips^1_\delta(X_n), \I)
}
\end{equation}
The vertical inclusions become equalities when there are no
{\em intersection-crossing edges} in the Rips graph, defined as
follows:
\begin{definition}\label{def:crossing_edge}
An edge $[u,v]$ of the Rips graph is {\em interval-crossing} if there
is an interval $I\in \I$ such that $I\subseteq (\min\{f(u), f(v)\},\,
\max\{f(u), f(v)\})$. It is {\em intersection-crossing} if there is a
pair of intervals $I,J\in \I$ such that $\emptyset\neq I\cap
J\subseteq (\min\{f(u), f(v)\},\, \max\{f(u), f(v)\})$.
\end{definition}
Indeed, in the absence of intersection-crossing edges,
each connected component in the preimage of an interval intersection in the second
construction contains a vertex and therefore
deform retracts onto the corresponding connected component in the first
construction. Hence:
\begin{lemma}\label{lem:cons_equ}
If there are no intersection-crossing edges in $\Rips^1_\delta(X_n)$, then
$\Mapperv_f(\Rips^1_\delta(X_n), \I)$ is isomorphic to $\Mappere_f(\Rips^1_\delta(X_n),\I)$  
as combinatorial multigraphs. The same is true for 
$\MMapperv_f(\Rips^1_\delta(X_n), \I)$ and $\MMappere_f(\Rips^1_\delta(X_n), \I)$.
\end{lemma}

Concerning the relation between the first two constructions and the third one, we have:
\begin{lemma}\label{lem:cons_equ2}
If there are no interval-crossing edges in $\Rips^1_\delta(X_n)$, then
$\Mapper_{\funcPL}(|\Rips_\delta(X_n)|, \I)$ is isomorphic to 
$\Mappere_f(\Rips^1_\delta(X_n), \I)$ as combinatorial multigraphs.
The same is true for $\MMapper_{\funcPL}(|\Rips_\delta(X_n)|, \I)$ 
and $\MMappere_f(\Rips^1_\delta(X_n), \I)$.
\end{lemma}
\begin{proof}
Note that $\Mapper_{\funcPL}(|\Rips_\delta(X_n)|, \I)$ and
$\MMapper_{\funcPL}(|\Rips_\delta(X_n)|, \I)$ are the same as
$\Mapper_{\funcPL}(|\Rips^1_\delta(X_n)|, \I)$ and 
$\MMapper_{\funcPL}(|\Rips^1_\delta(X_n)|, \I)$ respectively, since only the connected component of the
preimages of intervals are involved in the construction of the
(MultiNerve) Mapper. Hence, for the rest of the proof we set the
domain of $\funcPL$ to be $|\Rips^1_\delta(X_n)|$.  Every
connected component in the preimage through $\funcPL$ of an interval of~$\I$ must
contain a vertex, therefore it deform retracts onto the
corresponding preimage through~$f$. Hence the vertex sets of the
aforementioned simplicial posets are the same.
Every connected component in the preimage through $\funcPL$ of an interval
intersection $I\cap J$ either contains a vertex, in which case it
deform retracts onto the corresponding preimage through~$f$ in
the vertex-based connectivity, or it does not contain any vertex, in
which case the edge of the Rips graph that contains the connected component creates an
edge in the (MultiNerve) Mapper in the edge-based connectivity.
\qed\end{proof}

\begin{figure}[h!]\centering
\includegraphics[width=16cm]{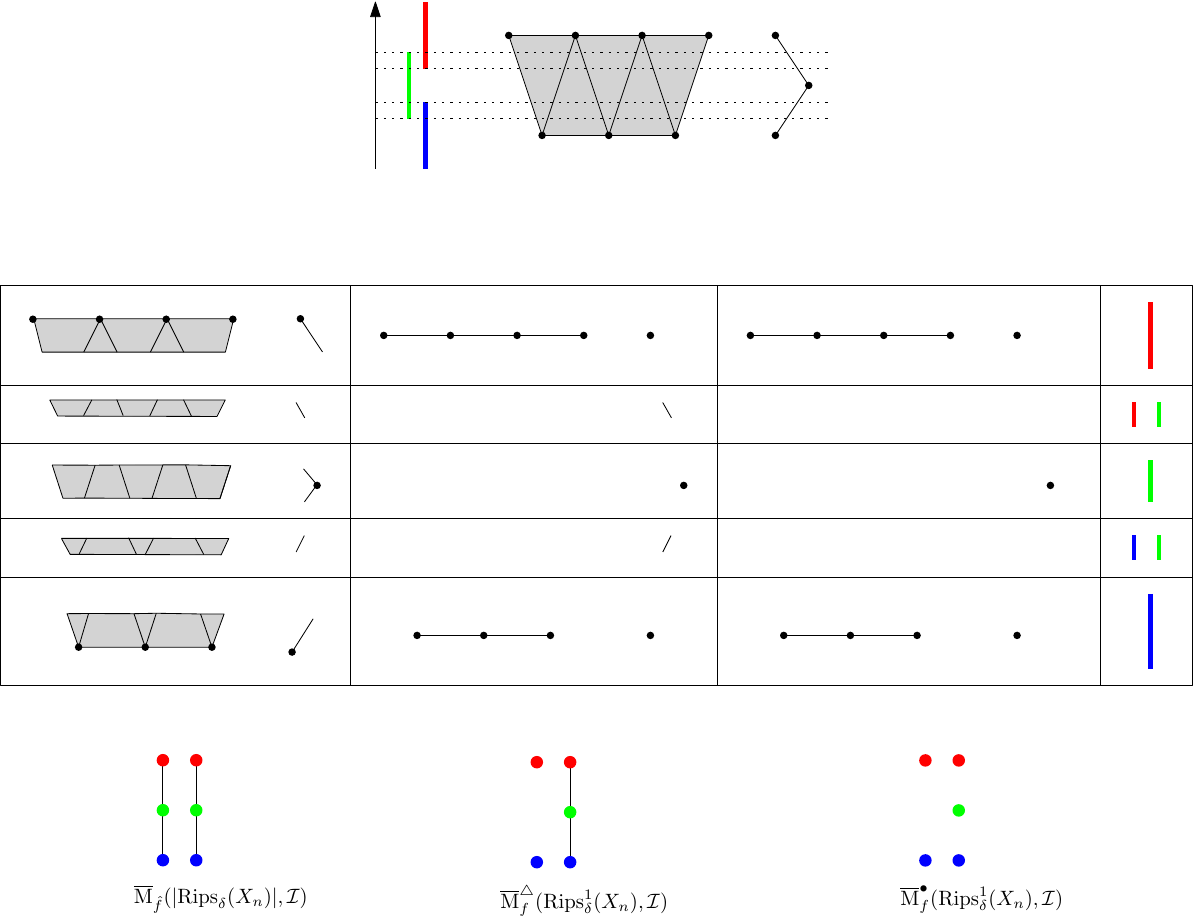}
\caption[Interval- and intersection-crossing edges]{
We study a Rips complex $\Rips_\delta(X_n)$ built on top of a point cloud $X_n$ with ten points. 
This complex has two connected components. We also compute (MultiNerve) Mappers with the height function, whose image is covered by three intervals.
We display the preimages of the intervals and their intersections for  $\MMapper_{\funcPL}(|\Rips_\delta(X_n)|,\I)$,  $\MMappere_f(\Rips_\delta^1(X_n),\I)$ 
and  $\MMapperv_f(\Rips_\delta^1(X_n),\I)$.
The edges of the right connected component are intersection-crossing
but not interval-crossing, so $\MMappere_f(\Rips_\delta^1(X_n),\I)$ recovers it correctly
while $\MMapperv_f(\Rips_\delta^1(X_n),\I)$ fails to.
The edges of the left connected component are interval-crossing, so both $\MMapperv_f(\Rips_\delta^1(X_n),\I)$ and
 $\MMappere_f(\Rips_\delta^1(X_n),\I)$ fail to recover the connected component.}
\label{fig:CounterExampleEdge}
\end{figure}

See Figure~\ref{fig:CounterExampleEdge} for an example showing the importance of the hypothesis in the lemma.

\subsection{Relationships between the signatures}
\label{sec:approx-sigs}

\paragraph{Relationships between the (MultiNerve) Mapper constructions.} 
The following diagram summarizes the relationships between the various
(MultiNerve) Mapper constructions:
\begin{equation}\label{eq:continuous-discrete}
\xymatrix@C=40pt{ (X, f) \ar@{<->}[d] & (|\Rips_\delta(X_n)|), \funcPL) \ar@{<->}[d] \\ 
(\Reeb_f(X), \tilde f) \ar@{<->}[d]
\ar@{<-->}[r] & (\Reeb_{\funcPL}(|\Rips_\delta(X_n)|), \tilde f^{\rm PL}) \ar@{<->}[d]
\\
\MMapper_f(X, \I) \ar@{<->}[d] & \MMapper_{\funcPL}(|\Rips_\delta(X_n)|, \I) \ar@{<->}[d] \ar@{<.>}[r] & 
\MMappere_{f}(\Rips^1_\delta(X_n), \I)\ar@{<.>}[r] \ar@{<->}[d] & 
\MMapperv_{f}(\Rips^1_\delta(X_n), \I) \ar@{<->}[d] \\
\Mapper_f(X, \I) & \Mapper_{\funcPL}(|\Rips_\delta(X_n)|, \I) \ar@{<.>}[r] & 
\Mappere_{f}(\Rips^1_\delta(X_n), \I)\ar@{<.>}[r] & 
\Mapperv_{f}(\Rips^1_\delta(X_n), \I) 
}
\end{equation}
The vertical arrows between the first and second rows are provided by
Theorem~\ref{thm:pdreeb}. The ones between the second, third and fourth rows
are given by Eqs.~(\ref{eq:sign1}) and~(\ref{eq:sign_Mapper}). The
dotted horizontal arrows are provided by Lemmas~\ref{lem:cons_equ}
and~\ref{lem:cons_equ2}. Finally, the dashed horizontal arrow is 
given by Theorem~\ref{th:Reebapprox}. 

\paragraph{Approximation of the (MultiNerve) Mapper.} We then derive from~(\ref{eq:continuous-discrete}) the
following approximation guarantee:
\begin{theorem}\label{th:sig-approx}
Under the assumptions of Theorem~\ref{th:Reebapprox1}, and given a gomic $\I$, we have:
%
\[ 
d_{\I}\left(\Dg(\MMapper_{f}(X, \I)), \, \Dg(\MMapper_{\funcPL}(|\Rips_\delta(X_n)|, \I))\right)\leq 2\omega(\delta).
\]
%
If furthermore there are no interval-crossing edges, then
$\MMapper_{\funcPL}(|\Rips_\delta(X_n)|, \I)$ and
$\MMappere_{f}(\Rips^1_\delta(X_n), \I)$ are isomorphic as combinatorial multigraphs. 

If there are no
intersection-crossing edges either, then $\MMapper_{\funcPL}(|\Rips_\delta(X_n)|, \I)$ and
$\MMapperv_{f}(\Rips^1_\delta(X_n), \I)$ are also isomorphic as combinatorial multigraphs.
\end{theorem}
The same result holds for $\Mapper_{\funcPL}(|\Rips_\delta(X_n)|, \I)$,
$\Mappere_{f}(\Rips^1_\delta(X_n), \I)$ and
$\Mapperv_{f}(\Rips^1_\delta(X_n), \I)$, provided $d_{\I}$ is replaced by
the bottleneck distance with the appropriate extended
staircase~$Q_{E}^{\I}$.  Thus, we can construct discrete (MultiNerve)
Mappers whose signatures approximate the ones of the corresponding
continuous structures $\Mapper_{f}(X, \I)$ and $\MMapper_{f}(X, \I)$.

\paragraph{Approximation of the (MultiNerve) Mapper signature.} In some situations, one is merely interested in
approximating the signatures of $\Mapper_{f}(X, \I)$ and
$\MMapper_{f}(X, \I)$ without actually building corresponding discrete
(MultiNerve) Mappers. In such cases, one can simply apply the scalar
fields analysis approach of~\cite{Chazal09a} to approximate
$\Dg(f)$, then remove the points from $(\Ext_1^+(f)\cup\Ord_1(f))$ as
well as the points lying in their corresponding staircases, to get an
approximation of the signatures:
\begin{theorem}\label{th:sig-approx-general}
Under the assumptions of Theorem~\ref{th:Reebapprox1}, and given a gomic $\I$, let
$\Dg$ denote the extended persistence diagram computed by the
algorithm of~\cite{Chazal09b}, and then pruned by removing the points of
the $\Ext_1^+$ and $\Ord_1$ subdiagrams as well as the points located
in the staircase corresponding to their type. Then this diagram approximates the
signature of $\MMapper_{f}(X, \I)$ as follows: 
\[ 
d_{\I}\left(\Dg(\MMapper_{f}(X, \I)), \, \Dg\right)\leq 2\omega(\delta).
\]
The same bound applies for the approximation of $\Dg(\Mapper_{f}(X,
\I))$, provided the staircase $Q_{E^-}^{\I}$ is replaced by
its extended version~$Q_{E}^{\I}$ in the definitions of~$d_{\I}$ and~$\Dg$.
\end{theorem}
Note that this result holds much more generally than
Theorem~\ref{th:sig-approx}, however there may be no discrete
(MultiNerve) Mapper construction associated with the approximate
diagram~$\Dg$.

\section{Conclusion}

Here we proposed a theoretical framework for the analysis of the
structure and stability of the Mapper.  In particular, we showed that
the topological structure of the (MultiNerve) Mapper is 
nothing but a simplification of the one of the Reeb graph.
From this we derived a bag-of-features type
signature for the (MultiNerve) Mapper, which gives a complete
description of the features present in the graph (trunks, branches,
holes). We related this signature to the persistence of the Reeb
graph, then to the persistence of the pair $(X, f)$ itself, using
various staircases induced by the interval cover~$\I$. This allowed us
to determine the structure of the Mapper from the corresponding
persistence diagrams. We also gave stability guarantees to the
signature, with respect to perturbations of the function, of its
domain, or of the interval cover. These guarantees make it possible to
predict the structural changes that may occur in the Mapper under such
perturbations.  
Then, we refined the analysis by showing that the MultiNerve Mapper is 
actually the Reeb graph of a perturbed pair $(X',f')$ that can be obtained from the
original pair $(X,f)$ through a sequence of controllable elementary perturbations.  
These perturbations allowed us to bound the functional distortion distance between
the (MultiNerve) Mapper and the Reeb graph, and are of an independent interest~\cite{Carriere17c,Carriere17a}.
Finally, we discussed the construction of the Mapper
and the approximation of its signature from point cloud data, relating
the structure of the discrete Mapper to the one of its continuous
version.

We plan to investigate two directions for future work:

\begin{itemize}

\item The first direction deals with the use of
the (MultiNerve) Mapper in Machine Learning.
Is it possible to feed such graphs to Machine Learning methods 
(using e.g. graph kernels) while taking into account the specificities
of these graphs? 


\item The second direction concerns vector-valued functions. Indeed, the
Mapper's definition is not bound to the codomain being $\R$ in any
fundamental way.  However, the stratification of the domain~$X$ by the
critical level sets in the general setting is less well understood.
Chattopadhyay et al.~\cite{Chattopadhyay16} studied the local
configuration of the Reeb space around singular values for bivariate
functions. In a different context, assuming $X$ is a combinatorial
manifold and $f$ is piecewise-linear, Edelsbrunner et
al.~\cite{Edelsbrunner08} showed that the corresponding Reeb
space is a stratified space, and they characterized its coarsest
stratification. We believe part of our framework could be generalized
to these settings, leading to convergence and stability results that
would refine the current state of knowledge~\cite{Munch16}.

\end{itemize}

\bibliography{biblio}
\bibliographystyle{plain}

\setcounter{section}{0}
\renewcommand{\thesection}{\Alph{section}}

\section{Proof of Lemma~\ref{lem:ZZ0}}
\label{sec:proof_pdreeb}


In this proof, we use the notations of Definition~\ref{def:Morse-type}.
Let $0<\epsilon<\frac{1}{2}\min_{k=1,...,n}\ \min\{s_k-a_k,a_k-s_{k-1}\}$. 
The idea of the proof is to replace the right inverse of the
projection $\pi:X\to\Reeb_f(X)$ by a continuous map
$\sigma:\Reeb_f(X)\to X$ such that the composition $\pi\circ\sigma$ is
homotopic to the identity of $\Reeb_f(X)$. In order to make our new
$\sigma$ compatible with the function~$f$, we need to perturb $f$ to
some other function $g$ whose
preimages of intervals $[s_i,s_j]$, $i\leq j$, 
are equal to 
the ones of $f$.

Let $g:X\rightarrow\R$ be defined by:
$$\forall x\in X,\ g(x)=\left\{ \begin{array}{l} f(x)\text{   if }\displaystyle\min_{k=1,...,n} |f(x)-a_k|>2\epsilon \\
						a_i\text{   otherwise, where }i=\text{argmin}_k|f(x)-a_k| \end{array} \right.$$

As $g$ is constant on equivalence classes of $\sim_f$,  there is an induced quotient map $\tilde{g}:\Reeb_f(X)\rightarrow\R$.
Moreover, for any $i\leq j$, we have $g^{-1}([s_i,s_j])=f^{-1}([s_i,s_j])$ by definition of $g$ and $\epsilon$. The
same holds for $\tilde f$ and $\tilde g$.

\begin{figure}[htb]
\begin{center}
\includegraphics[width=13cm]{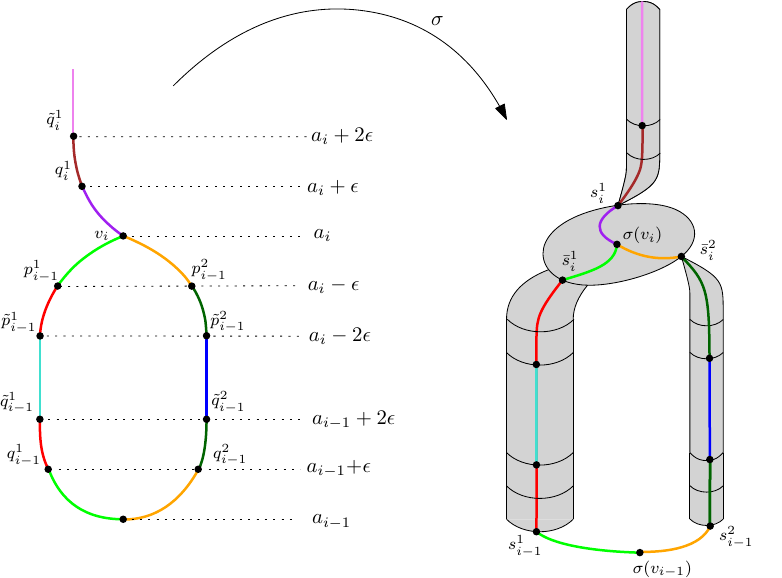}
\caption{\label{fig:contsect}
The left panel displays the Reeb graph and the right panel displays
the space~$X$ itself.
$\sigma$ sends an arc of the Reeb graph to the path with the same color in $X$. }
\end{center}
\end{figure}

Now we want to define a continuous map $\sigma:\Reeb_f(X)\rightarrow X$ 
such that the composition with the projection $\pi\circ\sigma$ is homotopic to $\text{id}_{\Reeb_f(X)}$.
For any node $v_i$, if $Y_{i-1}$ has $k_i$ connected components $Y_{i-1}^1,...,Y_{i-1}^{k_i}$ and 
$Y_i$ has $l_i$ connected components $Y_i^1,...,Y_i^{l_i}$, we let $\{(\tilde{p}_{i-1}^k,p_{i-1}^k)\ |\ k=1,...,k_i\}$ 
and $\{(q_i^l,\tilde{q}_i^l)\ |\ l=1,...,l_i\}$ denote points in $\Reeb_f(X)$ located at levelsets 
$a_i-2\epsilon,a_i-\epsilon,a_i+\epsilon,a_i+2\epsilon$.
See Figure~\ref{fig:contsect}.
For any $i=1,...,n$ and any $l=1,...,l_i$,
%
we select an arbitrary point $y_i^l\in Y_i^l$
%
and we let $s_i^l=\phi_i(y_i^l,a_i)$
%
and $\bar{s}_{i+1}^l=\psi_i(y_i^l,a_{i+1})$. \\

For any critical value $a_i$ and any vertex $v_i$ of $\Reeb_f(X)$ at that
level, we let $\sigma(v_i)$ be an arbitrary point in $\pi^{-1}(v_i)$,
$\sigma(q_i^l)=s_i^l$, and $\sigma(p_{i-1}^k)=\bar{s}_{i}^k$.
Moreover, as there exists a path
$\gamma_k^{i,-}:[a_i-\epsilon,a_i]\rightarrow X$ from $\bar{s}_{i}^k$ to
$\sigma(v_i)$, $\sigma$ sends the arc $[p_{i-1}^k,v_i]$ to this path
$\gamma_{k}^{i,-}$.  Similarly, it sends the arc $[v_i,q_i^l]$ to a path
$\gamma_l^{i,+}:[a_i,a_i+\epsilon]\rightarrow X$ from $\sigma(v_i)$ to
$s_i^l$.  Finally, $\sigma$ also monotonically reparametrizes the
arcs $[\tilde{p}_i^k,p_i^k]$ and $[q_i^l,\tilde{q}_i^l]$.  Let
$\text{param}_i^+:[a_i+\epsilon,a_i+2\epsilon]\rightarrow[a_i,a_i+2\epsilon]$,
and
$\text{param}_i^-:[a_i-2\epsilon,a_i-\epsilon]\rightarrow[a_i-2\epsilon,a_i]$
be these reparametrizations. Again, see Figure~\ref{fig:contsect}. 
More formally, let $x\in X$ and assume that $a_i\leq f(x)\leq a_{i+1}$
and that $\pi(x)$ belongs to the $l$-th edge of the Reeb graph between
these two critical values. Then:

\begin{itemize}
\item $\sigma\circ\pi(x)=\mu_i(y_i^l,f(x))$ if $a_i+2\epsilon\leq f(x)\leq a_{i+1}-2\epsilon$;

\item $\sigma\circ\pi(x)=\mu_i(y_i^l,\text{param}_i^+\circ f(x))$ if $a_i+\epsilon\leq f(x)\leq a_i+2\epsilon$;

\item $\sigma\circ\pi(x)=\mu_i(y_i^l,\text{param}_{i+1}^-\circ f(x))$ if $a_{i+1}-2\epsilon\leq f(x)\leq a_{i+1}-\epsilon$;

\item $\sigma\circ\pi(x)=\gamma_l^{i,+}(f(x))$ if $a_i\leq f(x)\leq a_i+\epsilon$;

\item $\sigma\circ\pi(x)=\gamma_l^{i+1,-}(f(x))$ if $a_{i+1}-\epsilon\leq f(x)\leq a_{i+1}$.
\end{itemize}

By construction we have $g\circ\sigma=\tilde{g}$ and $\tilde{g}\circ\pi=g$ (note that this is not true for $f$). 
%
%

Let $i\leq j$ and $I=[s_i,s_j]$.
Then we have $\pi(g^{-1}(I))\subseteq \tilde g^{-1} (I)$. 
Hence, $\pi$ induces a morphism between $H_0(g^{-1}(I))$ and $H_0(\tilde g^{-1}(I))$. 
Let us show that this morphism is
    an isomorphism. Since $\pi$ is surjective, this boils down to
    showing that $x,y$ are connected in $g^{-1}(I)$ iff
  $\pi(x),\pi(y)$ are connected in $\tilde g^{-1}(I)$.
%
\begin{itemize}
\item If $x,y$ are connected in $g^{-1}(I)$, then so are $\pi(x),\pi(y)$ in $\tilde g^{-1}(I)$, by
  continuity of $\pi$ and the fact that $\tilde{g}\circ\pi=g$.
\item If $\pi(x),\pi(y)$ are connected in $\tilde g^{-1}(I)$, then choose a path $\gamma$
  connecting $\pi(x)$ and $\pi(y)$. Now by definition of $\sigma$,
  there exists a path $\gamma_x$ connecting $x$ and
  $\sigma\circ\pi(x)$ in $g^{-1}(I)$.  Indeed, $\sigma$ can send
  $\pi(x)$ to five different locations in $g^{-1}(I)$
  according to the value of $f(x)$, as seen above.  
  Assume $f(x)\notin \Crit(f)$. Since there is a path $\tilde \gamma$ between
  $x$ and $\mu_i(y_i^l,f(x))$, one can always find a path $\gamma_x$
  between $x$ and $\sigma\circ\pi(x)$ in $g^{-1}(I)$ with an
  appropriate combination of $\tilde \gamma$, $\mu_i(y_i^l,\cdot)$ and
  $\gamma_l^{(i,+)/(i+1,-)}$.  
  Now, assume $f(x)\in\Crit(f)$, and let $v_i=\pi(x)$. Then $\sigma(v_i)$ and $x$
  both belong to $\pi^{-1}(v_i)$, so they belong to the same connected component of 
  the $g^{-1}(g(x))$ and one can find a path between them in $g^{-1}(I)$.
  Similarly, there exists a path $\gamma_y$
  connecting $\sigma\circ\pi(y)$ and $y$ in $g^{-1}(I)$. Then
  $\gamma_y\circ\sigma(\gamma)\circ\gamma_x$ is a path
  between $x$ and $y$ in $g^{-1}(I)$
  by continuity of $\sigma$ and the fact that $g\circ\sigma=\tilde{g}$.
  So $x,y$ are connected in $g^{-1}(I)$.
\end{itemize}
Since $g^{-1}(I)=f^{-1}(I)$ and $\tilde g^{-1}(I)= \tilde f^{-1}(I)$, we have that 
$\pi_*$ is an isomorphism between $H_0(f^{-1}(I))$ and $H_0(\tilde f^{-1}(I))$, and the proof is complete. 

\section{Proof of Lemma~\ref{lem:dfdmerge}}
\label{sec:proof_dfd}

This proof uses the same construction as the proof of Proposition~3.1 in~\cite{Carriere17a}.

\begin{figure}[htb]\centering
\includegraphics[width=10cm]{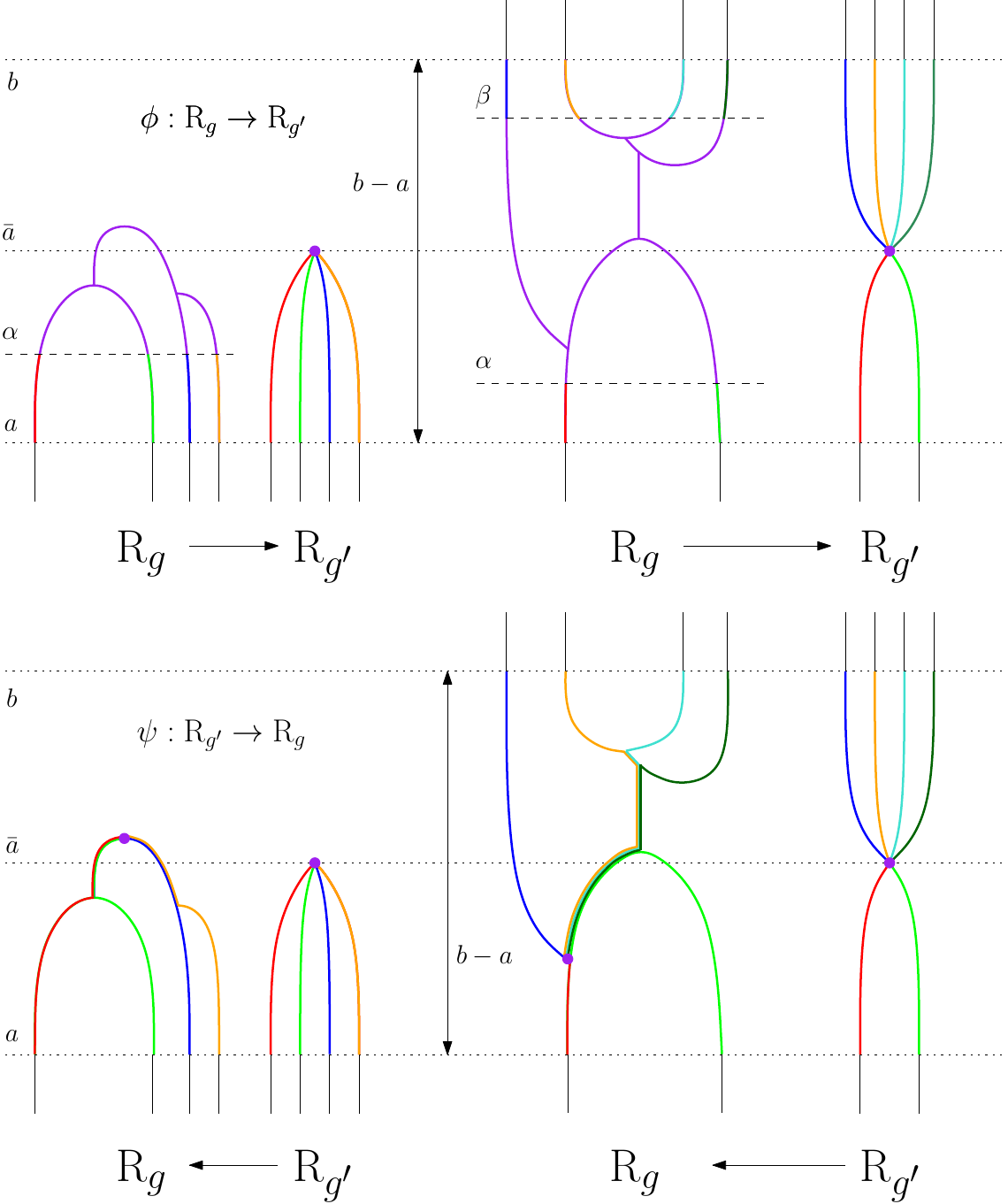}
\caption{\label{fig:fdmerge}
The effects of $\phi$ and $\psi$ around a specific critical value $a_i$ of $f$. Segments are matched according
to their colors (up to reparameterization).
}
\end{figure}

We first prove the result with the stronger assumption that critical values of $\Reeb_g$
are all located inside the interiors of the intervals of $S$, that is $\Crit(\Reeb_g)\subseteq \bigcup_{I\in S} {\rm int}\ I$.
We define explicit continuous maps $\phi:\Reeb_g\rightarrow\Reeb_{g'}$ and
$\psi:\Reeb_{g'}\rightarrow\Reeb_g$ as depicted in
Figure~\ref{fig:fdmerge}.  More precisely, since all critical values of $\Reeb_g$ belong
to $\bigcup_{I\in S} {\rm int}\ I$, we only need to specify $\phi$ and $\psi$
inside each interval $I\in S$ and then ensure that the piecewise-defined maps are assembled consistently.
Let $I=[a,b]$ be such an interval. Since $a,b\not\in \Crit(g)$, there exist two levels 
$a < \alpha\leq\beta < b$ such that $\Reeb_g$ is only composed of arcs in
$[a,\alpha]$ and $[\beta,b]$ (dashed lines in Figure~\ref{fig:fdmerge}).  
For any connected component~$C$ of $g^{-1}([a,b])$, the map $\phi$ sends all points of
$C\cap g^{-1}([\alpha,\beta])$ to the corresponding critical
point $y_C\in\Reeb_{g'}$ created by $\Merge_{a,b}$, and it extends the
arcs of $C\cap g^{-1}([a,\alpha])$ (resp. $C\cap g^{-1}([\beta,b])$) 
into arcs of $(g')^{-1}([a,\bar a])$ (resp. $(g')^{-1}([\bar a,b])$). 
In return, the map $\psi$ sends the critical
point~$y_C$ to an arbitrary point of $C$.  Then, since the Merge
operation preserves connected components, for each arc $A'$ of
$(g')^{-1}([a,b])$ connected to~$y_C$, there is at
least one corresponding path $A$ in $\Reeb$ whose endpoint in
$g^{-1}(a)$ or $g^{-1}(b)$ matches with the one of $A'$
(see the colors in Figure~\ref{fig:fdmerge}).  Hence $\psi$ sends
$A'$ to $A$ and the piecewise-defined maps
are assembled consistently.

\noindent Let us bound the
three terms in the $\max\{\cdots\}$ in~(\ref{eq:dfd})
with this choice of maps $\phi, \psi$: 
\begin{itemize}

\item We first  bound $\|g'-g\circ\psi\|_\infty$. Let $x\in\Reeb_{g'}$. Either $g'(x)\not\in \bigcup_{I\in S} I$,  
and in this case we have $g'(x)=g(\psi(x))$ by definition of~$\psi$;
or, there is $I=[a,b]\in S$ such that $g'(x)\in \mathring{I}$ and then $g(\psi(x))\in [a,b]$.
In both cases $|g'(x)-g\circ\psi(x)| < b-a$. Hence, $\|g'-g\circ\psi\|_\infty < \max_{I\in S}|I|$.

\item Since the previous proof is symmetric in $g$ and $g'$, one also has $\|g-g'\circ\phi\|_\infty < \max_{I\in S}|I|$.

\item We now bound $D(\phi,\psi)$: 

Let $(x,\phi(x)),(\psi(y),y)\in C(\phi,\psi)$ (the cases $(x,\phi(x)),(x',\phi(x'))$
and $(\psi(y),y),(\psi(y'),y')$ are similar). Let $\pi_{g'}:[0,1]\rightarrow\Reeb_{g'}$ be a continuous path from $\phi(x)$ to $y$  
which achieves $d_{g'}(\phi(x),y)$.

\begin{itemize}

\item Assume $g(x)\not\in \bigcup_{I\in S} I$. Then one has $\psi\circ\phi(x)=x$.
Hence, $\pi_g:=\psi\circ\pi_{g'}$ is a path from $x$ to $\psi(y)$. Moreover, since $\|g'-g\circ\psi\|_\infty < b-a$,
it follows that
\begin{equation}\label{eqs:noname}
\begin{array}{l}\max\ \im(g\circ \pi_g) < \max\ \im(g'\circ \pi_{g'}) + \max_{I\in S}|I|,\\[0.5ex]
\min\ \im(g\circ\pi_g) > \min\ \im(g'\circ\pi_{g'}) - \max_{I\in S}|I|.\end{array}
\end{equation}
Hence, one has 
$$\begin{array}{l}d_g(x,\psi(y))\leq \max\ \im(g\circ\pi_g)-\min\ \im(g\circ\pi_g) < d_{g'}(\phi(x),y) + 2\max_{I\in S}|I|,\\
-d_g(x,\psi(y))\geq \min\ \im(g\circ\pi_g)-\max\ \im(g\circ\pi_g) > -d_{g'}(\phi(x),y)- 2\max_{I\in S}|I|.\end{array}$$
This shows that $|d_g(x,\psi(y))-d_{g'}(\phi(x),y)| < 2\max_{I\in S}|I|$.

\item Assume that there is $I=[a,b]\in S$ such that $g(x)\in (a,b)$.
Then, by definition of $\phi,\psi$, we have  $g'(\phi(x))\in(a,b)$, 
and, since $\phi$ and $\psi$ preserve connected components, there is a path
$\pi'_g:[0,1]\rightarrow\Reeb_g$ from $x$ to $\psi\circ\phi(x)$
within the interval $[a,b]$, which itself is
included in the interior of the offset ${\rm im}(g'\circ\pi_{g'})^{b-a}$. Let now
$\pi_g$ be the concatenation of $\pi'_g$ with $\psi\circ\pi_{g'}$, which goes from $x$ to
$\psi(y)$.
Since $\|g'-g\circ\psi\|<\max_{I\in S}|I|$, it follows that ${\rm im}(g\circ\psi\circ\pi_{g'})\subseteq{\rm int}\ {\rm im}(g'\circ\pi_{g'})^{\max_{I\in S}|I|}$, 
and since ${\rm im}(g\circ\pi_g)={\rm im}(g\circ\pi'_g)\cup{\rm im}(g\circ\psi\circ\pi_{g'})$ by
concatenation, one finally has ${\rm im}(g\circ\pi_g)\subseteq {\rm int}\ {\rm im}(g'\circ\pi_{g'})^{\max_{I\in S}|I|}.$ Hence, the 
inequalities of~(\ref{eqs:noname}) hold, implying that $|d_g(x,\psi(y))-d_{g'}(\phi(x),y)| < 2\max_{I\in S}|I|$.
\end{itemize} 
Since these inequalities hold for any couples $(x,\phi(x))$ and $(\psi(y),y)$, we deduce that $D(\phi,\psi) \leq 2\max_{I\in S}|I|$.
\end{itemize}

In the general case where an endpoint of an interval of $S$ can be a critical value of $g$, 
one needs to slightly change the definition of $\phi$. More precisely, $\alpha$ and $\beta$
have to be taken outside but arbitrarily close to $[a,b]$. Then, the previous proof extends verbatim.

\end{document}